\documentstyle[12pt,diagram,amssymb]{article}

\newtheorem{thm}{Theorem}[section]
\newtheorem{assumption}[thm]{Assumption}
\newtheorem{prop}[thm]{Proposition}

\newtheorem{thm-defi}[thm]{Theorem/Definition}
\newtheorem{example}[thm]{Example}
\newtheorem{cor}[thm]{Corollary}

\newtheorem{new-lemma}[thm]{Lemma}
\newtheorem{defi}[thm]{Definition}
\newtheorem{rem}[thm]{Remark}

\newtheorem{condition}[thm]{Condition}

\newcommand{\PP}{{\Bbb P}}

\newcommand{\Integers}{{\Bbb Z}}
\newcommand{\ComplexNumbers}{{\Bbb C}}
\newcommand{\RationalNumbers}{{\Bbb Q}}
\newcommand{\linsys}[1]{{\mid}#1{\mid}}
\newcommand{\Order}[1]{{\mid}#1{\mid}}

\newcommand{\IsomRightArrow}{\stackrel{\cong}{\rightarrow}}
\newcommand{\LongIsomRightArrow}{\stackrel{\cong}{\longrightarrow}}
\newcommand{\RightArrowOf}[1]{\stackrel{#1}{\rightarrow}}

\newcommand{\LongRightArrowOf}[1]{\stackrel{#1}{\longrightarrow}}
\newcommand{\LongIsomRightArrowOf}[1]{
\stackrel
{\stackrel{#1}{\cong}}
{\longrightarrow}
}

\newcommand{\StructureSheaf}[1]{{\cal O}_{#1}}
\newcommand{\EndProof}{\hfill  $\Box$}
\newcommand{\restricted}[2]{#1_{\mid_{#2}}}

\newcommand{\Pic}{\rm Pic}
\newcommand{\Prym}{{\Bbb P}r}
\newcommand{\Sym}{\rm Sym}
\newcommand{\Ext}{\rm Ext}
\newcommand{\Hom}{\rm Hom}

\newcommand{\Aut}{\rm Aut}

\newcommand{\SheafHom}{{\cal H}om}

\newcommand{\Wedge}[1]{\stackrel{#1}{\wedge}}
\newcommand{\Contract}{\rfloor}

\newcommand{\LieAlg}[1]{{\frak #1}}

\newcommand{\A}{{\cal A}}
\newcommand{\C}{{\cal C}}
\newcommand{\D}{{\cal D}}
\newcommand{\JJ}{{\Bbb J}}
\newcommand{\LB}{{\cal L}}
\newcommand{\M}{{\cal M}}
\newcommand{\T}{{\cal T}}
\newcommand{\bba}{{\Bbb A}}

\newcommand{\bbc}{{\Bbb C}}

\newcommand{\bbi}{{\Bbb I}}
\newcommand{\bbj}{{\Bbb J}}
\newcommand{\bbk}{{\Bbb K}}

\newcommand{\bbp}{{\Bbb P}}

\newcommand{\bbs}{{\Bbb S}}

\newcommand{\bbz}{{\Bbb Z}}
\newcommand{\Z}{\bbz}
\newcommand{\J}{\bbj}
\newcommand{\Jchi}{\bbj\otimes_\bbz\chi}
\newcommand{\SS}{\bbs}
\newcommand{\Got}[1]{{\frak #1}}
\newcommand{\Gotg}{{\Got g}}

\newcommand{\gm}{\Got}
\newcommand{\gn}{{\Got n}}
\newcommand{\gh}{{\Got h}}

\newcommand{\gb}{{\Got b}}

\begin{document}
%\rightline{April 6, 98}
\begin{Large}
\centerline{\bf Rank $2$ integrable systems of Prym varieties}
\end{Large}

\vspace{0.3in}

\begin{large}
\centerline{J.C. Hurtubise\footnote{
The first author of this article would like 
to thank NSERC and FCAR for their support}
\ \ \ and  \ \ \ E. Markman\footnote{
The second author would like to thank the Max-Planck Institut f\"{u}r 
Mathematik for its support}}
\end{large}
%*********************************************************************
% Introduction
%*********************************************************************
\section{Introduction}
\label{sec-introduction}

Algebraically integrable systems occur in a wide variety of contexts in 
algebraic geometry. They are algebraic Poisson manifolds, along with a ring of 
Poisson commuting functions, and local identifications of the symplectic leaves
with  open subsets of    Lagrangian fibrations 
\begin{equation}
\label{1.1}
\Pi: \bba\rightarrow U
\end{equation}
where $U$ is some open set in $\bbc^g$, $\bba$ is $2g$ dimensional, 
holomorphically
 symplectic, and the 
fibers of $\Pi$ are Abelian varieties. The identifications are such that the 
Hamiltonian flows are
linear in the Abelian varieties. The local structure of an algebraically 
integrable system,
 at least when the manifold is symplectic, is then that of a fibration 
(\ref{1.1}).
 Given such a local structure, when $U$, say, is a small open ball,
 one could ask 
if there is some suitable invariant that distinguishes it, and allows an 
identification
at least in the neighborhood of a leaf with some known integrable system.

When the Abelian varieties are Jacobians, this question was examined in 
\cite{Hu-algebraic-surfaces}.
Associated to a Lagrangian fibration $\bbj\rightarrow U$ of Jacobians as in 
(\ref{1.1}), one has a family 
$\bbs\rightarrow U$ of genus $g$ Riemann surfaces, and with some choices, one 
has an Abel map $A:\bbs\rightarrow\bbj$. 
A key invariant turns out to be the rank of the 
pull-back $A^*(\Omega)$ of the symplectic form  on $\J$. 
In the minimal case when the rank is two, so that $A^*(\Omega)$
has square zero, there is a null foliation of dimension $g-1$ on $\bbs$ which 
one can quotient 
out to obtain a symplectic surface $Q$. In examples, this surface often 
extends to an algebraic surface
$\tilde Q$, and the symplectic form to a global meromorphic form. One can 
show that the
 original $\bbj$ is  ``birationally'' symplectomorphic 
to the Hilbert scheme $Hilb^g(Q)$ of degree $g$ zero-cycles on $Q$, 
in such a way that the 
Lagrangian leaves of the foliation are the degree $g$ Hilbert schemes 
$Hilb^g(S_u)$ of zero-cycles on the images 
in $Q$ of the fibers $S_u$ of $\bbs$;  birationality here means that
there is an  analytic variety $Y$ of dimension $2g$ which maps in a 
generically bijective way both to $\bbj$ and to $Hilb^g(Q)$.
Thus,  the surface encodes the 
integrable system. There is then a correspondence:
$$
\pmatrix{{\rm local\ rank\!-\!2\  integrable}\cr 
{\rm systems\ of\ Jacobians}}\quad
\leftrightarrow\quad
\pmatrix{\rm symplectic\ surfaces}
$$
Our aim here is to examine the analogous constructions which can be carried 
through
when one has a Lagrangian fibration of generalised Prym varieties over $U$. 
These
Prym varieties are associated to a finite group $W$ acting simultaneously on 
Riemann surfaces and on a lattice $\chi $ of rank $v$. In this case,
the analogue of the rank two condition for the Jacobian case yields us a
$v\!+\!1$-dimensional variety $X$, equipped with a $\chi^*\otimes\bbc$-valued 
two-form. 
At least in the algebro-geometric case,  
we show that giving such an $X$ satisfying appropriate
conditions allows us to reconstruct the integrable system. Our correspondence 
will then be:
$$
\pmatrix{{\rm local\ rank\!-\!2\  integrable}\cr {\rm systems\ of\ Pryms}}
\quad\leftrightarrow\quad
\pmatrix{ 
(v\!+\!1)\!-\!{\rm folds\ X,\ with\ a}
\cr 
(\chi^*\otimes\bbc)\!-\!{\rm valued\ two\ form}
}$$

Both the Jacobian case and this  Prym case of the rank two 
condition have the advantage of 
on the one hand, encompassing
a wide selection, if not most, of the known classical examples, while at 
the same time giving a fairly restrictive constraint on the 
integrable systems  (there are considerable constraints on the variety 
$X$). The condition thus seems to pick out the natural examples,
and provides an intrinsic characterization of them.

 The example  which served as motivation  
is the  class of systems studied by Hitchin 
\cite{hitchin,hitchin-integrable-system}, 
and generalised by Bottacin \cite{botachin} and Markman 
\cite{markman-spectral-curves},  
that has recently provoked a great deal of 
interest (e.g.,  
\cite{donagi-acihs-in-seiberg-witten-theory,
donagi-witten,faltings,kanev,
scognamillo-isogenies}). 
One chooses  a  reductive complex 
linear algebraic group $G$ with Lie algebra $\Gotg$, a compact Riemann 
surface $\Sigma$. Given an effective 
divisor $D$ on $\Sigma$, one considers the moduli space ${\cal M}(D)$ of 
stable pairs $(E,\phi)$ where  $E$ is a holomorphic $G$-bundle over $\Sigma$,
and  $\phi$ is a holomorphic section of $ad(E)(D)$. 
This space has a natural Poisson structure, and with respect to this structure
one has an integrable Hamiltonian system whose Hamiltonians are defined in 
terms of the invariant polynomials on $\Gotg$ applied to the section $\phi$.

This system has two interesting specialisations.
 The first is obtained when $D=0$;
this is Hitchin's original case. The moduli space is then symplectic. 
The other is obtained by choosing the Riemann sphere as the base curve 
$\Sigma$. One can then identify  suitable open sets of the symplectic leaves 
of the moduli space with a reduction of products of  coadjoint orbits for 
the Lie algebra $\Gotg$. Many classical integrable systems can be formulated 
as special cases of these second systems. (\cite{AvM,FRS,RS,AHH1,HH})

 For $G= Gl(r,\bbc)$, the Hitchin systems and their
generalisations are integrable systems of Jacobians and 
give examples of rank two systems \cite{Hu-algebraic-surfaces}. 
To a pair $(E,\phi)$ one can associate
a ``spectral curve'' $S$, the spectrum of $\phi$, in the   
 surface $Q$ obtained from the total space of the canonical line bundle 
$K_\Sigma(D)$ over $\Sigma$ 
by blowing up at some points lying
 over the  divisor $D$.  One
also has a line bundle $L$ over $S$ which is essentially the (dual of) the 
eigenline 
bundle of $S$. These spectral curves $S$, and the line bundles $L$ over them,
 encode the pairs $(E,\phi)$. The Lagrangian fibration is then the projection 
which associates to a pair $(S,L)$ the curve $S$. One shows that 
this system has rank two,     and that
the surface associated to it is locally isomorphic to  $Q$ \cite{Hu-algebraic-surfaces}.

When the group $G$ is a reductive group other than $Gl(r,\bbc)$, 
there is a more  general  spectral curve $S$ which does  not 
lie in $K_\Sigma$, but rather in the vector bundle $K_\Sigma\otimes \gh$, 
where $\gh$ is the Cartan  subalgebra of $\Gotg$, with corresponding group $H$. 
The generalisation of the line bundle in this case becomes an $H$-bundle over 
$S$. 
Both curve and 
bundle are invariant under the action of the Weyl group $W$, and the local 
structure of the integrable system is that of a fibration
$$\Prym \rightarrow U$$
by generalised Prym varieties over a set $U$ parametrising a family of 
$W$-invariant curves. These Prym varieties parametrise $W$-invariant 
$H$-bundles over the 
spectral curves $S$  
(see \cite{donagi-spectral-covers,faltings,scognamillo-isogenies}). 
We will see that these systems satisfy 
our generalised rank two condition, and that the variety $X$ associated 
to them is $K_\Sigma\otimes \gh$.

More generally, then,  we can consider the action of any finite group $W$ on 
the
curves $S_u$ in a family $U$ of $W$-invariant curves, and a 
representation of the group on a finite-dimensional lattice 
$\chi\simeq\bbz^v$,
inducing a representation on the complex vector 
space $V=\chi\otimes_\bbz\bbc$. (In the Hitchin case, $\chi$ was the co-root 
lattice). 
 One then has an action on the Jacobians $J$ of the curves,
and so a diagonal action on $J\otimes_\bbz\chi$. Our generalised Prym 
varieties 
will be connected components of the fixed point sets $(J\otimes_\bbz\chi)^W$.
Now let us consider the associated fibration $\Prym \rightarrow U$,
and assume that the foliation has a symplectic form $\omega$, such that the 
projection
to $U$ is Lagrangian. We will prove:

\begin{thm} 
\label{thm-quotient-of-universal-curve}
Let the system $\Prym \rightarrow U$ have rank two,
and assume that it satisfies genericity conditions A and \^B 
given below.
Restricting U if necessary,

(i) There is a $v+1$ dimensional 
complex manifold $X$  into which the curves $S_u$ all embed.
It is equipped with a generically non-degenerate $V^*$-valued two form $\Omega_V$.
The group $W$ acts
on $X$, preserving $\Omega_V$. 

(ii) Let $v>2$. the manifold $X$ comes equipped with a codimension 1 $W$-invariant foliation
 (``$\phi_0$-foliation''). The form defines a bundle map between the tangent spaces to the 
leaves, and the tensor product of the conormal bundle to the leaves with $V^*$.
\end{thm}

This is proven in section \ref{sec-prym-varieties} of the paper. 
A more detailed version, with a 
stronger genericity 
assumption, is given in Theorem \ref{2.16}. 
In section \ref{sec-rigidity}, we prove  a rigidity theorem: under reasonable 
assumptions, the simultaneous quotient 
of $X$ by the foliation  and by the group $W$ is a fixed curve.

\begin{thm} 
\label{thm-rigidity-in-introduction}
Assume that the system  satisfies genericity conditions
A  and B$'$, as well as regularity condition \ref{cond-regularity},
 given below. 
The manifold $X$ 
admits a $W$-invariant fibration to a closed curve $\Sigma$. 
The quotient curves $S_u/W$ are sections of $X/W \rightarrow \Sigma$. 
In particular, all the quotient curves $S_u/W$ are isomorphic to $\Sigma$. 
%fibers over a closed curve $\Sigma$, in such a way that 
%the  restriction of the  projection to the curves $S_u\subset X$ to $\Sigma$
% is the quotient map under the action of the group $W$: $\Sigma = S_u/W$.
The  projection from $X$ to $\Sigma$ is $V$-Lagrangian 
(Definition \ref{def-V-Lagrangian-fibration})
%and identifies the tangents to the fibers with $T_\Sigma\otimes V^*$. 
which, roughly, means that each fiber 
is $\Omega_V$-isotropic and the fiber over a generic point 
$a \in \Sigma$ has an affine structure modeled after the vector space
$T_{a}\Sigma\otimes V^*$. 
\end{thm}

In section \ref{sec-genericity-conditions} we reformulate 
the genericity and regularity assumptions we have made along the 
way in terms of a single algebro-geometric assumption (Assumption 
\ref{assumption-strong-generation}); 
this assumption is a natural one on the curves, 
and it is equivalent to the assumptions A, B$'$ and condition 3.1 of
Theorem \ref{thm-rigidity-in-introduction}; 
it implies the conditions  A and \^B of Theorem 
\ref{thm-quotient-of-universal-curve}. 
We then develop a simple geometric criterion for proving this assumption
(Theorem \ref{thm-SBPF-holds-if-many-branch-pts-exist}).
In section \ref{sec-from-hilbert-schemes-to-acihs}  
we prove a converse to Theorem \ref{thm-quotient-of-universal-curve}:

\begin{thm}
\label{thm-from-hilbert-schemes-to-acihs-in-introduction}
Let $X$  be a $v+1$ dimensional 
complex manifold, with a submersion onto a closed curve $\Sigma$.
Let $X$ be equipped with a minimally non-degenerate $V^*$-valued two form 
$\Omega_V$, such that the group $W$ acts
on $X$, preserving $\Omega_V$, and preserving the fibers of the map to 
$\Sigma$.
 Assume that there is a smooth $W$-invariant curve $S_0$ in $X$,
on which $W$ acts generically freely with quotient $\Sigma$. 
Then, deforming $S_0$ in $X$,
 the family of smooth $W$-invariant curves $S_u$
defines a rank-2 integrable system of Prym varieties.
\end{thm}

It should be noted that, in the large, several integrable systems
can correspond to the same variety $X$. 
The systems are classified by extra data such as cohomology classes 
(see section \ref{sec-symplectic-structure-on-twisted-pryms}).
Thus, if one starts with an integrable system of Prym 
varieties, applies the first two theorems to get the variety 
$X\rightarrow\Sigma$, and then applies theorem 
\ref{thm-from-hilbert-schemes-to-acihs-in-introduction}, 
the systems one ends up with are not necessarily
globally isomorphic to the original system. The isomorphism 
holds only in a neighbourhood of one of the Prym varieties.

Section \ref{sec-examples} 
is devoted to examples: the classical case of the Prym varieties of 
double covers is considered, 
then coadjoint orbits in loop algebras and the Hitchin systems.
We close with an example of current relevance: the moduli 
space of principal bundles over an elliptic K3 surface, which we discuss in
section \ref{sec-elliptic-K3-systems}. 

%********************************************************************
% Prym varieties
%********************************************************************
\section{Prym varieties}
\label{sec-prym-varieties}

%********************************************************************
% Definitions
%********************************************************************
\subsection{Definitions}

 We begin this section by summarising the  notation we will 
work with. Let
\begin{trivlist}
\item[-] $W$ be a finite group,
\item[-] $\chi$, a free $\Z$-module of rank $v$ on which $W$ acts linearly,
and faithfully,
\item[-] $V=\chi\otimes_\Z\ComplexNumbers$, a complex representation of $W$.
The representation $V$ decomposes as a sum $\sum_im_iV_i$ of  irreducibles 
$V_i$.
\item[-] Let $U\subset \ComplexNumbers^d$  be an open ball,
\item[-] $\sigma:\bbs\rightarrow U$, a holomorphic submersion whose fiber at 
$u\in U$
is a  compact Riemann surface $S_u$ of fixed genus $g$. 
We suppose that $W$ acts on $\bbs\rightarrow U$, inducing a 
trivial action on the base $U$, so that $W$ acts on each $S_u$. We will 
suppose that these
actions are generically free: $W$ embeds into $Aut(S_u)$.
\item[-] Let  $\rho:\J\rightarrow U$ be the corresponding fibration of 
Jacobians, on which $W$ then also acts. $\J$
has an invariant $0$-section, given by  the trivial line bundle. 
Throughout the paper,
we identify Jacobians and the degree 0 component of the Picard variety.
\item [-] For any of the Jacobians $J_u$, we can consider the tensor product
  $J_u\otimes_\Z\chi$, which is isomorphic to the Cartesian product $(J_u)^v$.
Let 
$\rho_\chi: \Jchi\rightarrow U$ denote the associated fibrewise tensor product.
\item [-]  We can associate to the group action a {\it generalised Prym 
variety} 
$$Pr_u = (J_u\otimes_\Z\chi)^W_0,$$
the connected component of the identity of the fixed point set of the diagonal
 $W$-action on $J_u\otimes_\Z\chi$. (The classical
Prym varieties correspond to taking $W=\Z/2$ and $\chi=\Z$ the canonical sign
representation).   Let 
$\Prym\rightarrow U$ be the corresponding fibration; it is a component of 
the fixed point set
of the $W$-action on $\Jchi$.
\end{trivlist}

{\sc Remark:} 1) The condition that $W$ be represented faithfully on $\chi$ is
not unduly restrictive. If the representation has a kernel $G$, one can 
reduce the problem up to isogeny to studying  $W/G$ invariant systems on 
Pryms defined over the curves $S_u/G$. In fact, if $S\rightarrow S/G$ is ramified
then $J({S/G})\rightarrow J({S})$ is injective, with image the identity component 
of the $G$-invariant subgroup $J(S)^G$. Thus,
$$
[J(S)\otimes\chi]^W_0=([J(S)\otimes\chi]^G)^{W/G}_0=[J(S/G)\otimes\chi]_0^{W/G}$$
2) Similarly, if $V$ is faithful but the map from $W$ to $Aut(S)$ has kernel 
$G$,
 one can replace (up to isogeny)
$W$ by $W/G$, and $\chi$ by $\chi^G$, the $G$ invariant submodule.
$$
[J(S)\otimes\chi]^W_0=
([J(S)\otimes\chi]^G)^{W/G}_0=[J(S)\otimes(\chi)^G]_0^{W/G}
$$

%**************************************************************
% Symplectic structures
%**************************************************************
\subsection{Symplectic structures}
\label{sec-symplectic-str}

We will be considering three different types of closed two-forms:
\begin{trivlist}
\item [1.] A $V^*$-valued two-form  
$\Omega_V\in H^0(\J, (\Lambda^2T^*\J)\otimes_\bbc V^*)$
 on the fibration $\J\rightarrow U$.
The fibration is supposed to be isotropic, and the zero section is also 
isotropic.

\item [2.] An ordinary two-form $\Omega$ on the fibration $\Jchi\rightarrow U$.
The fibration and zero-section are again supposed to be isotropic.

\item [3.] The restriction $\omega$ of $\Omega$ to the fibration 
$\Prym\rightarrow U$.
Again, the fibration and zero-section will be isotropic.
\end{trivlist}
We will mostly be interested in the 
case when the form $\omega$ is non-degenerate and the fibration and 
zero-section  are Lagrangian.
We will see that cases (1.) and (2.) are equivalent, and that they are both 
equivalent to (3.) in the invariant case.

We first  consider the action of $W$ on the tangent space 
of $\J$ at a point $p$ in the zero-section. We can decompose $T\J_p$ 
in a $W$-invariant way as follows:
\begin{eqnarray*}
T\J_p &=&  T({\rm fiber}) \oplus T({\rm zero\!-\!section})
\\
T({\rm fiber}) &=& {\cal V}^* \oplus  ({\rm other}).
\end{eqnarray*}
Here ${\cal V}^*$  contains all the $V_i^*$-summands of $T({\rm fiber})$,
say $n_i$ summands for $V_i^*$. The 
``other'' summand represents a $W$-invariant complement to ${\cal V}^*$. 
The space
$T({\rm zero\!-\!section})$ is trivial as a representation of $W$.

The corresponding decomposition of $T(\Jchi)$ along the zero-section
is given by
\begin{eqnarray*}
T(\Jchi)_p &=& T({\rm fiber}\otimes_\Z\chi) \oplus 
T({\rm zero\!-\!section})
\\
T({\rm fiber}\otimes_\Z\chi) &=& 
                       ({\cal V}^*\otimes V) \oplus  (({\rm other})\otimes V).
\end{eqnarray*}
We note that $({\cal V}^*\otimes V)$ contains $\sum_i(n_im_i )$ trivial  
summands. 

Now let $e_1,..,e_v$, $e^1,...,e^v$ be arbitrary dual bases of $\chi$ and 
$\chi^*$ 
respectively. We have the contractions
$$
\matrix {\Jchi&&{\buildrel{\pi_{e^i}}\over
{\longrightarrow}}&&\J\cr
&\searrow&&\swarrow& \cr
&&U&&}
$$
and the tensoring maps
$$
\matrix {\J&&{\buildrel{\eta_{e_i}}\over
{\longrightarrow}}&&\Jchi\cr
&\searrow&&\swarrow& \cr
&&U&&}
$$

\begin{prop}
\label{2.6}
\begin{enumerate}
\item
\label{prop-item-one-to-one-correspondence-of-V-valued-forms}
There is a canonical one-to-one correspondence 
between the  forms $\Omega_V$ and $\Omega$ considered in 1. and 2. above. 
For any choice 
of dual bases, this correspondence is given by:
\begin{equation}
\label{2.7}
\Omega = \sum_{i=1}^v (\pi_{e^i})^* (I( e_i) \Omega_V ),\quad
\Omega_V = \sum_{i=1}^v [\eta_{e_i}^* (\Omega)]\otimes e^i 
\end{equation}
where $I$ is the contraction. The form $\Omega_V$ is invariant with respect 
to the joint $W$-action
on $\J$ and $V^*$ if and only if the form $\Omega$ is invariant
 under the diagonal action on  $\Jchi$.
\item
\label{prop-item-one-to-one-correspondence-of-forms}
There is a canonical one-to-one correspondence given by restriction
between the $W$-invariant forms $\Omega$ of 2. and the forms $\omega$ of  3. 
above.
\end{enumerate}
\end{prop}

{\sc Proof:}
(\ref{prop-item-one-to-one-correspondence-of-V-valued-forms})  
For a fibration $\rho:\J\rightarrow U$ by Abelian varieties, 
one can define 
the vertical tangent bundle $\nu_\J$ over the base $U$ as 
the push-forward $R^0_{\rho^*}(T_{\J/U})$ of the relative tangent bundle;
one checks that the full-back of $\nu_\J$ to $\J$ is 
the tangent bundle $T_{\J/U}$ to the fibres. 
%its lift to $\J$  is indeed  the tangent bundle to the fibres.
Two-forms which are isotropic on both the fibers  
and the zero-section of $\J\rightarrow U$ are then determined by a linear 
bundle map
$TU\rightarrow \nu^*_\J$ that they induce over $U$ \cite {donagi-markman-cime}. As a first application 
of this, we have the lemma:

\begin{new-lemma}
\label{lemma-pull-back-by-linear-comb-is-a-linear-comb}
Consider two fibrations $\bbk\rightarrow U$, 
$\J\rightarrow U$
by abelian varieties, and fiber preserving maps 
$\psi_i:\bbk\rightarrow \J,i=1,..,n$
which are homomorphisms on the fibers.
 Let $\phi$ be a holomorphic
two-form on $\J$, isotropic on the fibers and the zero-section.  
For any choice of integers
$k_i$, the forms 
\begin{equation}
\label{2.9}
\sum_{i=1}^n k_i(\psi_i^*\phi),\quad (\sum_{i=1}^n k_i\psi_i)^*\phi
\end{equation}
are identical.
\end{new-lemma}

To see this, one simply remarks that the maps $TU\rightarrow \nu^*_\bbk$
 that they induce along the zero-section are the same. 

For part \ref{prop-item-one-to-one-correspondence-of-V-valued-forms}
of Proposition \ref{2.6}, 
we first check that the formulae (\ref{2.7}) give the same 
result, independently of the choice of basis.
Let $f_j= \phi_j^i e_i$ be a new basis for $\chi$, with 
$f^j = \tilde\phi^j_i e^i$ the new dual 
basis (we use the
summation convention).
We then have two-forms  $\Omega^e$, $\Omega^f$ defined from $\Omega_V$, with
\begin{eqnarray*}
\Omega^e = (\pi_{e^i})^*(I(e_i)\Omega_V)&=& (\tilde\phi^i_j\pi_{f^j})^*
(I (\phi_i^kf_k)\Omega_V )
\\
&=& 
\tilde\phi^i_j [(\pi_{f^j})^*(I(\phi_i^kf_k)\Omega_V)] \ \ \ 
({\rm using\ the\ lemma})\\
&=& 
\tilde\phi^i_j\phi_i^k[(\pi_{f^j})^*(I(f_k)\Omega_V)]\\
&=& 
\Omega^f
\end{eqnarray*}
Similarly, one checks that the  forms  $\Omega^e_V$, $\Omega^f_V$ defined from
 $\Omega$ are the same.

One next checks that the two formulae of (\ref{2.7}) are inverses of each 
other,
for example:
$$
\eta_{e_i}^*[(\pi_{e^k})^*(I(e_k)\Omega_V)]\otimes e^i= 
(\pi_{e^k}\circ \eta_{e_i})^*
(I(e_k)\Omega_V)\otimes e^i= \Omega_V,
$$
as $(\pi_{e^k}\circ \eta_{e_i})^* = \delta^k_i \cdot \bbi$.

For part \ref{prop-item-one-to-one-correspondence-of-forms}, 
we note that since the map $TU\rightarrow \nu^*_{\Jchi}$ determining the form
$\Omega$, considered as 
a map along the zero section, is equivariant, and the action of $W$ on $TU$ is trivial,
the image of the map must lie in the trivial summand of $\nu^*_{\Jchi}$, that is the vertical
cotangent bundle $\nu^*_{\Prym}$ of the Prym fibration. Nothing is lost, therefore, under restriction from
$\Jchi$ to $\Prym$: the restriction map taking $\Omega$ to $\omega$ is injective . To see that it is
surjective, we note that there is a natural averaging map $Av: \Jchi
\rightarrow \Prym$ given by $x \mapsto \sum_{g\in W}gx$. The restriction of $Av$
to $\Prym$ is simply multiplication by $|W|$, and so setting $\Omega =
Av^*(\omega)/|W|$ gives the correct inverse. This completes the proof of
Proposition 
\ref{2.6}. 
\EndProof

{\sc Remark:} As we have seen,  the forms $\Omega_V$, $\Omega$
are determined by bundle maps:
\begin{eqnarray*}
\hat\Omega_V:  TU & \rightarrow & \nu^*_\J\otimes V^* 
\\
\hat\Omega:  TU & \rightarrow & \nu^*_{\Jchi} 
\end{eqnarray*}
Under the identification 
$\nu^*_\J\otimes V \simeq \nu^*_{\Jchi}$, one has $\hat\Omega_V = 
\hat\Omega$. 

We are  interested in the case when $\Prym\rightarrow U $ is a Lagrangian 
fibration, that is, an integrable system.
 We therefore
want  $\omega$ to be nondegenerate, so that the map 
$TU\rightarrow \nu^*_{\Prym}$ is 
a bundle isomorphism. Correspondingly, $\Omega$ will be non-degenerate
when the map $TU\rightarrow \nu^*_\Jchi$ is an isomorphism on the trivial 
summand.
For the form $\Omega_V$, the corresponding condition is more elaborate; 
when the representation $V$ is irreducible, the condition is that  $\Omega_V$ 
defines  an isomorphism $TU\otimes V\rightarrow (\nu^*_\J)_V$ of the
 $V$-summands.

%*********************************************************
% Compatibility of 2 form with group structure 
%*********************************************************
\subsubsection{Compatibility of a 2-form with the  group structure}
\label{sec-compatibility-of-2-form-with-group-structure}
Let $h : \A^0 \rightarrow U$ be a family of abelian varieties
over a smooth base $U$ with a zero section. Let 
$\Gamma \rightarrow U$ be a discrete group scheme over $U$
and 
\[
0 \rightarrow \A^0 \rightarrow \A \LongRightArrowOf{c} \Gamma \rightarrow 0
\]
an extension. 
Given a $2$-form  $\omega$ on $\A$ and a section
$\gamma$ of $\Gamma$, denote by $\omega^{\gamma}$ the restriction of 
$\omega$ to $A^\gamma$. 

\begin{defi}
{\rm
\label{def-compatability-of-2-form-with-group-structure}
A $2$-form $\omega$ on $\A$ is {\em compatible with the group structure}
if it satisfies the following condition:
Given local (analytic or etale)
sections $\beta$ and $\gamma$ of $\Gamma$ and  a section $\xi$ of $\A$ 
with $c(\xi) = \gamma$, the translation
\[
t_{\xi} : \A^\beta \rightarrow \A^{\beta+\gamma}
\]
satisfies 
\begin{equation}
\label{eq-compatibility-with-group-structure}
(t_\xi)^*(\omega^{\beta+\gamma}) = 
\omega^\beta - h^*(\xi^*(\omega^\gamma)). 
\end{equation}
}
\end{defi}

\noindent
In particular, the sub-sheaf $\A_\omega$ of isotropic sections of $(\A,\omega)$
is a sheaf of sub-groups. 
It acts on $\A$ by translations leaving $\omega$ invariant. 
It also follows that we have the equality of 
$2$-forms on $U$
\[
\xi^*(\omega) + (-\xi)^*(\omega) = 0
\]
for any local section $\xi:U\rightarrow \A$. 
Conversely, given 1) a $2$-form $\omega^0$ on $\A^0$ 
with respect to which the zero section is lagrangian and 2) a sub-group-sheaf 
$\A_\omega$ of $\A$ extending  $\A^0_\omega$ by $\Gamma$, then there exists a
unique $2$-form $\omega$ on $\A$ which is 
compatible with the group structure and restricts to $\omega^0$. 

\begin{example}
\label{example-acihs-of-jacobians-with-symplectic-str-compatible-with-group}
{\rm 
We construct  a compatible  extension of a symplectic structure from  
$\A^0$ to $\A$ in the two equivalent cases: 
\begin{trivlist}
\item[1)] 
$\A = \JJ^{\bullet}$ is the relative Jacobian of all degrees, $\Gamma$ 
the trivial group scheme $\Integers_{U}$, and the $2$-form $\Omega_V$ 
is $V^*$-valued, and 
\item [2)] 
$\A = \JJ^{\bullet}\otimes_{\Integers} \chi$,  $\Gamma$ 
the trivial group scheme $\chi_{U}$, and  $\Omega$ is an ordinary $2$-form. 
\end{trivlist}
%The extenstion of $\Omega^0_V$ from $\JJ^0$ to $\JJ^{\bullet}$ 
%amounts to 
Let $\D$ be a $W$-invariant section of  $\JJ^d_{\SS/U}\rightarrow U$ 
of positive degree $d$. For example, $\D$ could be the relative canonical
line-bundle if the genus of $S_u$ is at least $2$. 
We get a subsheaf $\LB_\D$ of $\JJ^\bullet_{\SS/U}$ 
of sections of which a sufficiently large 
tensor power is a multiple of $\D$.
\begin{equation}
\label{eq-extension-of-torsion-sheaf-to-non-zero-degrees}
\LB_\D := 
\{
s \mid \ \exists n,m, \ \ \ n\neq 0, \ \mbox{such that} \ s^{\otimes n} = 
\D^{\otimes m} 
\}.
\end{equation}
The subsheaf $\LB_\D$ is an extension of the subsheaf $\LB^0$ 
of torsion sections of
$\JJ^0_{\SS/U}$ by the trivial sheaf $\Integers_U$ (or, more canonically,
the local system over $U$ of 2-nd integral cohomology of 
fibers of $\SS\rightarrow U$).
Lemma \ref{lemma-pull-back-by-linear-comb-is-a-linear-comb} implies that 
the pull-back of $\Omega_V$ by the homomorphism
$[k]:\JJ^0_{\SS/U}\rightarrow \JJ^0_{\SS/U}$ (raising to the $k$-th 
tensor power) is equal to $k\cdot \Omega_V$. Hence, the subsheaf 
$\LB^0$ is $\Omega_V$-isotropic. It follows that there is a unique 
$W$-invariant extension of $\Omega_V$ to $\JJ^\bullet_{\SS/U}$ 
with respect to which 
$\LB_\D$ is isotropic and such that  the extended $\Omega_V$
is invariant under translation by a local section of $\LB_\D$. 
By construction, $\Omega_V$ is compatible with the group structure. 
%The two sheaves $\JJ^0_\Omega_V$  and 
%$\LB_\D$ together generate the sub-group-sheaf $\JJ^{\bullet}_{\Omega_V,\D}$
%of $\Omega_V$ isotropic sections. 
}
\end{example}

%************************************************************
% Foliations
%************************************************************
\subsection{Foliations}

In \cite{Hu-algebraic-surfaces}, 
when dealing with integrable systems of Jacobians, we considered  
Abel maps 
\[
A: \SS\rightarrow\J^1, 
\]
and it was shown that if 
the symplectic structure on $\J^1$ satisfied 
%one had an $A$ such that 
$A^*(\Omega)\wedge A^*(\Omega)=0$, then one could define a 
symplectic surface, which, in essence, encoded the integrable system. 
For our Prym varieties,
we will still consider the Abel map, and pull back the $V^*$-valued form 
$\Omega_V$.
Let us then 
% take an equivariant Abel map $A:\SS\rightarrow\J,$ and 
consider the rank of $A^*(\Omega_V)$, i.e., the smallest $k$ such that
$A^*(\Omega_V)^{\wedge k} = 0$ as a $(V^{* })^{\otimes k}$-valued $2k$-form.
In analogy to the Jacobian case, we will ask that 
this rank be everywhere minimal (without $A^*(\Omega_V)$ vanishing):

\begin{defi}
We say that our system has {\it rank two} if   
$$A^*(\Omega_V)\wedge A^*(\Omega_V)=0$$
as a section of $\Lambda^4(T^*(\SS))\otimes V^*\otimes V^*$. In other words, 
for a basis 
$e_i$ of $V$, if $\Omega_i$ denotes the contraction of $\Omega_V$ with $e_i$,
then 
$$A^*(\Omega_i)\wedge A^*(\Omega_j)=0$$
for all $i,j$.
\end{defi}

If the $A^*(\Omega_i)$ are non-vanishing,  a theorem of Darboux 
\cite{BCGGG} tells us that locally there exist
functions $x_i,y_i$ with non-vanishing differentials such that $A^*(\Omega_i)= dx_i\wedge
dy_i$; globally, for each index $i$ one has a codimension two foliation
(the {\it ``$\Omega_i$-foliation'' }, whose leaves
correspond locally to constant values of $x_i$ and $y_i$. 

We will make the corresponding genericity
assumption, and see in section \ref{sec-relaxing-the-genericity-conditions} 
how they can be relaxed. 

\medskip
\noindent {\sc Genericity Assumption A}:\ {\em The pull-back
$A^*(\Omega_V)$ is nowhere vanishing on $\SS$, and its
null-space is everywhere transverse to the curves $S_u$.}

\medskip
If $\Omega_V$ is non-zero, one can choose a basis of $V$ such that the $\Omega_i$ are non-zero. Next, we note:

\begin{new-lemma}
The span of $dx_1,..,dx_v, dy_1,...dy_v$ is at most  $(v+1)$-dimensional.
\end{new-lemma}

\noindent
{\sc Proof:} 
One can proceed inductively. For $v=1$, the proposition is trivial.
If the proposition holds for $(v-1)$, then, adding $dx_{v}, dy_{v}$ increases 
the span
by at most two dimensions. If the dimension increases by two, however, one 
finds 
that $dx_1\wedge dy_1\wedge dx_{v} \wedge dy_{v}\neq 0$, 
a contradiction.
\EndProof

\medskip
\noindent 
{\sc Genericity Assumption B:}
{\em 
The span of $dx_1,..,dx_v,$ 
$dy_1,...dy_v$ is everywhere $(v+1)$-dimensional. 
More invariantly, the null-space of $\Omega_V$ 
has codimension $(v+1)$ everywhere.
}

\medskip
Note that one has from the proof of the lemma that the spans of the subsets
 $dx_1,..,dx_s,$ $ dy_1,...dy_s$ are then $s+1$-dimensional.

With this assumption, one has a codimension $(v+1)$ foliation, given by intersecting
the null leaves of the foliations of the $\Omega_i$.
 The leaves of this foliation are then 
isotropic for $\Omega_V$, and we refer to it as the $\Omega_V$ foliation. 
 The fact that it is transverse to the curves
means that by restricting $U$ if necessary one obtains a global leaf space $X$ of dimension
$v+1$ into which 
the curves embed.  The form $\Omega_V$  descends  to $X$. 
Also, since the action of $W$ permutes the forms $\Omega_i$, 
the action of $W$ on $\bbs$ also descends to $X$.
More can be said about the structure of the form $\Omega_V$ on $X$: 

\begin{prop}
(a) Locally, under genericity assumptions A and B, there
exist forms $\phi_0, \phi_1,$ $ \phi_2,...,\phi_v$ on $X$ such that 
$A^*(\Omega_i)= \phi_0\wedge \phi_i$.

(b)  For $v>2$, the form $\phi_0$ defines a  codimension one
  foliation (the {\em ``$\phi_0$-foliation''})
which contains the leaves of the other foliations.
Locally, there are coordinates $\lambda_0, \lambda_1,...,\lambda_v$ such that 
the leaves of the $\phi_0$-foliation are given by level sets of $\lambda_0$,
and the forms $A^*(\Omega_i)=d\lambda_0\wedge d\lambda_i$.
\end{prop}

{\sc Proof:} Again, one proceeds inductively. For $v=2$, the fact that
$dx_1\wedge dy_1\wedge dx_2 \wedge dy_2=0$ implies that the $dx_1,dy_1$ plane
and the $dx_2,dy_2$ plane intersect in a line. We can then choose $\phi_0$ in that line,
and the result follows.
Suppose then that it holds for $(v-1)$. By the genericity assumption, 
$\phi_0, \phi_1,...\phi_{v-1}$ are independent, and $Span(\phi_0,..\phi_{v-1}, dx_{v}, dy_{v})$
is $(v+1)$-dimensional. Permuting $x$ and $y$ if necessary, one can then write
$dy_v = \sum_{i=0}^{v-1} a_i\phi_i + bdx_{v}$. One  has for all $j$ that
$$0= A^*(\Omega_j)\wedge A^*(\Omega_v)= \sum_{i\neq 0,j}a_i \phi_0\wedge\phi_j\wedge
\phi_i\wedge dx_v$$
forcing $a_i=0$ for all $i\neq 0$, so that $ dy_v = a_0\phi_0 +bdx_v$. 
The genericity
tells us that $a_0\neq 0$ and  one sets $\phi_v=-a_0 dx_v$.

For part (b), the form $\phi_0$ vanishes by construction on the leaves of the 
codimension two
foliations and on the leaves of the codimension $(v+1)$ foliation. 
It remains to be 
shown that it is integrable on $X$, that is that 
$d\phi_0 = 0\ ({\rm mod }\phi_0)$.
This follows, for $v>2$,  from 
$$0=dA^*(\Omega_i) = d\phi_0 \wedge \phi_i - \phi_0 \wedge d\phi_i,$$
as then $d\phi_0\wedge\phi_i\wedge\phi_0=0$ for all $i$, forcing 
$d\phi_0\wedge\phi_0=0$.

 We then note that the integrability of the $\phi_0$ foliation tells
us that there is a function $\lambda_0$ such that $\phi_0 = f\cdot d\lambda_0$
for some nowhere zero function $f$; this then allows us to write
$A^*(\Omega_i)= d\lambda_0\wedge \psi_i$, for some one-form $\psi_i$.
In a coordinate system $x_0=\lambda_0, x_1, x_2,..$, setting 
$\psi_i= \Sigma_{j=0}^va_{ij}(x) dx_j$, the fact that $dA^*(\Omega_i)=0$ 
implies that
$$ \partial_ja_{ik}= \partial_ka_{ij},\  j,k >0.$$
We would like to modify the $d\lambda_0$ term of $\psi_i$ so that it is closed;
such a modification leaves $A^*(\Omega_i)$ invariant. 
This involves solving simultaneously for $a_{i0}$ the equations
$$\partial_j a_{i0} = \partial_0a_{ij}, i = 1,..,v.$$
This is possible since 
$\partial_j\partial_0a_{ik}= \partial_k\partial_0a_{ij}$,
by the above. Once $\psi_i$ is closed, we can integrate it to obtain
the appropriate $\lambda_i$.
\EndProof

\medskip
We note that the $\phi_0$-distribution is only intrinsically determined when 
$v>1$. Intrinsically, the tangent spaces to the leaves are given by:
$$ Ker\left[ T^*\bbs {\buildrel \wedge A^*(\Omega_V) \over \longrightarrow} 
\Lambda^3T^*\bbs\otimes V^* \right].$$

Summarising, we have obtained:

\begin{thm}
\label{2.16} 
\begin{enumerate}
\item Under genericity assumptions A and B, 
restricting U if necessary, the $\Omega_V$ foliation has a 
global leaf space $X$  which  is a $v+1$ dimensional 
complex manifold, and into which the curves $S_u$ embed. 
\item The form $\Omega_V$ descends to $X$, and the group $W$ acts
on $X$, preserving $\Omega_V$. 
\item
\label{thm-item-folliations-are-mutually-transverse}
For any $a\in V$, one has a closed two-form $\Omega_a$ on $X$ obtained by contracting
$\Omega_V$ with $a$. The forms $\Omega_a$ have codimension two null foliations on $X$, and
the action of $W$ permutes these foliations in a natural way. Choosing a basis $e_i$ 
of $V$, so that one has forms $\Omega_i$, one has that the $\Omega_i$-foliations are 
all mutually transverse, and there are forms $\phi_0,...\phi_v$ such that
$\Omega_i = \phi_0\wedge \phi_i$. 
\item
For $v> 2$, the form $\phi_0$ defines a distribution which is integrable.
For $v>1$, whenever the distribution is integrable,   
the ``$\phi_0$''-foliation it
defines descends to $X$, and is $W$-invariant. 
One can write $\Omega_i = d\lambda_0\wedge d\lambda_i$, and the 
functions $\lambda_0,.., \lambda_v$ are  local coordinates for $X$, with 
$\lambda_0= cst.$
defining the leaves of the foliation..
\item
The two form $\Omega_V$ 
equivariantly identifies $V^*$ with the tensor product 
$T_L\otimes N_L$ of the kernel of $\phi_0$ (the tangent space to the leaves, in the 
integrable case)  
with its normal space (the normal space of the leaves).
\end{enumerate}
\end{thm}

%*********************************************************
% Fixed points and tangencies to the $\phi_0$-foliation.}
%*********************************************************
\subsection{Fixed points and tangencies to the $\phi_0$-foliation.}

\begin{prop}
 The $\phi_0$-foliation is generically transverse to
the generic curve $S_u$.
\end{prop}

{\sc Proof:} 
Let us suppose not. Then the curves $S_u$ are all contained in the 
leaves of the foliation, and the foliation is thus lifted from $U$. 
Locally, there is then a function on $U$, say $H$, such that $\phi_0=dH$. 
The forms $A^*(\Omega_i) = I(e_i)\circ \tilde A^*(\omega)$ then satisfy 
\begin{equation}
\label{2.18}
0=A^*(\Omega_i) \wedge dH = I(e_i)\circ \tilde A^*(\omega\wedge dH)
\end{equation}
This contradicts the non-degeneracy of $\omega$ on $\Prym$. Indeed, 
$\omega $ can be written as a sum $\Sigma dt_\mu\wedge dH_\mu$, with 
the linear coordinates $t_\mu$ on $\Prym$ corresponding to a basis 
$\psi_\mu$ for 
the $V^*$-valued invariant one-forms on $\SS$, and the $H_\mu$ are a basis 
for the 
cotangent space to the base. The condition (\ref{2.18}) forces $\psi_\mu=0$ 
for any $\mu$ such that $dH_\mu\wedge dH\neq 0$, a contradiction.
\EndProof

\medskip
 We define a {\it tangency point} of the curve $S_u$ in $X$ to be a point 
where it 
is tangent to the $\phi_0$-foliation, so that the form $\phi_0$ has no $d\rho$
component. The preceding proposition tells us that for the generic curve, 
the tangency points are isolated. 

A symplectic form gives an isomorphism between the normal bundle of a
Lagrangian submanifold and its cotangent bundle. In our case, we have 
a $V^*$-valued form, and this will allow a map between the 
normal bundle of the spectral curves and the tensor product of 
their canonical bundles with $V^*$.

\begin{eqnarray} 
N_S &\rightarrow & K_S\otimes V^*
\nonumber
\\
a   &\mapsto &  r(i(a)\Omega ),
\label{2.19}
\end{eqnarray}
where $r$ denotes restriction to the curve.   With respect to a local basis
of one-forms  as above, this is 
\begin{equation}
\label{2.20}
a\mapsto (\phi_i(a)r(\phi_0)-\phi_0(a)r(\phi_i)) \otimes e^i.
\end{equation}

Away from the tangency points, the map (\ref{2.19}) is an isomorphism, as one can 
use the isomorphism between the normal bundle and the tangents to the 
leaves of the foliation. At a tangency point of a curve $S$, this is no longer the case.
Indeed, at a tangency point, we have an exact sequence for the tangent
space $TL$ to the leaf: 
\begin{equation}
\label{2.21}
0\rightarrow T_S\rightarrow TL \rightarrow N_{S|L}\rightarrow 0
\end{equation}
Corresponding to this, we have an exact sequence for $V^*$, by using the
identification of Theorem \ref{2.16}:
$$0\rightarrow T_S\otimes N_L\rightarrow V ^*\rightarrow N_{S|L}\otimes N_L\rightarrow 0.
$$

\begin{new-lemma}
\label{2.23}
At a tangency point, the image of the map
(\ref{2.19}) lies in $K_S\otimes (T_S\otimes N_L)$.
\end{new-lemma}

{\sc Proof}: 
The proof is simply a  computation in local coordinates, 
using the formula (\ref{2.20}).
\EndProof

\medskip
The map (\ref{2.19}) then fails to be an isomorphism of bundles for $v>1$. 
On the other hand, on the 
level of global $W$-invariant sections over a generic curve $S_u$, we have:

\begin{prop}
\label{prop-isomorphism-of-vector-spaces-of-global-sections}
The map (\ref{2.19}) induces an isomorphism
$$H^0(S_u, N_S)^W= H^0(S_u,K_S \otimes V^*)^W.$$
and both of these spaces are isomorphic to the tangent space $T_uU$ of the
 base.
At any point, the evaluation map $H^0(S_u, N)^W\rightarrow N_p$ is surjective.
\end{prop}

{\sc Proof}:
We recall that the variety $\Prym$ is supposed to be an integrable system. 
This implies
that the dimension of the fibers of the Lagrangian fibration, which is the dimension of 
$H^1(S_u, {\cal O}\otimes V)^W$, is the same as the dimension of $T_uU$.
The tangent space $T_uU$ of the base is identified with a subspace of 
$H^0(S_u, N_S)^W$.
In turn, the map (\ref{2.19}) gives an identification of $H^0(S_u, N_S)^W$
with a subspace of $H^0(S_u,K_S\otimes V^*)^W$. The latter space is however 
Serre dual to $H^1(S_u,{\cal O}\otimes V)^W$, which forces each 
subspace to be the full space. The surjectivity of the evaluation map follows from the fact that
$X$ is a quotient of the space $\bbs$, and on this space, the evaluation 
map is trivially surjective.
\EndProof

\medskip
The geometry of the curves with respect to the $\phi_0$-foliation is
thus tightly linked to the vanishing of forms in $H^0(S_u,K\otimes V^*)^W$:
for points $p$ which are not tangent, one has that for any non-zero $e$ in $V$,
the  map $H^0(S_u,K\otimes V^*)^W\rightarrow K_p$ given by 
contraction with $e$ and evaluation at $p$ is surjective.

The group invariance will force certain components of the form to vanish 
at the points with non-trivial stabiliser. We first examine the nature of 
these stabilisers.

\begin{new-lemma}
\label{2.26}
Assume genericity conditions A and B, and let $v\ge 2$. 
For  the generic curve in the fibration,
 the stabilisers in $W$ of 
the points on the curve are either trivial or $\bbz/2$. At the points
$S_u\subset X$ where
the stabiliser is $\bbz/2$, the action of $\bbz/2$ on the tangent space
has weights $(1,0,0,0...,0)$, and the weight $1$ space  is tangent to the 
curve. 
\end{new-lemma}

{\sc Proof}: 
Let  $p\in S_u\subset \bbs$, and
let $G\subset W$ be its  stabiliser. The fact that the $W$-action on 
$S_u$ is generically free
tells us that $G$ is cyclic, of order, say,  $n$; let $g$ be the generator 
which acts on 
the curve by $z\mapsto exp(2\pi i/n)z$, in an appropriate coordinate system 
centered at $p$. A simple argument (a contour integral of $dln(g(z)-z)$)
 tells us that if $g$
has a fixed point at $p$ on $S_u$, it has a fixed point near $p$ 
for nearby curves $S_{u'}$.
We can then choose a generic curve, for which the order of stabilisers are 
minimal; the 
  stabilisers are then constant for nearby curves, and the generator $g$ has  
 as local fixed point locus a codimension one submanifold
 transverse to the curves.

 At $p$, the tangent space to $\bbs$ is a  $G$-representation, and we can 
decompose 
it  as a sum of
 weight spaces $R_{j}, j= 1,..,v$, on which $g$ acts by 
$exp(\alpha_j2\pi i /n)$,
with $\alpha_j\in\{0,...,n-1\}$. We note that by restricting to the 
curve $S_u$ we see that at least one  of the 
$\alpha_i$, say $\alpha_1$,
 is one.
For the generic curve, the vector of weights $\alpha_j$ 
is, up to ordering, $(1,0,...,0)$, since there is a codimension one fixed 
locus.

The action of $G$ preserves the null-foliation and descends to the space $X$, 
and,
again, acts with weights $(1,0,..,0)$, as the curve embeds in $X$.
 By taking suitable averages, one can choose
 coordinates $z_0,..,z_v$  so that the action
of $g$ is given by $(z_0,z_1,..,z_v)\rightarrow (
exp(2\pi i/n)z_0,z_1,..,z_v)$; when the foliation is integrable, the 
coordinates 
$z_i$ can be chosen to coincide with the coordinates $\lambda_i$, in a 
suitable order.
 Now let $v\ge 2$. The kernel of  $\phi_0$ 
at  $p$ (tensored with the normal line)
is isomorphic to $V$ as a $G$-representation. There are two 
possibilities: the first is that the action of $G$ on the normal line to the 
foliation
has weight zero, and the second is that it has weight one. In the first 
case, the weights on the tangent space to the leaf are $(1,0,0,...,0)$; in the 
second, $(0,0,..,0)$. This in turn gives for the representation on 
$V$ the two possibilities $(1,0,..,0)$, $(1,1,...,1)$.
 Now $V$ is an {\it integral} 
representation, and so weight spaces with weight $i$ must have
counterpart weight spaces with weight $n-i$, for $i\neq (n/2)$.
This forces $n=2$, and $G=\bbz/2$.
\EndProof

{\sc Remark}: We note that there are two types of points at which the 
action is not free on the generic curve: the first (type I) is tangent to the 
$\phi_0$-foliation, and the second (type II) is transverse to it. 
For fixed points of the second type, one has that $g$ is mapped to $-1$ in 
$Aut(\chi)$; as our representation is faithful, there can only be one such 
$g$, and  $g$ is then in the center of $W$.
% Let $Z=\bbz/2$ be the subgroup of $W$ spanned by $g$ if there is such a 
% type II point, and  let $Z=0$ otherwise. 

%**************************************************************
% Relaxing the genericity conditions
%**************************************************************
\subsection{Relaxing the genericity conditions} 
\label{sec-relaxing-the-genericity-conditions}

The conditions that we gave above are not strictly necessary to have a nice 
quotient space,
and indeed we can relax them a bit. One modification which allows Theorem
\ref{2.16} to go through 
essentially unchanged is as follows: 

\medskip
{\sc Genericity condition \^B}: 
{\em
The null-space of $\Omega_V$ is  of 
codimension $v+1$ over a dense open set $O$. Over the set $O$, this null 
space defines
a vector bundle which extends to a globally defined $W$-invariant subbundle 
of the tangent bundle.
similarly, the sub bundle of the tangent bundle defined over $O$ as the kernel
of $\phi_0$ extends to all of $\bbs$ as a $W$-invariant subbundle of the 
tangent bundle.
}

\medskip
With this modified condition, Theorem \ref{2.16} goes through, 
(apart from part \ref{thm-item-folliations-are-mutually-transverse},
which is only valid over the image of the open set $O$), with the change 
that the form $\Omega_V$ on $X$ is then allowed to degenerate away from the 
image of the set $O$ in $X$.

While this new condition is sufficient to give us $X$, it allows a wide range
of degeneracies in the group action, in contrast to the rather 
strict constraints imposed on it by condition $B$, as in lemma \ref{2.26}.
These degeneracies are analysable in a straightforward way in any 
particular case; however, the variety of possibilities that they open up is 
rather large and intractable when considered in full generality. 
We have opted for a partial relaxation of condition $B$, 
which allows some interesting examples (see Example 
\ref{example-regular-coadjoint-orbits-with-sub-regular-semisimple-part}), 
while maintaining some control over the group action.

To motivate our choice, note that at the points $p$ 
(with a $\bbz/2$-stabiliser) of type I,
one can split $V^*$ into $\bbz/2$-weight spaces as $V^*_0\oplus V_1^*$.
The space $V^*_1$ is one dimensional, and is the subspace $T_S\otimes N_L$ of 
(\ref{2.21}) 
corresponding to $TS$. The space $ V_0^*$ has dimension $v-1$.
In a corresponding fashion, any $W$-invariant $V^*$-valued 1-form $\phi$
on a curve $S_u$ 
can be split near $p$ as  a sum $ \rho_0\otimes e_0 + \rho_1\otimes e_1$,
with $\rho_i$ one-forms and $e_i\in V_i$, so that
 $g^*(\rho_i) = (-1)^i\rho_i$, with $g\neq 0, g\in \bbz/2$. In particular,
$\rho_0$ vanishes at $p$.

Examples of integrable systems indicate that in addition to types I and II
there is a third type of points $p$ on the generic curve with stabilizer 
$\bbz/2$ (their closure is a divisor).
For these points, the weights of the $G$ action on the tangent space to the 
$\phi_0$-foliation are  $(0,0,..,0)$, so that the fixed point locus of 
$G$ coincides with a leaf, as for type II, with
the action on the normal space having weight $1$.
On $V$, however, the action  has arbitrary weights 
$(1,,..,1,0,..,0)$, with $k>1$ non-zero weights. 
(write the  corresponding basis as $e_1,..,e_v$).  The 
form $\Omega_V$ can then be written as $\sum \phi_0\wedge\phi_i\otimes e_i$.
The group invariance then forces $\phi_i, i = k+1,..,v$ to vanish along the 
components of the divisor of this new type. 
We say that such points are {\em of type III} \ if 
this vanishing is minimal, i.e., if $\phi_0,.., \phi_j$ are non-zero and 
independent, while 
$\phi_i$, $i=j+1,..,v$ have a simple zero, but are independent at first order. 
More precisely, if $\Delta$ is the divisor of type III points, 
splitting $V$ into $V_0\oplus V_1$ along $\Delta$, we want that the 
contraction map by $\Omega_V$ 
$$T\SS  \otimes V_0
\rightarrow T^*\SS$$
vanish along $\Delta$, but that the map
$$T\SS  \otimes \left[
V_1\ \ \oplus \ \  
\left(
V_0\otimes 
{\cal O}(\Delta)
\right)
\right]
\rightarrow T^*\SS$$
be of constant rank $v\!+\!1$.

\medskip
{\sc Regularity condition B$'$:} 
{\em 
We first ask that the stabilisers $G$ of point on the 
generic curve in the system all be  either $0$ or $\bbz/2$. At points with 
trivial stabiliser, we ask that the null-space of $\Omega_V$ have 
codimension $v\!+\!1$. Equivalently, one asks that the contraction by 
$\Omega_V$:
$$T\SS  \otimes V
\rightarrow T^*\SS$$
have constant rank $v\!+\!1$.
At points with $\bbz/2$-stabiliser, we ask that 
either the null-space of $\Omega_V$ be again $v\!+\!1$-dimensional 
(type I or II points), or be of type III.
}

\medskip
We can summarise the situation for points with non-trivial stabiliser as 
follows:
$$\matrix{ {\rm Type}& {\rm Weights\ on\ }TL& {\rm Weights\ on}\ V& {\rm Curve\
 and\ leaf}&{\rm  Fixed\ point\ locus\ and\ leaf}\cr\cr
 {\rm I} & (1,0,...,0)& (1,0,...,0)& {\rm Tangent}& {\rm Transverse}\cr
{\rm II} & (0,0,...,0)& (1,1,...,1)& {\rm Transverse}& {\rm Coincide}\cr
{\rm III} & (0,0,...,0)& (1,.., 1,0,...,0)& {\rm Transverse}&{\rm Coincide}\cr}$$

\begin{example}
\label{example-regular-coadjoint-orbits-with-sub-regular-semisimple-part}
{\rm
Prym integrable systems with points of type III
arise from certain reduced coadjoint orbits of the positive half of the 
loop algebra $L\LieAlg{g}$ (see section \ref{sec-coadjoint-orbits}). 
Specifically, consider the coadjoint orbit of 
$N(\lambda)\cdot p(\lambda)^{-1}$ 
(\ref{eq-element-in-coadjoint-orbit}) 
such that all roots $\{\lambda_i\}$ of the polynomial 
$p(\lambda)$ are simple and the values $N(\lambda_i)\in \LieAlg{g}$
are either 1) regular and semi-simple, or 2) regular with a sub-regular 
semi-simple part (in case $\LieAlg{g}=\LieAlg{sl}_n$, $n\geq 3$, the latter
$N(\lambda_i)$ have one eigenvalue with multiplicity $2$ but an $n\!-\!1$
dimensional centralizer). 
In this case, $\Sigma$ is $\PP^1$ and all spectral curves $S_u$ will have 
simple ramification points over 
$\lambda\in \PP^1$ where $N(\lambda)$ is regular with a 
sub-regular semi-simple part. It is a ramification point of type III if 
$\lambda$ is either a root of $p$ or $\lambda=\infty$ and of type I otherwise. 
A further relaxation of condition B, which includes more general 
coadjoint orbits in $L\LieAlg{g}$,  is considered in section
\ref{sec-from-hilbert-schemes-to-acihs}.
}
\end{example}

%*************************************************************
% Rigidity
%*************************************************************
\section{Rigidity}
\label{sec-rigidity}

This section is devoted to proving that the $\phi_0$-foliation allows us to define a uniform 
quotient of the curves $S_u$ under the action of the group $W$. We will need a regularity 
condition to prove this. While the condition seems rather artificial, we will see in the next section
that it is implied by some rather more natural conditions on the curves.

\begin{condition}
\label{cond-regularity}
{\sc (Regularity Condition)}
For a curve $S_u$, at each point $p$ with trivial stabiliser, 
there is a two-dimensional subspace of $V$ such that for each non zero vector 
$e$ in the subspace, the map $H^0(S_u,K\otimes V^*)^W\rightarrow K_p$
given by evaluation and contraction with $e$ is surjective.
For points $p$ of type I with non-trivial stabiliser $G =  \{1,g\}$, 
there is a one-dimensional subspace of the $g$-invariant subspace $V_0$ such
 that contraction with a non zero element of the subspace with 
elements of $H^0(S_u,K\otimes V^*)^W$ gives a form with a zero of order one.
\end{condition}

\medskip
In some sense this regularity condition tells us that forms only vanish
to the order forced by group invariance.
We note that the  condition is open, so the fact that it holds for one 
curve implies that it also holds for nearby curves.

\begin{prop} 
\label{prop-leaf-intersects-curve-in-a-single-orbit}
Suppose that our rank two system satisfies 
genericity conditions A and B$'$, as well as the regularity condition 
\ref{cond-regularity} at the
curve $S_u$. Suppose that the curve $S_u$ is generic, in that the points of 
$S_u$ have only trivial or $\bbz/2$ stabilisers.
There is an open set $U'$ containing $u$ and a covering of $S_u$ in $X$ by 
open sets
$B_\rho$ such that the intersection of each leaf of the $\phi_0$-foliation in 
$B_\rho$ with $S_{u'}, u'\in U'$ is non-empty and consists of elements 
belonging to a single $W$-orbit on $S_{u'}$.
\end{prop}

{\sc Proof:} 
Away from the tangency points, there is no difficulty in finding appropriate
open sets $B_\rho$. The problem occurs at the tangency points; our aim will 
be to show 
that tangency occurs only to the order forced by group invariance, that is 
only at type  points with $\bbz/2$ stabiliser.

Our regularity condition \ref{cond-regularity}
tells us that at points with trivial stabiliser, 
the image of $N$ in $K_S\otimes V^*$ is two dimensional, and so by lemma 
\ref{2.23}, this cannot be a tangency point. 
Thus, the points
of tangency all have non-trivial stabiliser, which then must be a $\bbz/2$.
If $z$ is a coordinate on the curve, let  the embedding of the curve at the 
tangency point be given  by 
$$z\mapsto 
(\lambda_0(z),..,\lambda_v(z)) = (z^{k_0}, a_1z^{k_1},...,a_vz^{k_v}) + 
({\rm higher \ order)},$$
with the $\phi_0$-foliation cut out by level sets of the first coordinate.
Tangency means that $k_0>1$, and $k_0$ is the order of the tangency.
We will suppose that the coordinates in $X$ correspond to a decomposition of 
the tangent space into weight spaces.   
 The  discussion of the weights given above tells us that 
$k_i $ is even and at least two, for all of the $i>0$ except one, say $k_1$, 
which must be one.

\begin{new-lemma}
For the generic curve,  let $s$ be the image in
$K\otimes V^*$ under the map (\ref{2.19}) 
of a local section of the normal bundle $N$ , and let $s$ 
be of the form $(z^adz+..)\otimes e + ({\rm higher\ order})$, where 
$e$ is a vector on which $G$ acts trivially. If the order of the tangency 
locus is $k_0$, then $a\ge k_0-1$.
\end{new-lemma}

{\sc Proof:} 
The proof is a simple computation: one has that the image of the normal vector
${\partial\over \partial \lambda_0}$ in $K\otimes V^*$ is 
$-\Sigma_{i>0} a_ik_iz^{k_i-1}dz\otimes
e^i+$ \ higher order, and that of ${\partial\over \partial \lambda_i}$ is 
$k_0z^{k_0-1}dz\otimes e^i+$ \ higher order. The leading order term of the 
image of the 
of ${\partial\over \partial \lambda_0}$ then corresponds to the vector $e^1$, 
which has odd
weight; for the leading term to have even weight, one must consider the image 
of the vectors ${\partial\over \partial \lambda_i}, i>1$. 
\EndProof

\medskip
As a consequence of this and the regularity condition \ref{cond-regularity},  
one has $k_0=2$
at a tangency point. This then tells
us that  the intersection of the curve with nearby leaves of the 
$\phi_0$-foliation 
consists of two points which lie in the same $\bbz/2$-orbit, which is what we 
needed to complete the proof of proposition 
\ref{prop-leaf-intersects-curve-in-a-single-orbit}.
\EndProof

\begin{thm}
\label{thm-rigidity}
For rank two systems satisfying the genericity conditions
A and B$'$ as well as the regularity condition \ref{cond-regularity}, 
restricting $U$ if necessary,  
the corresponding space $X$ has $S_u/|W|$
 as a simultaneous quotient for the $\phi_0$-foliation and the group action. The curve
$\Sigma = S_u/|W|$ is then independent of $u$.

The quotient $Y = X/W$ is smooth  near the curve $S_u/|W|$, and $Y$ fibers 
 over the base $\Sigma$. The curves 
$S_u/|W|$ are then sections of this fibering, and 
this realises $U$ as a space of sections of $Y\rightarrow\Sigma$. 

The form  $\Omega_V$ is isotropic on the fibers of the projection of $X$ to $\Sigma$,
and defines an isomorphism between the tensor product of the tangent bundle of the 
base with $V$ and the cotangent bundle to the fiber, away from the points with 
non trivial stabiliser of type II or III. 
\end{thm}

{\sc Proof:}
One defines the quotient on the union of the open sets $B_\rho$, 
mapping each leaf to the orbit of $W$ on $S_u$ that it intersects. 
This simultaneously quotients
the nearby curves $S_{u'}$ by the action of $W$, giving a uniform quotient 
for the curves.

The smoothness of $Y$ follows from the fact that 
%the only stabilisers are $\bbz/2$.
at any point $p$ with non-trivial stabilizer, the non-trivial character of 
the stabilizer appears in $T_pX$ with multiplicity $1$.
\EndProof

\medskip
We close by noting the following infinitesimal version of the theorem above, 
whose proof would, regrettably, take us too far afield.

\begin{prop}
Let $q: X\rightarrow Y= X/W$ be the quotient map, and let 
$C_u\subset Y$ be the smooth  quotient of a smooth curve $S_u\subset X$. 
Let $V$ be a vector bundle on $X$ of rank ${\rm dim}(X)-1$ and $\Omega_V$ 
a $V^*$-valued two-form on $X$ satisfying

1. $\Omega_V\wedge \Omega_V = 0$,

2. The cokernel of the contraction mapping 
$TS_u\otimes V\rightarrow N^*_{S|X}$ is 
supported on a subscheme of the ramification divisor $\Delta$ of 
$S\rightarrow C$.

Then the exact sequence
$$0\rightarrow TC\rightarrow TY_{|C} \rightarrow N_{C|Y}\rightarrow 0$$
splits canonically.
\end{prop}

One virtue of this proposition is that it highlights the role played in 
rigidity by the rank-2 condition. For example, one can have situations in 
which the same Prym variety can correspond to a family of different curves, 
with non-constant quotient.
The integrable system one could build from such a family would contradict 
rigidity, but for the fact that it does not have rank 2.

%******************************************************************
% Genericity Conditions
%******************************************************************
\section{Genericity Conditions}
\label{sec-genericity-conditions}
In section \ref{sec-weak-and-strong-BPF} we reformulate all our
genericity and regularity conditions in terms of a single algebro-geometric 
assumption  -
the Strong Base Point Freeness Assumption 
\ref{assumption-strong-generation}. 
The normal bundle $N_{S_u}$ of $S_u$ in $X$ turns out to 
be an intrinsic vector bundle 
(\ref{4.6}) which depends only on the $W$-action (and is independent of $X$). 
%There is an intrinsic vector bundle 
%(\ref{4.6}) on each curve $S_u$, which depends only on the $W$-action,  
%and is a candidate to become the normal bundle $N_{S_u}$ of $S_u$ in $X$
%if $X$ can be constructed. 
Assumption \ref{assumption-strong-generation} requires this  intrinsic 
vector bundle to be generated by its global sections. 
In section \ref{sec-geometric-conditions} we prove that  
the Strong Base Point Freeness Assumption 
holds provided the covers $S_u\rightarrow S_u/W$
have sufficiently many branch points (Theorem 
\ref{thm-SBPF-holds-if-many-branch-pts-exist}). 

%******************************************************************
% Genericity, regularity  and invariant forms
%******************************************************************
\subsection{Genericity, regularity  and invariant forms}
\label{sec-weak-and-strong-BPF}

The 2-forms we are considering on the variety $\SS$ all arise as pull-backs
of forms via an Abel map to a family of Jacobians, and eventually, to a Prym 
variety, whose cotangent bundle is trivial, and is generated by global 
invariant $V^*$-valued
one-forms on the curve. Furthermore, we are looking at a map into an integrable system: the coordinates 
normal to the leaves of the integrable system are pulled back to coordinates on $\SS$, 
while the coordinates on the leaves correspond to Abelian integrals of 
global one forms on the curves
in  $\SS$. It is then not surprising then that the non-degeneracy of the 
forms on $\SS$ is linked to the non-vanishing of global ($V^*$-valued) one-forms
on the curves.

For concreteness, let us write out some of our forms in coordinates. 
Let $H_1,...H_d$ be coordinates 
lifted from $U$ to $\SS$, and complete them with a coordinate $\rho$  along the curves. 
On $\bbj$, we can we can complete the $H_j$ with coordinates
 $t_i^j$, linear along the 
fibers, so that $\Omega^i = \sum_jdt^j_i\wedge dH_j$. The pull-back to the 
curves $S_u$ of these linear forms will be global holomorphic one-forms
$f_i^j d\rho$. Therefore, on $\bbs$
\begin{equation}
\label{4.1}
A^*(\Omega_i) = \sum_j (f_i^j d\rho\wedge dH_j) + 
\sum_{j,k} ({{\partial t_i^j}\over {\partial H_k}}-{{\partial t_i^k}\over {\partial H_j}})
dH_k\wedge dH_j = \phi_0\wedge \phi_i
\end{equation}
The genericity condition A is that $A^*(\Omega_V) =
 A^*(\Sigma_i(\Omega_i\otimes e^i))$ have a non-vanishing $d\rho$-component, so
 that, 
at each point, at least
one of the forms $f_i^j d\rho$ is non-vanishing. 
(We note in passing that our rank 2 condition essentially tells us that we can 
choose our coordinates $\rho, H_j$ in such a way that the $dH\wedge dH$ terms in (\ref{4.1})
vanish) 

Now we use the fact that our form $\Omega_V$ arises from a form on $\Prym$: 
recall that the cotangent bundle to the Prym variety is spanned by
linear forms $dt_\alpha$ corresponding to the $W$-invariant $V^*$-valued 
1-forms $\phi_\alpha$ on $S_u$. Choosing suitable coordinates $H_\alpha$
on $U$, one has that
$$A^*(\Omega_V) = \sum_\alpha\phi_\alpha\wedge dH_\alpha  + (dH\wedge dH),$$
so that condition A is then equivalent to the non-vanishing of one of the 
$\phi_\alpha$. 

We can then rewrite this as 

\medskip
\noindent
{\sc Weak base point freeness assumption}:
{\em 
At any point $p\in S_u$, the evaluation map
$$H^0(S_u, K_{S_u}\otimes V^*)^W \otimes V \rightarrow (K_{S_u})_p$$
is surjective.
}

\medskip
Geometrically, if we define the Abel map:
$$
\tilde A: \bbs{\buildrel A\over \longrightarrow} 
\bbj{\buildrel \otimes Id \over 
\longrightarrow}\bbj\otimes_\bbz\chi\otimes_\bbz\chi^* {\buildrel Av\over
\longrightarrow}\Prym\otimes_\bbz \chi^*,
$$
where $Av$ is some multiple of averaging under the group action (one takes a 
multiple to land in the connected component of the identity),
then genericity condition A is equivalent to the Abel map $\tilde A$ being an 
immersion.

\medskip 
{\sc Remark:} In the case studied in \cite{Hu-algebraic-surfaces}, 
that of a trivial group, the genericity conditions are 
automatically satisfied, simply as a consequence of the fact that at any 
point on a Riemann surface of genus $\geq 1$ there is a global holomorphic 
one-form that does not vanish. Equivalently, the Abel map is an immersion.

We now turn to our other genericity assumptions: we will give a condition 
which implies them all.
We will assume that all the branch points of $\SS \rightarrow\SS/W$ 
are of order two, with stabilisers $G_x=\bbz/2$.

Recall that the set of branch points forms a smooth divisor $\Delta$,
intersecting $S_u$ in a divisor $\Delta_u$. At $x\in \Delta$,
we split $V$ into weight spaces $(V_x)_0\oplus (V_x)_1$ for the $G_x$ action.

\begin{assumption}
\label{assumption-strong-generation}
{\sc (Strong base point freeness)}:
The vector bundle 
\begin{equation}
\label{4.6}
{\rm ker}[K_{S_u}\otimes V^*\rightarrow 
\oplus_{x\in \Delta_u} (K_{S_u}\otimes (V_x^*)_0)_x]
\end{equation}
is generated by its $W$-invariant global sections.
\end{assumption}

\medskip
In terms of the preceding coordinates, writing 
$A^*(\Omega_V) = \sum_\alpha\phi_\alpha\wedge dH_\alpha  + (dH\wedge dH),$
the strong base point condition tells us that away from branch points,
the evaluation map 
$${\rm span}\{\phi_1,...\phi_k\}\rightarrow (K_{S_u}\otimes V^*)_x$$
is surjective. At the branch points, recalling that $K_{S_u}$
has weight one, and  choosing a basis so that $\phi_1,...,\phi_j$
has weight zero, and $\phi_{j+1},...,\phi_k$ weight one, the evaluation 
map 
$${\rm span}\{\phi_{ 1},...\phi_j\}\rightarrow 
(K_{S_u}\otimes (V^*_1)_x)_x$$
is surjective, while, as  $\phi_{j+1},...\phi_k$ have a simple zero at the 
branch point,
choosing a coordinate $z$ so that the branch point is $z=0$, the composition
of the evaluation map with division by $z$ 
$${\rm span}\{\phi_{j+1},...\phi_k\}\rightarrow 
(K_{S_u}\otimes (V^*_0)_x)_x$$
is also surjective. Another way of saying this is as follows: we have a map of 
contraction by $\Omega_V$
\begin{equation}
\label{4.10}
N_{S_u}\rightarrow K_{S_u}\otimes V^*
\end{equation}
inducing an isomorphism on global sections, by Proposition
\ref{prop-isomorphism-of-vector-spaces-of-global-sections}. 
The strong base point freeness
assumption is that the image of these sections only vanish at a point 
$p$ to the order that is forced by invariance under the isotropy group. 

The statement, that invariant sections of $K_{S_u}\otimes V^*$ vanish only to 
the order imposed by their invariance, admits the following sheaf-theoretic 
translation:
%A direct translation of this statement is most easily given if we 
Push forward $K_{S_u}\otimes V^*$ via the quotient map 
$q:S_u \rightarrow S_u/W$ and consider its $W$-invariant sub-bundle. 
Then the strong base point freeness assumption translates to: 
{\em The vector bundle }

\noindent
\hspace{1ex} 
\hfill
$(q_*[K_{S_u}\otimes V^*])^W$
\hfill
(\ref{4.6}')

\medskip
\noindent
{\em is generated by its global sections.} \ Note that the vector bundle 
(\ref{4.6}) is simply the pull-back of (\ref{4.6}') to $S_u$. 
The formulation in terms of  (\ref{4.6}') enables one to consider 
the case of cyclic stabilizers of higher order
(see section \ref{sec-from-hilbert-schemes-to-acihs}).

\begin{prop}
\label{prop-SBPF-equivalent-to-genericity-and-regularity-cond}
The Strong base point freeness assumption is 
equivalent to the conjunction of the following assumptions:
\begin{trivlist}
\item [1.] The Weak base point freeness assumption,
\item[2.] Genericity Condition B$'$,
\item [3.] The  Regularity Condition \ref{cond-regularity}.
\end{trivlist}
\end{prop}

{\sc Proof:} $\Rightarrow$: Let us assume Strong base point freeness.
We will prove the three conditions in two steps: In Step I we will prove 
them over a dense open subset $U_0$ of $U$. 
In Step II we will extend the result over the whole of $U$. 

\noindent
\underline{Step I}: 
The representation $V$ is supposed to be faithful, so at least one of the 
weights is non-zero. 
Thus, at every point, there is at least one non-vanishing $\phi_\alpha$,
and so one has Weak base point freeness.  

In turn, this tells us that the image of contraction with $\Omega_V$
\begin{equation}
\label{4.12}
T\SS  \otimes V
\rightarrow T^*\SS
\end{equation}
is transverse to the hyperplane $N^*_{S_u} =T^*U\subset T^*\SS$, 
so that the image has dimension
$r$ iff  its intersection with $N^*_{S_u}$ has dimension $r\!-\!1$. 
The restriction of  (\ref{4.12}) to the vertical subbundle $TS_u\otimes V$ 
has image in $N^*_{S_u}$: in
 the coordinates above, the image of $x\in TS_u\otimes V$ is 
$\sum_\alpha \phi_\alpha(x) dH_\alpha$. Referring to (\ref{4.6}) above, the image 
of $TS_u\otimes V$ is
$v$-dimensional away from the branch points, giving condition B$'$ away from
the branch points.

Thus, away from the branch points, all the structure we have developped in 
section \ref{sec-prym-varieties} exists, 
in particular, 
%the null foliation, given by the subbundle of $T\SS$
%which vanishes under contraction with $A^*(\Omega_V)$, and 
we get the $\phi_0$-distribution,
given by the  sub line-bundle of $T^*\SS$ obtained as the kernel of 
$\wedge A^*(\Omega_V): T^*\SS\rightarrow\Lambda^3T^*\SS\otimes V^*$. 
This line-bundle extends as a subbundle 
to the complement of a co-dimension 2 subset of $\SS$, 
and so over a dense open subset $U_0$ of $U$. In this first step we
proceed to work over $U_0$. 
%to an open dense set 
%$\Delta_0$ of the divisor $\Delta$. The foliations also extend to the same 
%set. (Note that if $v=2$, the $\phi_0$- foliation is not 
%necessarily integrable; if it is not , we will choose locally an invariant 
%foliation whose tangent space coincides with the distribution at the point 
%under study, and call that the $\phi_0$ foliation.) 

 Let $TL$ be the co-rank $1$ subbundle of $T\SS$ corresponding to 
the $\phi_0$-distribution. Denote by
$\Delta_{reg}$ the points of 
% $\Delta_0$ 
$\Delta$ for which 
$\Delta$ and the $\phi_0$-distribution 
are transverse, and $\Delta_\infty$ the set for which they are not.
%Along $\Delta_\infty$, the group invariance forces $\Delta_\infty$ and a 
%leaf of the $\phi_0$-foliation to coincide. 

{\it Case i) $\Delta_{reg}$:} Let $x\in \Delta_{reg}$. 
%Let $L$ be the leaf of the $\phi_0$-foliation through $x$.  
As $TL$ is transversal to $\Delta_{reg}$,
the action of $G_x$ on $TL$ has weights $(1,0,..,0)$, and the curve is 
tangent to $TL$ at $x$. Consider the diagram along $S_u$:
\begin{equation}
\label{4.13}
\matrix {TS_u\otimes V&\rightarrow&N^*_{S_u}& \cr \downarrow& & \downarrow\cr
N_L\otimes V&\rightarrow&T^*L},
\end{equation}
where the horizontal maps are contraction by $\Omega_V$, and $N_L$
%, $T^*L$ are the normal bundles and the cotangent bundles of the foliation.
is the quotient line-bundle $T\SS/TL$ (the relative normal bundle to the 
$\phi_0$-foliation in case $TL$ is integrable). 

 The top horizontal map is the dual of (\ref{4.10}), 
and so is injective at $x$ on the $TS_u\otimes (V_x)_1$ summand, and vanishes
to order 1 on the  $TS_u\otimes (V_x)_0$-summand. 
The left hand vertical map vanishes at 
$x$, as the curve is tangent to $TL$. Thus, the induced map 
$TS_u\otimes V\rightarrow 
T^*L$ must vanish, and so, as the map $N^*_{S_u}\rightarrow T^*L$ has rank 
$v-1$ at $x$, the weight space
$(V_x)_1$ must be 1-dimensional. Furthermore, the map 
$TS_u\otimes V\rightarrow N^*_{S_u}\rightarrow  T^*L$ vanishes to order one, 
and so the map $TS_u\rightarrow N_L$ must also vanish to order one; 
the tangency of the curve with the $\phi_0$-foliation is of minimal order. 
This in turn implies that the determinant of the map  
$N^*_{S_u}\rightarrow T^*L$ vanishes to order one at $x$. 

One then has that the map $TS_u\otimes V\rightarrow T^*L$ vanishes to order 
one at $x$, and, if one divides the map by a coordinate function vanishing 
at $x$, one obtains a rank $v$ map. In turn this tells us that the map 
$N_L\otimes V\rightarrow T^*L$ at $x$ must have rank $v$ at $x$. 
From the diagram
$$
\matrix{& & N^*_L\cr& & \downarrow \cr T\SS\otimes V&\rightarrow& T^*\SS\cr
\downarrow& &\downarrow\cr N_L\otimes V&\rightarrow&T^*L},
$$
the image of $T\SS\otimes V$ in $T^*L$ is $v$-dimensional; but the line 
$N^*_L$ also lies in the image of $T\SS\otimes V$, giving $v+1$ dimensions in 
all in $T^*\SS$.

{\it Case ii) 
$\Delta_{\infty}$:} In this case 
%the leaf 
$TL$ and the tangent bundle of the component of $\Delta$
coincide, and the curve is transverse to both. The vertical maps in 
(\ref{4.13}) are isomorphisms,
and so the map $N_L\otimes V \rightarrow T^*L$ vanishes to order one on the 
$(V_x)_0$
summand, and is injective on the $(V_x)_1$ summand. 
This in turn gives us an $\Omega_V$
consistent with a type II or type III point in condition B$'$.

Finally, the regularity condition \ref{cond-regularity} is an 
automatic consequence of the Strong base 
point freeness, as it is simply a weakening of it.

\medskip
\noindent
\underline{Step II}: 
We have established the assumptions of 
Theorem \ref{thm-rigidity} over the open subset $\SS_0$ of $\SS$ 
which is the union of 1) the complement of $\Delta$ and 2) 
the inverse image of $U_0$. 
%where the $\phi_0$-distribution extends as a subbundle of $T\SS$. 
Let $\pi:\C\rightarrow U$ be the quotient family $\SS/W$. 
Theorem \ref{thm-rigidity} implies that the $\phi_0$-distribution $TL$ 
descends to give a splitting $f_{\Omega_V} : T\C \rightarrow T_{\C/U}$ of 
the short exact sequence
\[
0 \rightarrow 
T_{\C/U} \rightarrow
T\C \rightarrow 
\pi^*TU \rightarrow 0
\]
over $\SS_0/W$. The complement $\SS\setminus \SS_0$ has codimension $2$. 
Hence, the homomorphism $f_{\Omega_V}$ extends to a regular homomorphism 
over the whole of $\C$. By continuity, the extended $f_{\Omega_V}$ is 
also a splitting. Hence, $f^*_{\Omega_V}$ embeds 
$T^*_{\C/U}$ as a line-sub-bundle of $T\C$. 
It is now easy to see that $(q^*T^*_{\C/U})(\Delta_{\infty})$ 
embeds as a line-sub-bundle of $T^*\SS$ extending $N^*_L$. It follows
that  the $\phi_0$-distribution extends as a subbundle of $T\SS$
over the whole of $\SS$. Returning to Step I 
we establish the three assumptions over the whole of $\SS$.

\medskip
\noindent
$\Leftarrow$: 
We saw in the previous sections that the three conditions enabled us to 
prove that there was a well defined codimension one distribution $\phi_0=0$, 
as well as forms $\phi_i$ such that the form $\Omega_V$ could be written as 
$\Omega_V = \sum_i\phi_0\wedge\phi_i\otimes e^i$. Also, the curve was 
either transverse to the $\phi_0$ distribution, or simply tangent  to it at 
type I branch points. This, as well as the vanishing of components of 
$\Omega_V$ at the type III points tells us that there is an exact sequence
$$
N_{S_u}\rightarrow K_{S_u}\otimes V^*\rightarrow 
\oplus_{x\in \Delta_u} (K_{S_u}\otimes (V_x)_0)_{x}.
$$
On the level of $W$-invariant sections, we have seen that the map
$$H^0(S_u, N_{S_u})^W\rightarrow H^0(S_u,K_{S_u}\otimes V^*)^W$$
is an isomorphism. However the normal bundle of $S_u$ in $\SS$ is
trivial, with a trivial $W$-action. The evaluation map from invariant sections 
to the bundle is then surjective for the normal bundle at each point.
Translating this over to $K_{S_u}\otimes V^*$ gives the desired result.
\EndProof

%**********************************************************
% Geometric conditions
%**********************************************************
\subsection{Geometric conditions}
\label{sec-geometric-conditions}

We now give a geometric condition that ensures Strong base point freeness.
We first describe  some of the sheaves we have been working with 
in terms of the quotient curve $C = S/W$, $S=S_u$. 
Throughout, we suppose that $S$
has simple branch points. Let $q:S\rightarrow C$ be the 
projection, and let $\omega_{S/C}$ denote the
relative dualising sheaf; we define
\begin{equation}
\label{eq-E}
E:=[(q_*\omega_{S/C})\otimes V^*]^W. 
\end{equation}
Duality for finite flat morphisms tells us that 
$E^*= [(q_*\StructureSheaf{S})\otimes V]^W$. 
Note that decomposing $V$ into irreducibles 
induces a corresponding decomposition of $E$. For a trivial summand, 
$E$ is trivial.
One can identify $q^*[K_C\otimes E]$ with the kernel of (\ref{4.6}):

\begin{new-lemma}
There is an exact sequence
$$ 0\rightarrow q^*[K_C\otimes E]\rightarrow 
K_{S}\otimes V^*\rightarrow \oplus_{x\in \Delta} (K_{S}\otimes
 (V_x)_0^*)_x\rightarrow 0.$$
\end{new-lemma}

Thus, Strong base point freeness is  a statement about invariant 
sections of $q^*[K_C\otimes E]$.

 We note that if the quotient $C$ is a rational curve, the Prym variety
is the same for $\chi$ and for $\chi/\chi^W$. We will therefore make the 

\begin{assumption}
\label{4.19}
If $C$ is rational, the representation $V$ has no trivial summand. 
\end{assumption}

The failure of Strong base point freeness has strong implications 
about the structure of $E$:

\begin{new-lemma}
\label{4.20}
Let $p$ be a point in the support of the co-kernel
sheaf $Q$ in the sequence 
\begin{equation}
\label{4.21}
0 \rightarrow H^0(K_C\otimes E)\otimes {\cal O}_C\rightarrow
K_C\otimes E\rightarrow Q\rightarrow 0
\end{equation}
Denote by $t_p $ the dimension of the vector space $Q/Q(-p)$. 
If $q:S\rightarrow C$ is a ramified cover, there are  short exact
 sequences
\begin{equation}
\label{4.22}
0\rightarrow [Q/Q(-p)]^*\otimes {\cal O}_C(-p)\rightarrow 
E^*\rightarrow F\rightarrow 0
\end{equation}
and dually,
$$0 \rightarrow K_C\otimes F^* \rightarrow K_C\otimes E \rightarrow 
K_C(p)\otimes [Q/Q(-p)]\rightarrow 0,$$
where $F$ is a vector bundle of rank $[v-t_p]$, and degree $(2t_p-\tilde b)/2$,
with, if $\Delta$ is the ramification locus of $q$ in $S$,
$$\tilde b = {2\over |W|} \sum_{x\in \Delta} \dim((V_x)_1)).$$
\end{new-lemma}

{\sc Proof}: 
 From the definition of $t_p$, one has that
$$h^0(C, K_C\otimes E(-p)) = h^0(C, K_C\otimes E) +t_p - v.$$
Riemann-Roch and Serre duality then tell us that:
$$h^0(C, E^*(p)) = h^0(C,  E^*) +t_p.$$
Now, for a trivial summand $V_i$, we have that the corresponding
$E_i$ is trivial, and so $E_i$ does not contribute to the quotient 
$Q$, unless $C$ is rational. This is precisely the case we have 
excluded, and so  we can replace $V$ by $V/(V^W)$ and so suppose that $V^W=0$.
We then have that $H^0(C,E^*) = H^0(C,  
[(q_*\StructureSheaf{S})\otimes V]^W) = V^W = 0$, so then 
$H^1(C,K_C\otimes E)=0$. 
Extend (\ref{4.21}) to a diagram 
$$
\matrix{ 0& \rightarrow& H^0(K_C\otimes E)\otimes {\cal O}_C(-p) &
\rightarrow&
K_C\otimes E(-p)&\rightarrow& Q(-p)&\rightarrow &0\cr
& & \downarrow& &\downarrow& & \downarrow\cr
0& \rightarrow& H^0(K_C\otimes E)\otimes {\cal O}_C&\rightarrow&
K_C\otimes E&\rightarrow& Q&\rightarrow &0\cr
& & & &\downarrow& & \downarrow\cr
& & & &(K_C\otimes E)|_p&\rightarrow&[Q/Q(-p)]} 
$$
We find that  $(K_C\otimes E)|_p$ maps surjectively both to $[Q/Q(-p)]$ and to
$H^1(C,K_C\otimes E(-p))$ and both maps have the same kernel. 
We get an isomorphism 
between $H^1(C,K_C\otimes E(-p))$ and $[Q/Q(-p)]$. Dually, 
$$H^0(C, E^*(p)) \simeq [Q/Q(-p)]^*$$
 with both spaces injecting into $(E^*(p) )|_p$. The evaluation homomorphism 
$$H^0(C, E^*(p)) \otimes {\cal O}_C\rightarrow E^*(p)$$
is then injective as a  sheaf map. If this map fails to be 
injective as a bundle map at some points, we can take a rank one subsheaf, 
twist it, and obtain a sub-line-bundle of $E^*(p)$ of positive degree.   
This would contradict  the following lemma:

\begin{new-lemma}
\label{lemma-non-existence-of-positive-line-sub-bundles}
\begin{enumerate}
\item 
Assume that $V^W=(0)$. If $q$ is ramified then every line-sub-bundle of 
the vector bundle $[q_{*}{\cal O}({S})\otimes V]^W$
has negative degree. 
If $q$ is unramified then every such line-sub-bundle
of degree zero is a non-trivial line bundle which is a torsion point of 
$\J_{C}$ of order which divides the order of $W$. 
\item 
There is a canonical exact sequence
$$
0 \rightarrow 
q^*[q_*{\cal O}(S)\otimes V]^W 
\rightarrow
{\cal O}(S)\otimes V
\rightarrow
\oplus_{x\in \Delta}(V_x)_0
\rightarrow 0. 
$$
\end{enumerate}

\noindent
In particular, the degree of $q^*[q_*{\cal O}(S)\otimes V]^W $ is 
$- \tilde b|W|/2$. 
Consequently, the degree of $[q_*{\cal O}(S)\otimes V]^W $ is 
$-\tilde{b} / 2$. 
\end{new-lemma}

\noindent
{\sc Proof:}
We first remark that the natural homomorphism
\begin{equation}
\label{4.26}
q^*\left([q_{*}{\cal O}(S)\otimes V]^W\right)
\buildrel {q^{\natural}} \over \longrightarrow
{\cal O}(S)\otimes V
\end{equation}
is injective and surjective away from the ramification locus $\Delta $. 

Let $L$ be a line sub-bundle of $[q_{*}{\cal O}(S)\otimes V]^W$.
Then $q^*(L)$ is a subsheaf of 
${\cal O}(S)\otimes V$ and hence has degree $\leq 0$. Moreover,
the equality $\deg(q^*(L))=0$ holds only if $q^*(L)$ is the trivial 
line-bundle. If $q$ is a branched covering, the homomorphism 
$q^*:\J_{C}\hookrightarrow \J_{S}$
is injective. 
We conclude that $\deg(L)$ is $\leq 0$ and 
equality $\deg(L)=0$ holds only if $L$ is the trivial line-bundle.
As  $H^0([q_{*}{\cal O}(S)\otimes V]^W)$ vanishes
so does $H^0(L)$. Hence, $L$ has negative degree. 
If $q$ is unramified and $L$ has degree zero we saw that $L$ must be in the
kernel of $q^*:\J_{C}\hookrightarrow \J_{S}$. 
In particular, $L$ has finite order, say $n$. 
The line bundle  $L$ determine
a cyclic $n$-sheeted cover $a:\tilde{C}\rightarrow C$ 
which is minimal in the following sense:  If $q':S' \rightarrow C$ is an
unramified cover and the pullback 
$q'^*(L)$ is trivial then $q'$ factors through $a:\tilde{C}\rightarrow C$. 
In particular, $n$ divides the order of $W$.

Recall that all the ramification points of $q$ are simple. 
At each ramification point $x\in S$, the image of 
$q^*\left([q_{*}{\cal O}(S)\otimes V]^W\right)$ 
in ${\cal O}(S)\otimes V$ must be in $(V_x)_0$.
 It follows that the cokernel
of $q^{\natural}$ in 
(\ref{4.26}) 
has length $= \tilde b|W|/2$. 
\EndProof

\medskip
This lemma, and the exact sequence (\ref{4.22}), 
tells us what the degree of
the cokernel $F$ must be, finishing the proof of  lemma
\ref{4.20}.  
\EndProof

\medskip
We can now prove  a geometric criterion for Strong base point
freeness. We will make the following assumption:

\begin{assumption}
\label{assumption-compatibility}
(Compatibility Assumption)
For all ramification points $x$,
the stabiliser $G_x$ acts on $V$ with weights $(1,0,..,0)$, so that the weight
one space is one dimensional. 
\end{assumption}

This is the case, for example, for the Hitchin systems, 
where the group is the Weyl
group $W$ of a reductive group, and the representation is the standard one
on the Cartan algebra. The stabilisers are generated by the reflection
in the root planes. The compatibility assumption rules out the 
possibility of ramification points of type II. It tells us that
$$b=\tilde b,$$
where $b$ is the cardinality of the image of the branch locus in $C$.

Let $\alpha(W,V)$ denote the number of different orbits
in $\bbp V$ of the points $\bbp((V_x)_1)$, and let $\lambda(W,V)$ 
be the order of the longest orbit. For example, for the Hitchin systems, 
$\alpha(W,V)$ counts the number of orbits of the root vectors
and is either 1 or 2, depending on the group. In turn, $\lambda(W,V)$
satisfies the inequality:
$$[\dim(\LieAlg{g}) - \dim(\LieAlg{h})]/4\le \lambda(W,V) \le
[\dim(\LieAlg{g}) - \dim(\LieAlg{h})]/2.$$

\begin{thm}
\label{thm-SBPF-holds-if-many-branch-pts-exist}
The Strong base point freeness condition holds if

\begin{trivlist}
\item [1.] The map $q$ has simple branch points, and the compatibility 
assumption \ref{assumption-compatibility} holds. 

\item [2.] Assumption \ref{4.19} holds.

\item [3.] The number $b$ of branch points of $q:S\rightarrow C$ satisfies:
\begin{equation}
\label{4.29}
b> 2 \alpha (W, V) \lambda (W,V).
\end{equation}
\end{trivlist}
\end{thm}

{\sc Proof:} 
Assume that the Strong base point freeness assumption does not hold.
Lemma \ref{4.20}
implies that there exists a surjective homomorphism of vector bundles

\begin{equation}
\label{4.30}
E\otimes K_C  \rightarrow K_C (p).
\end{equation}
The pullback of $E\otimes K_C$ is isomorphic to 
$ker\left[V^*\otimes K_S 
\rightarrow 
\oplus_{x\in \Delta_{q }}(V_x^*)_0\right].$
Pulling (\ref{4.30}) back 
to $S $ we get a surjective homomorphism 
$$
ker\left[V^*\otimes K_S  \rightarrow 
\oplus_{x\in \Delta_{q }}(V_x^*)_0\right] \rightarrow q ^*(K_C (p)). 
$$
Let $R$ be the finite subset 
$$R:=\{(V_x^*)_1 \mid \ \ x \in \Delta_{q } \}$$
of $\bbp({V}^*)$ (i.e., of hyperplanes in $\bbp({V})$). 
Choose a line $\ell \subset V^*$ lying
in the $W$-orbit in $R$ which arise most 
frequently from the ramification points. 
Denote the cardinality of this orbit by $\lambda(\ell)$. 
Let $L'$ be the
line sub-bundle $\ell\otimes K_S  \subset V^*\otimes K_S $ and 
$L$ the line bundle $L'(-D)$ where $D$ is the effective divisor 
$$D := \sum \{x\in \Delta_{q }\mid \ \ \ell \neq (V_x^*)_1\}.
$$
Then $L$ is a subsheaf of the pullback of $q ^*[E\otimes K_C ]$. 
Hence, there exists  a non-trivial homomorphism from
$L$ to $q ^*(K_C(p))$ (perhaps after replacing $\ell$ by another 
line in the same $W$-orbit).
On the other hand, a simple computation  shows that 
$\deg(L)$ is larger than $\deg(q^*(K_C(p)))$. 
This will give the desired contradiction.
 
We proceed to compute the degree of $L$. 
If $H_i = (V_x^*)_0$ then the stabilizer $W_{H_i}$ has order
${{|W|}\over {\lambda(\ell)}}$ and it contains $W_x$. 
Hence, there are 
${{|W|}\over{2\lambda(\ell)}}$ ramification points
$x'$ in the fiber over $q(x)$ satisfying $\ell = (V_{x'}^*)_1$. 
By our choice of $\ell$, 
the number of branch points of this type is 
$\geq {{b}\over{\alpha(W,V)}}$ where $b$ is the total number of branch points.
By definition, we have the inequality $\lambda(\ell)\leq \lambda(W,V)$.  
Thus, the degree of $L$ satisfies the inequality
\begin{eqnarray*}
\deg(L) & \geq & (2g_{S }-2) - \deg(\Delta_{q }) +
{{b}\over{\alpha(W,V)}}\cdot{{|W|}\over{2\lambda(W,V)}}
%\\
%& = &  \left[ (2g_{C}-2)\cdot |W| + \deg(\Delta_{q })\right]
%- \left[1-{{1}\over{\alpha(W,V)\cdot\lambda(W,V)}}\right]\cdot
%\deg(\Delta_{q })
%\\
%& = & (2g_{C}-2)\cdot |W|
%+ {{\deg(\Delta_{q })}\over{\alpha(W,V)\cdot\lambda(W,V)}}
\\
& = &  \left[
(2g_{C }-2)
+ {{b}\over{2\alpha(W,V)\cdot\lambda(W,V)}}
\right]\cdot |W|. 
\end{eqnarray*}
On the other hand, we have
$$
\deg\left(q ^*(\omega_{C}(p))\right) = (2g_{C}-1)\cdot |W|. 
$$
Taking the difference we get
$$
{{\deg(L) - \deg(q ^*(\omega_{C}(p)))}\over{|W|} }
\geq \left[
{{b}\over{2\alpha(W,V)\cdot\lambda(W,V)}}
-1\right]. 
$$
This difference is positive if the inequality 
(\ref{4.29}) holds. 
\EndProof

%**************************************************************
% Test of efficacy of our geometric criterion:
%**************************************************************
\medskip
{\sc Example:}\ As a test of the efficacy of our geometric criterion,
we can show that the strong base point freeness 
Assumption \ref{assumption-strong-generation} 
holds in most cases for the Hitchin system
and its generalizations. For these cases, one chooses a reductive Lie algebra
$\LieAlg{g}$ and an effective divisor $D$ of degree $d$ on a fixed curve $C$. 
The group $W$ 
is then the Weyl group and   $V$ is the Cartan subalgebra of ${ \Got g}$. 
The curves $S$ embed equivariantly into $T^*C(D)\otimes V$. 
In these examples, Assumption \ref{assumption-strong-generation}
can be verified independently of theorem
\ref{thm-SBPF-holds-if-many-branch-pts-exist} because the vector bundle
$E\otimes T^*C$ is a direct sum of line bundles 
(see Lemma \ref{lemma-E-is-a-direct-sum-of-line-bundles}).
%Assumption \ref{assumption-strong-generation} holds if and only if
%each of these line bundles is generated by its global sections. 
Comparing both methods we will see that Theorem 
\ref{thm-SBPF-holds-if-many-branch-pts-exist} successfully establishes
Assumption \ref{assumption-strong-generation} in almost all cases. 

%If for example $W$ is the symmetric group ${\rm Sym}_n$, $n \geq 3$, 
%and $V$ its standard $n\!-\!1$ dimensional representation, then 
%we require more than $n(n\!-\!1)$ branch points. 
%More generally, 
If the rank of  $\LieAlg{g}$
is $\geq 2$, then we require that
\begin{eqnarray}
\label{eq-inequality-for-number-of-branch-points-lie-gp-case}
b & > & [\dim(\LieAlg{g}) - \dim(\LieAlg{h})] \ \ \ \ \ \ \
\mbox{if the group is simply laced and}
\\
\label{eq-inequality-for-number-of-branch-points-non-simply-laced-lie-gp-case}
b & > &  2\cdot[\dim(\LieAlg{g}) - \dim(\LieAlg{h})] \ \ \ 
\mbox{if the group is not simply laced.}
\end{eqnarray}

%In the Hitchin system of $G$-principal Higgs pairs the actual number of
%branch points is 
%$(2g_{\Sigma}-2)\cdot [\dim(\LieAlg{g}) - \rank(\LieAlg{g})]$. 
%Inequality (\ref{eq-inequality-for-number-of-branch-points-lie-gp-case}) 
%is satisfied in the Hitchin system case for $g_{\Sigma} \geq 2$. 
%In the $g_{\Sigma}=1$ case the spectral curve $S_u$ of the Hitchin system is 
%disconnected. The Hitchin system is empty in case $g_{\Sigma}=0$.
%In the non-simply-laced case, inequality 
%(\ref{eq-inequality-for-number-of-branch-points-non-simply-laced-lie-gp-case})
%is satisfied in the Hitchin system case for $g_{\Sigma} > 2$. 

In the example of integrable systems of reduced 
coadjoint orbits  (in section \ref{sec-coadjoint-orbits}) and 
in the generalized Hitchin integrable systems 
(see \cite{markman-spectral-curves}) the number of branch points 
$b$ is equal to 
$(2g_{C}-2+d)\cdot[\dim(\LieAlg{g}) - \dim(\LieAlg{h})]$ when
 $(2g_{C}-2+d)$ is positive. 
Theorem \ref{thm-SBPF-holds-if-many-branch-pts-exist}
implies Assumption \ref{assumption-strong-generation} 
in all cases except $(g_{C},d)=(0,3)$, 
$(g_{C},d)=(1,1)$, where equality holds in 
(\ref{eq-inequality-for-number-of-branch-points-lie-gp-case}) 
and the reversed inequality in
(\ref{eq-inequality-for-number-of-branch-points-non-simply-laced-lie-gp-case}).
Note that Assumption \ref{assumption-strong-generation} fails in both cases 
for the group $SL(n)$ 
(see Lemma \ref{lemma-E-is-a-direct-sum-of-line-bundles}).
In the non-simply-laced case  equality holds in 
(\ref{eq-inequality-for-number-of-branch-points-non-simply-laced-lie-gp-case})
also when $(g_{C},d)$ equals $(2,0)$, $(0,4)$ and $(1,2)$. Now let us 
check Assumption \ref{assumption-strong-generation} directly:

%****************************************
% Lemma: E is a direct sum in the generalized Hitchin systems
%****************************************
\begin{new-lemma}
\label{lemma-E-is-a-direct-sum-of-line-bundles} 
Assume that the algebra of invariant polynomials $\ComplexNumbers[V]^W$ 
is free. Let $(d_1\leq d_2 \leq \dots \leq d_r)$ be the degrees of the
$W$-invariant homogenous generators. Fix a line-bundle $L$ on $C$ and 
assume that $S$ admits a
$W$-equivariant embedding into the vector bundle
$X:=L \otimes V^*$. Then 
\begin{enumerate}
\item
\label{lemma-item-E-is-a-direct-sum-of-line-bundles}
the vector bundle $E$ in
(\ref{eq-E}) is isomorphic to
\begin{equation}
\label{eq-E-is-a-direct-sum-of-line-bundles}
E \cong \oplus_{i=1}^r L^{\otimes [d_i-1]}.
\end{equation}
\item
\label{lemma-item-strong-generation-assumption-for-spectral-covers}
If $L$ equals $T^*C(D)$ for some effective divisor $D\geq 0$, then
$E\otimes T^*C$ is isomorphic to
\[
E\otimes T^*C \cong 
\oplus_{i=1}^r (T^*{C})^{\otimes d_i}([d_i-1]\cdot D)
\]
and is generated by global sections precisely in the following cases:
\begin{enumerate}
\item
$g_{C} \geq 2$.
\item
$g_{C} = 0$, $\deg(D)\geq 4$ and $d_1\geq 2$,
\item
$g_{C}=0$, $\deg(D)=3$ and $d_1\geq 3$,
\item
$g_{C}=1$ and $\deg(D)\geq 2$,
\item
$g_{C}=1$, $\deg(D)=1$, and $d_1\geq 3$.
\end{enumerate}
\end{enumerate}
\end{new-lemma}

\noindent
{\sc Proof:}
\ref{lemma-item-E-is-a-direct-sum-of-line-bundles})
A choice of generators $\{f_1,\dots,f_r\}$ 
of $\ComplexNumbers[V]^W$ determines an isomorphism 
of vector bundles over $C$
\[
(L\otimes V^*)/W \IsomRightArrow  \oplus_{i=1}^r L^{\otimes d_i}.
\]
A relative version of 
Lemma \ref{lemma-invariant-polynomials-split} implies that we have an 
isomorphism of vector bundles over the total space of 
$\oplus_{i=1}^r L^{\otimes d_i}$
\begin{equation}
\label{eq-relative-isomorphism-of-invariant-cotangent-bundles}
q_*(\pi^*[L^{-1}\otimes V])^W \ \LongIsomRightArrow \ 
\bar{\pi}^* [\oplus_{i=1}^r L^{\otimes -d_i}]. 
\end{equation}
Thinking off $C$ as a section of $\oplus_{i=1}^r L^{\otimes d_i}$ and
restricting the isomorphism 
(\ref{eq-relative-isomorphism-of-invariant-cotangent-bundles}) to $C$ 
we get the isomorphism
\[
q_{*}(q^*(L^{-1})\otimes V)^W \IsomRightArrow 
\oplus_{i=1}^{r} L^{\otimes -d_i}.
\]
Tensoring both sides by $L$ we get
\[
q_{*}(\StructureSheaf{S}\otimes V)^W \IsomRightArrow 
\oplus_{i=1}^{r} L^{\otimes [1-d_i]}.
\]
The isomorphism $E^*\cong q_{*}(\StructureSheaf{S}\otimes V)^W$
implies the isomorphism
(\ref{eq-E-is-a-direct-sum-of-line-bundles}).

\medskip
\noindent
\ref{lemma-item-strong-generation-assumption-for-spectral-covers})
Follows easily from Part \ref{lemma-item-E-is-a-direct-sum-of-line-bundles}. 
\EndProof

\begin{new-lemma}
\label{lemma-invariant-polynomials-split}
Let $(V,W)$ be as in Lemma \ref{lemma-E-is-a-direct-sum-of-line-bundles},
$q:V\rightarrow V/W$ the quotient map. 
Then the codifferential $d^*q$ induces a canonical isomorphism of trivial 
vector bundles over $V/W$
\begin{equation}
\label{eq-codifferential-induces-isomorphism-of-cotangent-bundles}
q_*(d^*q): T^*(V/W) \ \LongIsomRightArrow \ [q_* T^*V]^W.
\end{equation}
%\begin{equation}
%\label{eq-invariant-summand-is-the-trivial-vb-of-invariants}
%q_*(d^*q): (V/W)
%\otimes_{\ComplexNumbers}\StructureSheaf{V/W} \ \ 
%\LongIsomRightArrow \ \ 
%[q_*\StructureSheaf{V}\otimes V^*]^W.
%\end{equation}
\end{new-lemma}

\noindent
{\sc Proof:}
That (\ref{eq-codifferential-induces-isomorphism-of-cotangent-bundles})
is an isomorphism is a general fact which holds 
whenever the quotient $Y:=X/W$ of a smooth variety $X$ by 
a finite group $W$ is smooth. 
%Let $\delta$ be a 
%codimension $1$ irreducible component 
%of the ramification divisor of the quotient map
%$q:X\rightarrow Y$. Assume that $W$ embedds in the automorphism group of 
%$X$ (the action is faithful).
%Endow $\delta$ with the reduced induced scheme structure.
%The stabilizing subgroup $W_\delta$ of the generic smooth point of
%$\delta$ is a cyclic group. The point is that the $W$ action can be linearized
%locally around such a point $p \in \delta$ and thus $W_\delta$ embedds 
%via the character $\rho_\delta$ of its action on the normal 
%fiber $[N_{\delta/X}]_p$ and acts trivially on the tangent hyperplane 
%$T_p\delta$. The isomorphism 
%(\ref{eq-codifferential-induces-isomorphism-of-cotangent-bundles}), 
%locally around the fiber $q^{-1}(q(p))$, is equivalent to 
%an isomorphism locally around $p$, replacing $W$ by $W_\delta$. 
%Proving it is now a simple calculation in local coordinates (which 
%is identical to that in the one-dimensional case). 
%We conclude that 
%(\ref{eq-codifferential-induces-isomorphism-of-cotangent-bundles}) 
%is an isomorphism over a Zariski open subset whose complement has
%codimension $\geq 2$. Hence it is a global isomorphism.
\EndProof

%********************************************************************
% From Hilbert schemes to integrable systems
%********************************************************************
\section{From Hilbert schemes to integrable systems}
\label{sec-from-hilbert-schemes-to-acihs}

We prove in this section the converse to Theorem 
\ref{thm-quotient-of-universal-curve}. 
We start with a family of curves on a variety $X$ 
which has the geometric properties expected from the quotient 
of a family of curves $\SS$ by a null foliation. We formulate this
setup in section \ref{sec-lagrangian-fibrations-over-a-curve} where 
we define the notion of a $(W,V)$-lagrangian fibration 
$(\pi:X\rightarrow\Sigma,W,V,\Omega_V)$ over a curve
(Definition \ref{def-V-Lagrangian-fibration}). 
We have made an effort to work with a sufficiently general
definition which includes the important examples of 
integrable systems coming from meromorphic Higgs pairs and reduced
coadjoint orbits of loop groups. In particular, we allow 
ramification points whose stabilizers have arbitrary order. 

In section \ref{sec-smoothness-of-hilbert-scheme}
we study the deformation theory of smooth $W$-Galois covers of $\Sigma$  
which are embedded in $X$. We prove that the 
Hilbert scheme of such curves is smooth. 
The main result, 
Theorem \ref{thm-compatibility-of-2-forms}, is proven in 
Section \ref{sec-symplectic-structure-on-relative-prym}. 
We construct a completely
integrable Hamiltonian system supported on the family of generalized pryms
of all smooth $W$-Galois covers of $\Sigma$  
which are embedded in $X$. 

%***************************************************************************
% $(W,V)$-Lagrangian fibrations over a curve
%***************************************************************************
\subsection{
$(W,V)$-Lagrangian fibrations over a curve
} 
\label{sec-lagrangian-fibrations-over-a-curve}

We describe in this section the fine structure of the 
quotient $X$ of the relative Galois cover $\SS \rightarrow U$ by the null
foliation. Especially delicate is the structure 
along the ``vertical'' ramification divisors, those which are leaves of the $\phi_0$-
foliation. 
We will see examples with vertical ramification divisors in section 
\ref{sec-the-lagrangian-fibration-of-a-coadjoint-orbit}. 

\medskip
Let $W$ be a finite group, $\chi$ a faithful integral representation
of rank $v$,
$V:=\chi_{\ComplexNumbers}$ its complexification, 
$X$ a smooth $n+1$ dimensional quasi-projective (or complex analytic) variety,
$W\hookrightarrow \Aut(X)$ a faithful action, $q:X\rightarrow Y$ the  
quotient,
$\Delta_q\subset X$ the ramification divisor,  
\[
\pi : X \rightarrow \Sigma
\]
a $W$-invariant morphism onto a smooth projective curve $\Sigma$,
and 
\[
\Omega_V \in H^0(X,[\Wedge{2}T^*X]\otimes V^*)
\]
a closed $W$-invariant $V^*$-valued $2$-form. 

%*****************************************************
% Assumption: regularity of $(W,V)$-Lagrangian fibration
%*****************************************************

\begin{assumption}
\label{assumption-regularity-of-W-V-lagrangian-fibration}
We require that the quotient $Y$ is {\em smooth} and the morphism
$\bar{\pi}: Y \rightarrow \Sigma$ is submersive. We further assume that 
the fibers of $\pi : X \rightarrow \Sigma$ are smooth (but may have
multiplicities). 
\end{assumption}

We denote by 
\[
\Delta_q = \Delta_q[1] \supset \Delta_q[2] \supset \cdots \supset 
\Delta_q[k] \supset \Delta_q[k+1] \cdots 
\]
the stratification where $\Delta_q[k]$ is the locus of points contained in 
the intersection of $k$ components of $\Delta_q$. We further set
\[
X[0] \supset X[1] \supset \cdots X[k] \supset X[k+1] \cdots
\]
to be the stratification defined by $X[0]=X$ and 
$X[k]=\Delta_q[k]$ for $k\geq 1$. Let 
\[
\Delta_q \ = \ \Delta_q^{reg} \ + \ \Delta_q^{\infty}
\]
be the decomposition into non-fiber ($\Delta_q^{reg}$) and   fiber 
($\Delta_q^{\infty}$) components with respect to $\pi$. 
Denote by $red(\Delta_q^{\infty})$ the reduced induced 
scheme structure of $\Delta_q^\infty$.
Given a component $\delta$ of $red(\Delta_q^{\infty})$, 
we denote by $W_\delta$ 
its (necessarily cyclic) stabilizer subgroup in $W$ and by $m_\delta$ 
the order of $W_\delta$. 
We have the equality
\[
\Delta_q^{\infty} = 
\sum_{\delta\subset red(\Delta_q^{\infty})} (m_\delta-1) \delta. 
\]
Observe that the smoothness of both $X$ and the quotient $Y$ implies that 
each component $\delta$ of $\Delta_q$, endowed with the reduced induced 
scheme structure, is a smooth subvariety of $X$. The point is that 
$\delta$ is the fixed locus of the finite group $W_\delta$. 
We denote by $\rho_\delta$ the primitive character of $W_\delta$ 
corresponding to the fibers of $N_{\delta/X}$.

\begin{defi}
\label{def-V-Lagrangian-fibration}
The quadruple $(\pi:X\rightarrow\Sigma,W,V,\Omega_V)$
is said to be a $(W,V)$-Lagrangian fibration if the following conditions
are satisfied:
\begin{enumerate}
%*******
% Item
%*******
\item
\label{def-item-pi-is-isotropic}
The morphism $\pi$ is $\Omega_V$-isotropic.
%*******
% Item
%*******
\item
\label{def-item-surjectivity-of-contraction-away-from-infinity}
The contraction
\begin{equation}
\label{eq-contraction-in-def-of-V-lagrangian}
\Contract\Omega_V \ : \ \Wedge{2}TX \rightarrow V^*\otimes \StructureSheaf{X}
\end{equation}
is surjective over $X \setminus \Delta_q^{\infty}$. 
%*******
% Item
%*******
\item
\label{def-item-minimality-of-degeneracy-of-2-form}
(Minimal degeneracy of the $2$-form $\Omega_V$ along
$\Delta_q^\infty$) The image subsheaf 
\begin{equation}
\label{eq-V-prime-is-the-image-of-contraction}
(V')^*\subset V^*\otimes \StructureSheaf{X}
\end{equation} 
of the contraction
homomorphism (\ref{eq-contraction-in-def-of-V-lagrangian})
satisfies the equality
\begin{equation}
\label{eq-characterization-of-V-prime}
(V')^*
\ \  = \ \ 
\left\{
q^*\left(
q_*\left[
\StructureSheaf{X}(-\sum\delta_i)\otimes V^*
\right]^W
\right) 
\right\}\left(\sum\delta_i\right)
\end{equation}
locally along $\Delta_q^\infty$.
Above, $\delta_i$ are the components of the reduced induced scheme
structure of $\Delta_q^{\infty}$.
\end{enumerate}
\end{defi}

\begin{rem}
\label{rem-after-def-of-W-V-lagrangian-fibration}
{\rm
\begin{enumerate}
\item
\label{rem-item-minimal-degeneracy-is-satisfied-by-quotient-of-prym-acihs}
Let $\Prym_{\SS/U}\rightarrow U$
be a rank 2 Prym integrable system. 
Assume that the Strong Base Point Freeness condition holds
(where we use the generalization of Assumption
\ref{assumption-strong-generation} 
to the case of general ramification points, i.e., we assume 
that the vector bundle $q_{u_*}(T^*S_u\otimes V^*)^W$ 
is generated by global sections). 
Let $(\pi:X\rightarrow \Sigma,W,V,\Omega_V)$ be the quotient 
of $\SS$ by the null foliation. Assume further that $Y:=X/W$ is smooth. 
Then the minimal degeneracy condition 
(part \ref{def-item-minimality-of-degeneracy-of-2-form} of 
Definition \ref{def-V-Lagrangian-fibration}) 
is satisfied.
The proof goes as follows. 
Proposition \ref{prop-isomorphism-of-vector-spaces-of-global-sections} 
together with Assumption \ref{assumption-strong-generation} 
imply that the contraction with 
$\Omega_V$ induces an isomorphism of vector bundles 
$f^W: (q_{u_*}N_{S_u/X})^W\rightarrow q_{u_*}(T^*S_u\otimes V^*)^W$ 
over the quotient curves $C_u:=S_u/W \subset Y$. 
But the surjectivity of $f^W$ 
is equivalent to the minimal degeneracy condition 
(\ref{eq-characterization-of-V-prime}) 
(see the proof of the first part of
Lemma \ref{lemma-identification-of-tangent-to-hilbert-scheme}). 
\item
It follows from the smoothness of the quotient $Y$
and condition \ref{def-item-surjectivity-of-contraction-away-from-infinity}
of the Definition 
that each component $\delta$ of $\Delta_q^{reg}$ is reduced and its stabilizer 
subgroup $W_\delta$ in $W$ is isomorphic to $\Integers/2\Integers$. 
Moreover, the multiplicity of the trivial character of $W_\delta$
in $V$ is precisely $n-1$. 
The point is that $\Omega_V$ induces an isomorphism 
\begin{equation}
\label{eq-integrable-trivialization-of-vertical-tangent-bundle-of-pi}
T_\pi \otimes \pi^*T\Sigma  \cong V^*
\end{equation}
over $[X\setminus \Delta_q^\infty]$. 
Hence, the multiplicity of the trivial character 
in $T_\pi$ is  equal to its multiplicity in $V$. 
As $T\delta$ and $T_\pi$ are transversal and $TX/T_\pi$ is the trivial 
character, $T_\pi$ contains the primitive character 
$\rho_\delta$ of $W_\delta$. 
Hence, so does $V$. 
On the other hand, $V$ is an integral representation.
Hence, if the stabilizer $W_\delta$ has order $\geq 3$ than the 
multiplicity of the trivial character in $V$ is $\leq n-2$. 
It would follow that the multiplicity of the trivial character 
of $W_\delta$ in $\restricted{TX}{\delta}$ is $\leq n-1$. This contradicts the
smoothness of $Y$ (the triviality of the $W_\delta$-representation
$T\delta$). 
\item
A fiber of $\pi : [X\setminus \Delta_q^\infty] \rightarrow \Sigma$ over
a point $a\in \Sigma$ 
has a canonical affine structure modeled over the vector space
$T^*_a\Sigma\otimes V^*$ via the integrable trivialization
(\ref{eq-integrable-trivialization-of-vertical-tangent-bundle-of-pi}).
In particular, smooth projective fibers are abelian varieties. 
\item
The character $\rho_{\delta}$ of the stabilizer sub-group $W_\delta$
of a component $\delta$ of $red(\Delta_q^\infty)$ may have zero multiplicity 
in $V$. 
%(see Section \ref{sec-structure-of-vertical-ramification-divisor}).
\end{enumerate}
}
\end{rem}

\medskip
The following lemma explains
why part \ref{def-item-minimality-of-degeneracy-of-2-form}
of Definition \ref{def-V-Lagrangian-fibration} is a minimal degeneracy 
condition for the $2$-form $\Omega_V$: 

\begin{new-lemma}
\label{lemma-characterization-of-V-prime-as-a-maximal-subsheaf}
{\rm
The right hand side $(V'')^*$ of
(\ref{eq-characterization-of-V-prime})
is characterized as the maximal 
(locally free) subsheaf
of $V^*\otimes \StructureSheaf{X}$ satisfying
\begin{enumerate}
\item
$(V'')^* = V^*\otimes \StructureSheaf{X}$ over $X\setminus \Delta_q^{\infty}$
and
\item
The restriction of $(V'')^*$ to each component 
${\displaystyle \delta\subset \left[red(\Delta_q^{\infty})\setminus 
\Delta_q^{\infty}[2]\right]}$
is locally isomorphic to 
$[N_{\delta/X}]^{\oplus n}$ as an 
$\StructureSheaf{\delta}[W_\delta]$-module. 
In other words, each fiber of the vector bundle $(V'')^*$ over
$\delta$ is a direct sum of $n$ copies of the same primitive character 
$\rho_\delta$ 
of the cyclic stabilizer subgroup $W_\delta$.
\end{enumerate}
}
\end{new-lemma}

\noindent
The proof follows immediately from the definition of $(V'')^*$. On the other
 hand, one
can check that, if $(\pi:X\rightarrow\Sigma,\Omega_V)$ 
satisfies all the properties of 
a $(W,V)$-Lagrangian fibration 
over $\Sigma$ except possibly property 
\ref{def-item-minimality-of-degeneracy-of-2-form}
of Definition \ref{def-V-Lagrangian-fibration}, then 
$(V')^*$  satisfies (1) and (2) in the lemma. The minimal degeneracy
then says that $(V')^*$ is maximal.

%***************************************************************************
% Smoothness of the Hilbert scheme of curves
%***************************************************************************
\subsection{
Smoothness of the Hilbert scheme of curves
} 
\label{sec-smoothness-of-hilbert-scheme}

Let $U_Y$ be an irreducible Zariski open subset of a Hilbert scheme of 
curves on $Y$ which satisfy 

\begin{enumerate}
\item
$C_u \subset Y$ is a section of $\bar{\pi}: Y \rightarrow \Sigma$, i.e.,
$\bar{\pi}_u : C_u \rightarrow \Sigma$ is an isomorphism.
\item
\label{item-transversality-of-intersection-of-section-u-with-branch-divisor}
(Transversality of the intersection of $C_u$ with the
branch locus) \ 
Let $\pi_u : S_u \rightarrow \Sigma$ be the $W$-Galois cover 
obtained as the inverse image $S_u := q^{-1}(C_u) \subset X$. 
Then $S_u$ is a smooth curve. 
\end{enumerate}

\noindent
%Note that assumption 
%\ref{item-transversality-of-intersection-of-section-u-with-branch-divisor}
%implies that 1) The stabilizer $W_x$ of a point $x$ in $S_u$ 
%is equal to
%$\pi_u : S_u \rightarrow \Sigma$

Denote by $U_X$ the ($W$-invariant) sub-scheme of the Hilbert scheme
of curves on $X$ which are inverse images of curves in  $U_Y$. 
Clearly $U_X$ is isomorphic to $U_Y$. We will abuse notation and
denote both by $U$.
Let $p : \SS \rightarrow U$ be the universal curve, 
$\bar{p} : \C \rightarrow U$ its $W$-quotient, 
$N_\SS := N_{\SS/[X\!\times\!U]}$  and $N_{\C}:= N_{\C/[X\!\times\!U]}$ 
their relative normal bundles. Note that $\C$ is isomorphic 
to $U\times \Sigma$.
Let $\ell : \SS \rightarrow X$ be the natural morphism. 

\begin{prop}
\label{prop-smoothness-of-hilbert-scheme}
\begin{enumerate}
\item
\label{prop-item-smoothness}
$U$ is smooth.
\item
\label{prop-item-identification-of-tangent-to-hilbert-scheme}
There is canonical isomorphisms of sheaves on $U$
\begin{equation}
\label{eq-identification-of-tangent-to-hilbert-scheme-as-an-hodge-bundle}
TU \ \ \IsomRightArrow 
p_*(N_\SS)^W \ \ \LongIsomRightArrowOf{f^W} 
p_*(\omega_{\SS/U}\otimes_{\ComplexNumbers}V^*)^W. 
\end{equation}
\end{enumerate}
\end{prop}

\noindent
{\sc Proof:}
The main task is to construct the isomorphism $f^W$ in 
(\ref{eq-identification-of-tangent-to-hilbert-scheme-as-an-hodge-bundle}).
This is achieved in Lemma 
\ref{lemma-identification-of-tangent-to-hilbert-scheme}. 
The isomorphism $f^W$ identifies the Zariski tangent sheaf 
of the Hilbert scheme as a
$W$-invariant sub-bundle of a Hodge-bundle. 
The smoothness in 
Part \ref{prop-item-smoothness}, as well as the isomorphism 
$TU \IsomRightArrow 
p_*(N_\SS)^W$ then follows from the infinitesimal 
$T^1$-lifting property \cite{ziv-ran-T-1-lifting,kawamata}
(see also \cite{donagi-markman-cime} Ch 8 Section 8.2). 
\EndProof

\medskip
The $2$-form $\Omega_V$ induces a sheaf homomorphism
(see (\ref{eq-commutative-diagram-defining-f}))
\[
f : q_{*}(N_{\SS}) \rightarrow 
q_{*}\left(T^*_{\SS/U}\otimes_{\ComplexNumbers} V^*\right).
\]
In general, $f$ is not an isomorphism. However, $f$ induces an
isomorphism on the invariant sub-bundles.

\begin{new-lemma}
\label{lemma-identification-of-tangent-to-hilbert-scheme}
The $2$-form $\Omega_V$ induces an isomorphism of sheaves on
$U\times \Sigma$. 

\begin{equation}
\label{eq-sheaf-isomorphism-leading-to-identification-of-tangent-to-hilb}
f^W : q_*(N_\SS)^W \ \ \LongIsomRightArrow 
q_*(T^*_{\SS/U}\otimes_{\ComplexNumbers}V^*)^W. 
\end{equation}
\end{new-lemma}

\noindent
{\sc Proof:}
Since each $S_u$ is $1$-dimensional, it is $\Omega_V$-isotropic
and we get the commutative diagram of vector bundles on $\SS$

\begin{equation}
\label{eq-commutative-diagram-defining-f}
{
\divide\dgARROWLENGTH by 4
\begin{diagram}
\node{0}
\arrow{e}
\node{T_{\SS/U}} 
\arrow{e}
\arrow{s}
\node{\ell^*(T_X)}
\arrow{e}
\arrow{s,r}{\Omega_V}
\node{N_{\SS}}
\arrow{e}
\arrow{s,r}{f}
\node{0}
\\
\node{0}
\arrow{e}
\node{N^*_{\SS}\otimes V^*}
\arrow{e}
\node{\ell^*(T^*X) \otimes V^*}
\arrow{e}
\node{T^*_{\SS/U}\otimes V^*}
\arrow{e}
\node{0.}
\end{diagram}
}
\end{equation}

\noindent
Pushing forward to $U \times \Sigma$ and taking $W$-invariant sub-bundles we 
get the homomorphism $f^W$ in 
(\ref{eq-sheaf-isomorphism-leading-to-identification-of-tangent-to-hilb}).
Away from the branch divisor $f^W$ is an isomorphism. Hence,
$f^W$ is an injective sheaf homomorphism. We prove surjectivity
over every $C_u$. 

\smallskip
\noindent
{\bf Surjectivity of $f^W$ over $\Delta_{q_u}^{\infty}$:} 
Over $\Delta_{q_u}^{\infty}$ the surjectivity is essentially 
the minimal degeneracy condition 
(part \ref{def-item-minimality-of-degeneracy-of-2-form}
of Definition \ref{def-V-Lagrangian-fibration}). 
At a point $x \in \Delta_{q_u}^{\infty}$ the curve $S_u$ is transversal to 
the divisor $red(\Delta_q^\infty)$ (otherwise $S_u$ would not be smooth). 
Denote by $y$ the point $q_u(x)$ in $C_u$. 
Let $\delta\subset red(\Delta_q^\infty)$ be the corresponding component
through $x$. Thus, $T{\delta}$ surjects onto $N_{S_u/X}$ at $x$. 
The fiber of 
$T_x{\delta}\wedge T_xS_u$ projects onto $N_{S_u/X}\otimes T_xS_u$. 
By the minimal degeneracy condition, it maps 
via $f$ onto the fiber of $(V')^*$. Taking $W$-invariant sections, 
we conclude that $[q_{u_*}(N_{S_u/X}\otimes TS_u)]^W$ maps onto
$[q_{u_*}((V')^*)]^W$. The latter is, by definition of $V'$, 
\begin{equation}
\label{eq-bla}
\left[q_{u_*}\left\{V^*(-red(\Delta_{q_u}^\infty))\right\}\right]^W.
\end{equation}

%Tensor both sides of 
%(\ref{eq-sheaf-isomorphism-leading-to-identification-of-tangent-to-hilb})
%by $T\Sigma(-y)$. 
%by the $W$-invariant line bundle $\pi^*(T\Sigma)$. 
Since $N_{S_u/X}$ is a $W$-invariant vector bundle 
(it is the pullback of $N_{C_u/Y}$), then we have the isomorphisms 
\[
[q_{u_*}N_{S_u/X}\otimes TS_u]^W
\cong
N_{C_u/Y}\otimes (q_{u_*}TS_u)^W
\cong
N_{C_u/Y}\otimes TC_u(-y)
\cong
(q_{u_*}N_{S_u/X})^W \otimes T\Sigma(-y). 
\]
Thus, locally around $y$, the restriction of the left hand side of 
$T\Sigma(-y)\otimes$\/
(\ref{eq-sheaf-isomorphism-leading-to-identification-of-tangent-to-hilb})
to $C_u$ is $[q_{u_*}(N_{S_u/X}\otimes  TS_u)]^W$. 
%or equivalently  
%$[q_{u_*}(T_x{\delta}\wedge T_xS_u)]^W$. 
On the other hand, the restriction of the right hand side of 
$T\Sigma(-y)\otimes$\/
(\ref{eq-sheaf-isomorphism-leading-to-identification-of-tangent-to-hilb})
is (locally around $y$) 
\begin{equation}
\label{eq-bla-bla}
\left[
q_{u_*}\left\{
V^*\otimes \StructureSheaf{S_u}((m_y-1)red(\Delta_{q_u}^\infty))
\right\}
\right]^W \otimes \StructureSheaf{C_u}(-y)
\cong
\left[
q_{u_*}\left\{
V^*(-red(\Delta_{q_u}^\infty))
\right\}
\right]^W.
\end{equation}
The surjectivity of $f^W$ around $y$ follows from the equality
of (\ref{eq-bla}) and (\ref{eq-bla-bla}).

\medskip
\noindent
{\bf Surjectivity of $f^W$ over $\Delta_{q_u}^{\rm reg}$:} 
This will follow
from  part \ref{def-item-surjectivity-of-contraction-away-from-infinity}
of Definition \ref{def-V-Lagrangian-fibration}. Indeed, we know 
that the quotiented curves are sections of $Y$, and that the 
points in $\Delta_{q_u}^{\rm reg}$ have stabiliser $\bbz/2\bbz$ (remark
\ref{rem-after-def-of-W-V-lagrangian-fibration}). It follows that the 
curves $S_u$, at their intersection with the divisor $\Delta_{q}^{\rm reg}$,
are simply tangent to the fibers of the projection $\pi: X\rightarrow \Sigma$.
The isomorphism 
(\ref{eq-integrable-trivialization-of-vertical-tangent-bundle-of-pi})
can be written as the top row of a diagram 

\begin{equation}
\label{fibers-and curves}
\matrix {
T_\pi                        & \cong       & \pi^*T^*\Sigma\otimes   V^* \cr
p \downarrow  \hspace{1ex}   &             & \hspace{3ex} \downarrow 
                                                    d^*q_u\otimes 1 \cr
N_S                          & \stackrel{f}{\rightarrow} 
                                           & \omega_S\otimes V^*}
\end{equation}
The horizontal maps are the contractions with $\Omega_V$, and the
vertical ones are the natural projections.
One can choose suitable bases so that the right vertical map $d^*q_u\otimes 1$
is given by $z\cdot \bbi$, where $z$ is a coordinate on the curve vanishing at 
the tangency point. The left vertical map $p$ 
has co-rank $1$ along $\Delta_{q_u}^{reg}$ with kernel the non-rtivial 
character $TS_u$. 
Moreover, $\det(p)$ vanishes to order $1$ along $\Delta_{q_u}^{reg}$. 
Decomposing both $V^*$ and 
the restriction $\restricted{TX}{\Delta_{q_u}^{reg}}$ into weight 
spaces of the stabiliser, we see that the horizontal map $f$ 
has rank $1$ over $\Delta_{q_u}^{reg}$ with image 
$\omega_{S}\otimes (V^*)_1$. Moreover, the sheaf-theoretic image of $f$
around $\Delta_{q_u}^{reg}$ is precisely the kernel of (\ref{4.6}). 
We have already observed that (\ref{4.6}) is the pull-back 
$q_u^*\{q_{u_*}[\omega_S\otimes V^*]^W\}$ of the sheaf of $W$-invariant 
sections. The map $f$ is $W$-equivariant, so we have an isomorphism
on the invariant sections, which is what we needed.

Since the map $f^W$ is an isomorphism away from the branching locus,
we are done. 
\EndProof

%***********************************************************************
% A symplectic structure on the relative Prym 
%***********************************************************************
\subsection{A symplectic structure on the relative Prym}
\label{sec-symplectic-structure-on-relative-prym}

Let $h: \Prym \rightarrow U$ be the relative generalized Prym variety.
Its fiber over $u\in U$ is the identity component of 
$[Pic^0(S_u)\otimes_\Integers \chi]^W$. 
In section \ref{sec-construction-of-2-form}
we construct the $V^*$-valued $2$-form $\Omega_V$  
on the relative Jacobian (of degree zero) $\JJ_{\SS/U}$, and show that it is
non-degenerate. 
In section 
\ref{sec-closedness} we prove that it is closed. 
Using the results of section \ref{sec-symplectic-str}
we get a symplectic structure on the relative prym $\Prym_{\SS/U}$. 
The construction is similar to the construction of a symplectic structure 
on the moduli space of Lagrangian sheaves in \cite{donagi-markman-cime}
Ch. 8. 

%***********************************************************************
% Construction of 2-form
%***********************************************************************
\subsubsection{Construction of the 2-form}
\label{sec-construction-of-2-form}

We want  a $W$-invariant $V^*$-valued $2$-form on the
relative Jacobian over the Hilbert scheme of all (smooth) curves 
in $X$ (including curves which are not $W$-invariant).  
It may happen that this Hilbert scheme is singular. We know however that
the locus $U$ of $W$-invariant curves is smooth (Proposition 
\ref{prop-smoothness-of-hilbert-scheme}). 
We will thus carry the construction only over $U$. 

Let 
\begin{equation}
\label{eq-exact-seq-of-normal-bundle}
0 \rightarrow
\tau \rightarrow
T \rightarrow
N \rightarrow 0
\end{equation}
be an exact sequence of vector bundles on an algebraic variety $M$, 
$V$ a vector bundle and
\begin{equation}
\label{eq-f}
f: \tau\otimes V \rightarrow N^*
\end{equation}
a sheaf homomorphism. 
Lemma \ref{lemma-koszul} introduces a symmetry
condition analogous to the Cubic
Condition of \cite{donagi-markman-cime,donagi-markman-cy-threefolds}. 

%**************************************************
% Lemma: Koszul Lemma
%**************************************************
\begin{new-lemma}
\label{lemma-koszul}
The following are equivalent
\begin{enumerate}
%********
% Item 
%********
\item
\label{lemma-item-f-lifts}
The homomorphism $f$ lifts to a $V^*$-valued $2$-form  
\[
\Omega_V \ \in \ H^0([\Wedge{2}T^*]\otimes V^*)
\]
with respect to which $\tau$ is isotropic (and inducing $f$).
%********
% Item 
%********
\item
\label{lemma-item-extension-class-is-symmetric}
The extension class $\epsilon \in H^1(N^*\otimes \tau)$ 
of (\ref{eq-exact-seq-of-normal-bundle}) is mapped via
$f\cup(\bullet): \tau \rightarrow N^*\otimes V^*$ to a class in 
$H^1(N^*\otimes N^*\otimes V^*)$ which is symmetric, i.e., satisfies 
\begin{equation}
\label{eq-extension-class-is-symmetric}
f\cup \epsilon \in H^1((\Sym^2N^*)\otimes V^*).
\end{equation} 
\end{enumerate}
\end{new-lemma}

\noindent
{\sc Proof:}
The proof is a simple Koszul cohomology argument. 
Let $B$ be the kernel 
\[
B := \ker[\Wedge{2}T^* \rightarrow \Wedge{2}\tau^*].
\]
$B$ parametrizes $2$-forms with respect to which $\tau$ is isotropic. 
We get the short exact sequence
\[
0 \rightarrow 
[\Wedge{2}N^*]\otimes V^* \rightarrow 
B \otimes V^* \rightarrow 
(N^*\otimes\tau^*)\otimes V^* \rightarrow 0.
\]
The connecting homomorphism $\delta_1$ in the long exact sequence
\[
0 \rightarrow 
H^0([\Wedge{2}N^*]\otimes V^*) \rightarrow 
H^0(B\otimes V^*) \rightarrow 
H^0((N^*\otimes\tau^*)\otimes V^*) \RightArrowOf{\delta_1}
H^1([\Wedge{2}N^*]\otimes V^*) \rightarrow \dots
\]
is induced by the composition of 
\begin{enumerate}
\item
cup product with the extension class $\epsilon\in H^1(N^*\otimes\tau)$
of (\ref{eq-exact-seq-of-normal-bundle})
\[
H^0(N^*\otimes \tau^* \otimes V^*) \LongRightArrowOf{\cup\epsilon}
H^1\left((N^*\otimes N^*)\otimes (\tau^*\otimes\tau)\otimes V^*\right),
\] 
\item
contraction
\[
H^1\left((N^*\otimes N^*)\otimes (\tau^*\otimes\tau)\otimes V^*\right) 
\RightArrowOf{\iota} H^1(N^*\otimes N^* \otimes V^*), \ \ \ \mbox{and}
\]
\item
wedge product 
\[
H^1(N^*\otimes N^* \otimes V^*)
\RightArrowOf{\wedge}
H^1([\Wedge{2}N^*]\otimes V^*).
\]
\end{enumerate}

\noindent
An element $f$ in $H^0(N^*\otimes \tau^* \otimes V^*)$ lifts to
a $V^*$-valued $2$-form in $H^0(B\otimes V^*)$
if and only if it is in the kernel of $\delta_1$, which is equivalent to 
the condition that $\iota(f\cup \epsilon)$ belongs to the kernel of
$\wedge$. The kernel of $\wedge$ is precisely
$H^1([\Sym^2 N^*] \otimes V^*)$.
\EndProof

\begin{new-lemma}
\label{lemma-kodaira-extension-class-is-symmetric}
The homomorphism
$f\in Hom(TS_u, N^*_{S_u/X}\otimes V^*)$ induced by the 2-form
$\Omega_V$ on $X$ (see (\ref{eq-commutative-diagram-defining-f}))
pairs with the extension class 
$\epsilon_{S_u} \in H^1(S_u,TS_u\otimes N_{S_u/X}^*)$ of
\[
0 \rightarrow 
N^*_{S_u/X} \rightarrow 
(T^*X\restricted{)}{S_u} \rightarrow
T^*S_u \rightarrow 0
\]
to give a class 
\[
\epsilon_{S_u}\cup f \in H^1(S_u,N^*_{S_u/X}\otimes N^*_{S_u/X}\otimes V^*)
\]
which is symmetric, namely, it is in fact in
\[
H^1(S_u,\Sym^2 [N^*_{S_u/X}]\otimes V^*).
\] 
\end{new-lemma}

\noindent
{\sc Proof:}
This is a direct consequence of the direction
\ref{lemma-item-f-lifts} $\Rightarrow$ 
\ref{lemma-item-extension-class-is-symmetric}
of Lemma \ref{lemma-koszul}.
\EndProof

\begin{new-lemma}
\label{lemma-construction-of-2-form-on-relative-Jac}
There exists a unique $W$-invariant section
\[
\Omega_V \ \ \in \ \ 
H^0(\JJ_{\SS/U}, \ [\Wedge{2}T^*\JJ_{\SS/U}]\otimes V^*)
\]
satisfying
\begin{enumerate}
\item
The fibration 
$\JJ_{\SS/U} \rightarrow U$ is $\Omega_V$-isotropic.
\item
The zero section $U\hookrightarrow \JJ_{\SS/U}$
is $\Omega_V$-isotropic, and
\item
Contraction with $\Omega_V$ induces the homomorphism 
$\varphi : h^*TU \rightarrow T^*_{\JJ/U}\otimes V^*$ 
which, as a homomorphism of trivial vector bundles  over 
each Jacobian $J_{S_u}$, 
is the composition of 
\end{enumerate}

\noindent
$
%T_uU\!=\! 
H^0(S_u,N_{S_u/X})^W \subset H^0(S_u,N_{S_u/X})
\LongRightArrowOf{H^0(f)} H^0(S_u,T^*S_u\otimes V^*) \cong 
H^0(J_{S_u},T^*J_{S_u}\otimes V^*).
$
\end{new-lemma}

\noindent
{\sc Proof:}
The uniqueness of $\Omega_V$ is a simple linear algebra exercise. 
Hence, existence locally over the base $U$ (and even 
fiberwize over $U$) implies global existence. 
The existence on each fiber $J_{S_u}$ is proven using the direction
\ref{lemma-item-extension-class-is-symmetric} 
$\Rightarrow$ 
\ref{lemma-item-f-lifts} 
of Lemma \ref{lemma-koszul}.
%We need only to prove that 
The extension class
\[
\epsilon_{J_u} \in H^1(J_{S_u},TJ_{S_u}\otimes T^*_uU) 
\cong Hom\left[H^0(S_u,N_{S_u/X})^W,H^{0,1}(S)^{\otimes 2}\right]
\]
of
\[
0 \rightarrow
[H^0(S_u,N_{S_u/X})^W]^*\otimes \StructureSheaf{J_{S_u}}
\rightarrow (T^*\JJ\restricted{)}{J_{S_u}}
\rightarrow T^*J_{S_u}
\rightarrow 0
\]
is induced by cup product with $\epsilon_{S_u}$

\noindent
$
H^0(S_u,N_{S_u/X})^W\subset H^0(S_u,N_{S_u/X}) 
\LongRightArrowOf{\epsilon_{S_u}}
H^1(S_u,TS_u) \LongRightArrowOf{\cup}
Hom\left[H^0(S_u,T^*S_u), \ H^{0,1}(S_u) \right]
$

\noindent
$
\cong \ H^{0,1}(S_u)^{\otimes 2}.
$

\noindent
The symmetry of $\epsilon_{S_u}\cup f$ implies that 
$\epsilon_{J_u}$ 
pairs with the homomorphism 
\[
H^0(f) \ \ \in \ \ H^0(S_u,N_{S_u/X})^*\otimes H^{1,0}(S_u)\otimes V^*
\] 
to a 
symmetric class. In other words, 
as a class
\[
\epsilon_{J_u}\cup H^0(f) \ \ \in \ \ 
H^{0,1}(S)\otimes [H^0(S_u,N_{S_u/X})^*]^{\otimes 2}\otimes V^*
\]
it is in fact in the image of 
$H^{0,1}(S)\otimes \Sym^2[H^0(S_u,N_{S_u/X})^*]\otimes V^*$. 
The latter is isomorphic to 
\[
H^1(J_{S_u},\Sym^2[H^0(S_u,N_{S_u/X})^*]\otimes V^*)
\] 
(after identifying the weight-1 Hodge structure
of $J_{S_u}$ with that of $S_u$). 
Replacing $H^0(S_u,N_{S_u/X})^*$ by its $W$-invariant subspace, 
we conclude that $\epsilon_{J_u}\cup \varphi$ is  symmetric.
By the direction
\ref{lemma-item-extension-class-is-symmetric} 
$\Rightarrow$ 
\ref{lemma-item-f-lifts} 
of Lemma \ref{lemma-koszul}, $\varphi$ lifts to a section
of $H^0(J_{S_u},(\Wedge{2}T^*\JJ\restricted{)}{J_{S_u}}\otimes V^*)$. 
\EndProof

\begin{cor}
\label{cor-existence-of-2-form-on-relative-prym}
There exists a canonical non-degenerate $2$-form $\omega$ on 
the relative Prym $\Prym^0_{\SS/U}$ with respect to which
$h:\Prym^0_{\SS/U} \rightarrow U$ is a Lagrangian fibration. 
\end{cor}

\noindent
{\sc Proof:}
The existence of $\omega$ on $\Prym^0_{\SS/U}$ is 
equivalent to the existence of $\Omega_V$ on $\JJ^0_{\SS/U}$ (see
Section \ref{sec-symplectic-str}). The non-degeneracy of $\omega$ follows from
the fact that $f^W$ is an isomorphism (Proposition
\ref{prop-smoothness-of-hilbert-scheme}).
\EndProof

%***********************************************************************
% Closedness
%***********************************************************************
\subsubsection{Closedness}
\label{sec-closedness}

The closedness of the $2$-form $\Omega_V$ on $\JJ^0_{\SS/U}$ 
is easier to prove once we extend the $2$-form 
$\Omega_V$ as a $2$-form on the Picard bundle of all degrees 
\[
\JJ^\bullet_{\SS/U} \ := \ \cup_{d\in \Integers} \ \JJ^d_{\SS/U}.
\] 
Denote by 
\[ 
\ell_d \ : \  {\prod^d}_{U}\SS \ \ \ \longrightarrow \ \ \ \prod^d X
\]
the natural morphism from the $d$-th fiber product of $\SS$ over $U$ to
the $d$-th cartesian product of $X$. 
Let 
\[
A_d \ : \  {\prod^d}_{U}\SS \ \ \ \longrightarrow \ \ \ \JJ^d_{\SS/U}
\]
be the Abel Jacobi morphism. Denote by  $\Omega_{V,d}$
the 2-form   on $\prod^d X$ given by 
$\sum_{i=1}^d p_i^*(\Omega_V)$ where $p_i : \prod^d X \rightarrow X$
is projection on the $i$-th factor. 

%**********************************
% Theorem: Closedness and Compatibility
%**********************************
\begin{thm}
\label{thm-compatibility-of-2-forms}
\begin{enumerate}
\item
\label{thm-item-extension-of-V-valued-2-form}
There exists a unique $W$-invariant extension of 
%the $V^*$-valued $2$-form 
$\Omega_V$ from $\JJ^0_{\SS/U}$ to $\JJ^\bullet_{\SS/U}$,
denoted also by $\Omega_V$, satisfying
\begin{enumerate}
\item
\label{thm-item-compatibility-with-group-structure}
The extended $2$-form is compatible with the group structure of 
$\JJ^\bullet_{\SS/U}$ (see Definition 
\ref{def-compatability-of-2-form-with-group-structure}). 
\item
\label{thm-item-symplectic-compatibility-of-abel-jacobi-map-with-X}
We have the equality for all $d\geq 1$ 
\begin{equation}
\label{eq-symplectic-compatibility-of-abel-jacobi-map-with-X}
A_d^*(\Omega_V) \ \ = \ \ \ell_d^*(\Omega_{V,d}).
\end{equation}
\end{enumerate}
\item
\label{thm-item-closedness-and-rank-2}
The extended $V^*$-valued 
$2$-form $\Omega_V$ on $\JJ^\bullet_{\SS/U}$ is closed and
its restriction to $\JJ^1_{\SS/U}$ satisfies the rank $2$ condition
$A_1^*(\Omega_V)\wedge A_1^*(\Omega_V)=0$. 
\item
\label{thm-item-extension-of-2-form-on-hom-sheaf-from-chi-to-pic}
There exists a canonical 
$W$-invariant extension of the $2$-form $\Omega$ from
$\JJ_{\SS/U}^{0}\otimes_{\Integers}\chi$ to 
$\JJ_{\SS/U}^{\bullet}\otimes_{\Integers}\chi$. The extended $2$-form 
$\Omega$ is compatible with the group structure of 
$\JJ^\bullet_{\SS/U}\otimes_{\Integers}\chi$ (see Definition 
\ref{def-compatability-of-2-form-with-group-structure}). 
\item
\label{thm-item-symplectic-structure}
The restriction of $\Omega$ to the relative prym variety $\Prym^0_{\SS/U}$ 
is a symplectic structure and the morphism $\Prym^0_{\SS/U}\rightarrow U$ 
is a Lagrangian fibration. 
\end{enumerate}
\end{thm}

\noindent
{\sc Proof:}
\ref{thm-item-extension-of-V-valued-2-form}) 

\noindent
{\bf Uniqueness:}
If $\Omega'$ and $\Omega''$ are two such extensions, then
the equality (\ref{eq-symplectic-compatibility-of-abel-jacobi-map-with-X})
would imply that $A^*_d(\Omega'-\Omega'')=0$ for all $d\geq 1$.
The dominance of $A_d$ for $d\geq g$ implies the vanishing of
$\Omega'-\Omega''$ on $\JJ_{\SS/U}^d$ for $d\geq g$. Compatibility with the 
group structure implies the equality over all components. 

\medskip
\noindent
{\bf $W$-invariance:} 
Holds for $\JJ^0_{\SS/U}$ by construction. It is also satisfies by 
$\Omega_{V,d}$ on $\prod_d X$. Hence, the unique extension 
satisfying the compatibility conditions 
\ref{thm-item-compatibility-with-group-structure} and 
\ref{thm-item-symplectic-compatibility-of-abel-jacobi-map-with-X}
must also be $W$-invariant. 

\medskip
\noindent
{\bf Existence:} 
Let $D \subset X$ be a divisor 
whose components are isotropic with respect to the $V^*$-valued $2$-form
on $X$. Assume further that the degree $d_D$ of the restriction 
of $\StructureSheaf{X}(D)$ to $S_u$ is not equal to zero. 
Note that the pull-back of a divisor on $\Sigma$ is isotropic. 
%as well as every component of a fiber of $\pi:X\rightarrow \Sigma$.
%In particular, every irreducible component of $red(\Delta_q^\infty)$  
%is isotropic.
%However, components of $\Delta_q^{reg}$ need not be isotropic. 
Denote by $\D$ the section of  $\JJ^\bullet_{\SS/U}\rightarrow U$ 
obtained by pullback of $\StructureSheaf{X}(D)$ to $\SS$. 
It follows from Example 
\ref{example-acihs-of-jacobians-with-symplectic-str-compatible-with-group}
that there exists a unique
extension of $\Omega_V$ to $\JJ^\bullet_{\SS/U}$ which is compatible 
with the group structure and with respect to which $\D$ is isotropic. 
The proof that the extended $\Omega_V$ satisfies the  compatibility
condition  
%\ref{thm-item-compatibility-with-group-structure} and
\ref{thm-item-symplectic-compatibility-of-abel-jacobi-map-with-X}
of the Theorem is similar to that of Lemma 8.12 and Corollary 8.13 
in \cite{donagi-markman-cime}. 

\medskip
\noindent
\ref{thm-item-closedness-and-rank-2}) 
Both closedness and the rank 2 condition are
immediate consequences of part  
\ref{thm-item-extension-of-V-valued-2-form}. 
On $X$ the $2$-form $\Omega_V$ has rank 2. Hence, on $\SS$
\[
A_1^*(\Omega_V)\wedge A_1^*(\Omega_V) = \ell^*_1(\Omega_V\wedge\Omega_V) = 0.
\]

\noindent
Since 
i) $\Omega_V$ is closed on $X$ and  ii) The Abel-Jacobi morphism
is dominant for $d\geq g$ (and submersive for $d\geq 2g-1$), 
then 
the $2$-form is closed on components of $\JJ^\bullet_{\SS/U}$
of degree $\geq g$.
For closedness on components of degree $\leq g$, 
use the fact that the extended form $\Omega_V$ on
$\JJ^\bullet_{\SS/U}$ is invariant under translations by 
sections of $\LB_\D$ 
(\ref{eq-extension-of-torsion-sheaf-to-non-zero-degrees}).

\medskip
\noindent
\ref{thm-item-extension-of-2-form-on-hom-sheaf-from-chi-to-pic}) 
We construct $\Omega$ on $\JJ_{\SS/U}^{\bullet}\otimes_{\Integers}\chi$  from 
the extended $V^*$-valued $\Omega_V$ on $\JJ_{\SS/U}^{\bullet}$
in the same manner as in Section
\ref{sec-symplectic-str}. 

\medskip
\noindent
\ref{thm-item-symplectic-structure}) 
Follows from the closedness and Corollary 
\ref{cor-existence-of-2-form-on-relative-prym}.
\EndProof

%*******************************************************
% Twisted pryms: 
%*******************************************************

\subsubsection{Symplectic structure on twisted pryms}
\label{sec-symplectic-structure-on-twisted-pryms}

In many applications, integrable systems of generalized prym varieties 
are not globally isomorphic to $\Prym^0_{\SS/U}$ but are only 
$\Prym^0_{\SS/U}$-torsors (regarding $\Prym^0_{\SS/U}$ as a group scheme).
In all the examples in sections 
\ref{sec-examples} and \ref{sec-elliptic-K3-systems} the relevant torsors 
naturally embed into $\JJ_{\SS/U}^{\bullet}\otimes_{\Integers}\chi$
and are, locally over the base $U$, cosets of $\Prym^0_{\SS/U}$. 
There is a {\em canonical} symplectic structure 
on any such torsor (independent of a local choice of a translation). 
Simply restrict the $2$-form 
$\Omega$ which was constructed on 
$\JJ_{\SS/U}^{\bullet}\otimes_{\Integers}\chi$ in part 
\ref{thm-item-extension-of-2-form-on-hom-sheaf-from-chi-to-pic}
of Theorem \ref{thm-compatibility-of-2-forms}. 
Equivalently, translate the symplectic structure on $\Prym^0_{\SS/U}$
to the torsor in $\JJ_{\SS/U}^{\bullet}\otimes_{\Integers}\chi$ 
via local sections which are $\Omega$-isotropic
(with $\Omega$ as in part 
\ref{thm-item-extension-of-2-form-on-hom-sheaf-from-chi-to-pic}
of Theorem \ref{thm-compatibility-of-2-forms}).

Recall that $\Prym^0_{\SS/U}$ is the identity component of the group 
$[\JJ^{\bullet}_{\SS/U}\otimes\chi]^W$. Hence, we can split the data
determining a $\Prym^0_{\SS/U}$-torsor into Datum 1)  determining a
$[\JJ^{\bullet}_{\SS/U}\otimes\chi]^W$-torsor and Datum 2) determining
a $\Prym^0_{\SS/U}$-torsor within 
$[\JJ^{\bullet}_{\SS/U}\otimes\chi]^W$. 
Torsors of $[\JJ^{\bullet}_{\SS/U}\otimes\chi]^W$ in 
$\JJ^{\bullet}_{\SS/U}\otimes\chi$ which are locally cosets of 
$[\JJ^{\bullet}_{\SS/U}\otimes\chi]^W$ are parametrized by 
elements in the group
\begin{equation}
\label{eq-group-of-global-cocycles}
Z^1[W, H^0(U,\JJ^{\bullet}_{\SS/U}\otimes\chi)]
\end{equation}
of global cocycles. 
When the integrable system comes from a $(W,V)$-lagrangian fibration $X$, 
global  cocycles commonly arise  from cocycles in 
$Z^1[W,\Pic(X)]$ via the homomorphism 
\[
Z^1[W,\Pic(X)] \ \ \longrightarrow \ \
Z^1[W, H^0(U,\JJ^{\bullet}_{\SS/U}\otimes\chi)].
\]
A cocycle $c$ in
(\ref{eq-group-of-global-cocycles}) is a map
\[
c \ : \ W \ \ \longrightarrow \ \ H^0(U,\JJ^{\bullet}_{\SS/U}\otimes\chi)
\]
satisfying the condition 
$w_1(c(w_2)) + c(w_1) = c(w_1w_2)$, for every $w_1, w_2$ in $W$ 
(recall that the $W$-action is the diagonal action).
We have the global coboundary map
\begin{equation}
\label{eq-coboundary}
\delta :  H^0(U,\JJ^{\bullet}_{\SS/U}\otimes\chi) \longrightarrow 
Z^1[W, H^0(U,\JJ^{\bullet}_{\SS/U}\otimes\chi)]
\end{equation}
where, given a section $\tau$ of $\JJ^{\bullet}_{\SS/U}\otimes\chi$,
the coboundary $\delta\tau$ is the map
\[
(\delta\tau)(w) \ = \ w\tau - \tau. 
\]
Given a cocycle $c$ in (\ref{eq-group-of-global-cocycles})
we get the $[\JJ^{\bullet}_{\SS/U}\otimes\chi]^W$-torsor
\begin{equation}
\label{eq-torsor-in-sheaf-hom-chi-pic}
[\JJ^{\bullet}_{\SS/U}\otimes\chi]^{W,c}   \ \ := \ \ 
\{\LB \in [\JJ^{\bullet}_{\SS/U}\otimes\chi]  \ \mid \ \ 
\delta\LB = c \}
\end{equation}
where $\LB$ is a local section of $\JJ^{\bullet}_{\SS/U}\otimes\chi$
and $\delta$ in (\ref{eq-torsor-in-sheaf-hom-chi-pic})  is
the local analogue of (\ref{eq-coboundary}). 
For example, given any {\em global} section $\tau$ of 
$\JJ^{\bullet}_{\SS/U}\otimes\chi$, the torsor of its coboundary 
$[\JJ^{\bullet}_{\SS/U}\otimes\chi]^{W,\delta(\tau)}$ 
is simply the coset $\tau + [\JJ^{\bullet}_{\SS/U}\otimes\chi]^W$.
%we get the equality
%\[
%\tau + [\JJ^{\bullet}_{\SS/U}\otimes\chi]^W = 
%[\JJ^{\bullet}_{\SS/U}\otimes\chi]^{W,\delta(\tau)}. 
%\]
%I.e., the cocycle of a torsor which is a global coset corresponds to 
%a coboundary and 
We see that isomorphism classes of 
$[\JJ^{\bullet}_{\SS/U}\otimes\chi]^W$-torsors are parametrized 
by classes in $H^1[W,H^0(U,\JJ^{\bullet}_{\SS/U}\otimes\chi)]$ 
which are locally trivial, i.e.,
by the kernel of the homomorphism 
$H^1[W,H^0(U,\JJ^{\bullet}_{\SS/U}\otimes\chi)] \rightarrow 
H^0[U,H^1(W,\JJ^{\bullet}_{\SS/U}\otimes\chi)]$. 

The relevant cocycle $c$ in (\ref{eq-group-of-global-cocycles}) for the 
Hitchin system and for most of the examples in sections 
\ref{sec-examples} and \ref{sec-elliptic-K3-systems}
is identified in \cite{donagi-spectral-covers} section 5.3. 
On the other hand, Datum 2 turns out to include a class in $H^2(W,T)$ 
where $T$ is the torus $\chi\otimes \ComplexNumbers^\times$ 
(\cite{donagi-spectral-covers} section 5.2).

%********************************************************************
% Examples
%********************************************************************
\section{Examples}
\label{sec-examples}

%********************************************************************
% Classical Pryms
%********************************************************************
\subsection{Classical Pryms}

Our first example will be that of classical Prym varieties,
for which $W=\bbz/2$, and $V$ is the sign representation. Let us then suppose 
that our curves $S_u, u\in U$, all have genus $g$, and let $\bbz/2$ act on each of these
curves, with the quotient curves $S_u/(\bbz/2)$ having genus $h$. The Prym varieties
then have dimension $g-h$, and then $U$ must also be of the same dimension. 
Let us suppose that the system of Prym varieties $\Prym\rightarrow U$ has rank two.
Assuming that the genericity conditions are satisfied, one then obtains a 
surface $X$, equipped with an involution $I$ and 
a non-degenerate $V$-valued $I$-invariant 2-form. This is the same as an ordinary 2-form whose
sign is changed by $I$: $I^*\omega= -\omega$. Note that this forces the fixed
points of the involution on $X$ to be curves.

The curves $S_u$ embed in $X$, with normal bundle $K_{S_u}$. 
As we have seen, the deformation theory
for such curves is unobstructed, giving a $g$-dimensional family of curves,
with a $g-h$-dimensional family of $I$-invariant curves. 
One can build \cite{Hu-algebraic-surfaces} the integrable system
of Jacobians corresponding to the $g$-dimensional family, 
giving the proposition:

\begin{prop}
Any rank 2 system of classical pryms has a canonical extension to
 a rank 2 system of Jacobians.
\end{prop}

We can also examine the genericity condition A (condition B is vacuous when 
$v=1$), and show that it holds 
unless $S$ is a hyperelliptic curve with an extra $\Integers/2$-symmetry. 
%except in some rather unusual cases. 
Indeed, the genericity condition is that for every point $p$ on the curve, 
there be a 1-form $\phi$ on the curve which does not vanish at the point with 
$I^*\phi = -\phi$. If there is a form $\psi$ which 
vanishes at $p$ but not at $I(p)$, one can build an anti-invariant 
$\phi$ by taking $\phi= \psi - I^*\psi$. 
Thus, for the genericity condition to fail, one needs:
$$h^0(S,K(-p-I(p)) = g-1,$$
and so 
$$h^0(S, {\cal O}(p+I(p)) = 2.$$
The curve is thus hyperelliptic, with $p$ and $I(p)$ mapping to the same 
point in $\bbp^1$. If the involution $I$ is the hyperelliptic involution, 
the genericity condition holds, for a trivial reason: all the forms are 
anti-invariant. Thus $I$ is not hyperelliptic, and then $I$ descends 
to a non-trivial automorphism of  $\bbp^1$, as it preserves the 
linear system $|{\cal O}(p+I(p)|$.
In any case, this shows that for the genericity condition to fail, 
the curve must be hyperelliptic, with an extra involution commuting with 
the hyperelliptic involution. 

%***************************************************************
% Hitchin systems
%***************************************************************
\subsection{Hitchin systems}
\label{sec-hitchin-systems}

 We now turn to our motivating  example: 
the Hitchin systems for arbitrary reductive
complex groups $G$ \cite{hitchin,hitchin-integrable-system}.  
Fixing such a $G$,  a compact Riemann surface of genus $g$, 
and a degree $k$, we consider the moduli space ${\cal M}$
of stable Higgs pairs $(P,\phi)$, where $P$ is a $G$-principal bundle over 
$\Sigma$ and $\phi$, the Higgs field, is an element of 
$H^0(\Sigma, ad(P)\otimes K_\Sigma)$.
A Zariski open set of this space can be identified with the cotangent bundle 
of the moduli space of stable principal $G$-bundles on $\Sigma$, 
and the symplectic structure on ${\cal M}$ is the cotangent structure. 
In this way, the tangent bundle to 
${\cal M}$ is described by an exact sequence.
\begin{equation}
\label{exact-seq}
0\rightarrow H^0(\Sigma, ad(P)\otimes K_\Sigma)\rightarrow T{\cal M}\rightarrow
H^1(\Sigma, ad(P))\rightarrow 0
\end{equation}

To each element of ${\cal M}$ we can associate a {\it spectral curve}:
one associates to each point $p$ of $\Sigma$ the Weyl group orbit of points in 
$K_\Sigma\otimes\gh|_p$
which lie in the closure of the $G$-orbit of $\phi(p)$. Doing this for all $p$ yields a 
spectral curve $S$ lying in $K_\Sigma\otimes \gh$. $S$ is a $W$-Galois cover of the original $\Sigma$. 
The spectral curves, by the way they are built, only have one point in the closure of any
Weyl chamber, and the branch points of the projection to $\Sigma$ only occur at the 
walls of the Weyl chambers, that is the points with non-trivial stabiliser. 

Let us fix a Cartan subgroup $H$, a Borel subgroup $B$ of 
$G$ which contains $H$, and let $\gh,\gb$ be the corresponding Lie algebras.
 One shows that the lift $q^*P$ of the bundle $P$ to $S$ 
has a natural reduction $P_B$ to $B$ in such a way that 
$q^*\phi$ lies in $ad(P_B)\otimes K_\Sigma$ and the image of
$q^*\phi$ via the characteristic map from 
$ad(P_B)\otimes K_\Sigma$ to $\gh\otimes K_{\Sigma}$ yields precisely the 
curve $S$ \cite{scognamillo-isogenies}. 
The projection on the level of
groups from $B$ to $H$ associates to $P_B$ a bundle $P_H$; a choice of 
theta divisor on the base curve then gives us a canonical way of twisting 
$P_H$ to a $\tilde P_H$ \cite{scognamillo-isogenies}, 
which is invariant under $W$. 
The data of $S, \tilde P_H$ is sufficient to reconstruct
the Higgs pair. 
The connected component containing the trivial bundle
of the set of $W$-invariant bundles
is parametrised by the Prym variety $Pr(S)$ of $S$.  Unfortunately,  
the bundle $\tilde P_H$ need not lie in $Pr(S)$ 
(the discrete data which determines the component of
$\tilde P_H$ includes, in addition, the extension class in
$H^2(W,H)$ of the normalizer $N(H)$, see \cite{donagi-spectral-covers}). 
Nevertheless, we do have
a local version of the parametrization, sufficient for our purposes.

\begin{thm}\cite{donagi-spectral-covers,faltings,scognamillo-isogenies,
hitchin,hitchin-integrable-system}
(a) Let $S_0$ be a spectral curve as above. 
Consider the variety ${\cal N}$ of pairs $(S, \tilde P_H)$,
where $S$ is a  $W$-invariant deformation in $K_\Sigma\otimes \gh$ of  
$S_0$, and  $\tilde P_H$  a $W$-invariant $H$-bundle over $S$, 
lying in $Pr(S)$. Under the correspondence given above, ${\cal N}$ is 
locally isomorphic to  ${\cal M}$.

(b) The projection $(S, P_H)\rightarrow S$ under the identification of (a) 
defines a Lagrangian foliation  of a Zariski open subset of ${\cal M}$. 
\end{thm}

Deformations of the spectral curve correspond as usual to sections of the 
normal bundle. The space $K_\Sigma\otimes \gh $ has a natural $\gh$-valued 
two-form on it, obtained from the cotangent structure on $K_\Sigma$. 
As above, this form gives a map from the
normal bundle to the curves to the bundle $K_S\otimes \gh$. If we impose 
some genericity  on the branch points (that the stabilisers be of order two),
  the regularity condition \ref{cond-regularity}
holds for these curves. In any case, there is an isomorphism 
$$H^0(S, N)^W\simeq H^0(S, K_S\otimes \gh)^W$$ 
Corresponding to the Lagrangian foliation, one thus has at a generic point of 
${\cal M}$, an exact sequence
\begin{equation}
\label{exact-sequence-Hitchin}
0\rightarrow H^1(S, {\cal O}\otimes \gh)^W\rightarrow T{\cal M}\rightarrow
H^0(S,  K_\Sigma\otimes \gh)^W\rightarrow 0
\end{equation}

At a generic point there are two Lagrangian leaves, one coming from the 
cotangent
structure and one from the integrable system of the spectral curves, which 
intersect
transversally. In particular, the leaves of the integrable system define a 
splitting
of the sequence 
(\ref{exact-sequence-Hitchin}) and the symplectic form on ${\cal M}$ is computed from 
Serre duality on the induced splitting 
\begin{equation}
\label{splitting-Hitchin}
T{\cal M}\simeq H^0(\Sigma, ad(P)\otimes K_\Sigma)^W\oplus
H^1(\Sigma, ad(P))^W
\end{equation}
Let us compute this form explicitly. We can assume generically that 
$\phi$ is everywhere a regular element of $\Gotg$. 
Let $(a_1,b_1), (a_2,b_2)$ represent two elements
of (\ref{splitting-Hitchin}). The symplectic form applied to these two elements is 
$$\Omega((a_1,b_1), (a_2,b_2))=<a_1,b_2>_\Sigma- <a_2,b_1>_\Sigma,$$
where $<,>$ denotes Serre duality on $\Sigma$. 
Explicitly, let us suppose that we have a two parameter family
$P(x_1, x_2), \phi(x_1,x_2)$ of elements of ${\cal M}$, so that 
$(a_j, b_j) = \partial _j (\phi, P)$ at $x_j = 0$. Let $\Sigma$ be covered 
by two open sets, $U_0$, the complement of a point $p$, and a disk $U_1$ around that point,
and let $T(x_1, x_2)$ be transition matrices for $P(x_1,x_2)$ with respect to 
this covering. The Higgs fields are then represented by Lie algebra-valued forms
$\phi_i$ on $U_i$ with $\phi_1=T\phi_0$ on the overlaps.  The symplectic form is given 
by 
$$\Omega((a_1,b_1), (a_2,b_2)) = res_p \left[ tr
\left[(\partial_1\phi)(\partial_2 T\cdot T^{-1})-
(\partial_2\phi)(\partial_1T\cdot T^{-1})\right]\right],$$
with $tr$ denoting the Killing form. 

Instead of computing on the base curve $\Sigma$,
we can lift to the spectral curve ($\pi:S\rightarrow \Sigma$), and compute there:
$$\Omega((a_1,b_1), (a_2,b_2))={1\over |W|}(<a_1,b_2>_S- <a_2,b_1>_S).$$
 Now on the spectral curve, we have a reduction to the Borel subgroup $B$
in such a way that $\phi$ lies in $\gb$.   We have exact sequences of groups
$$0\rightarrow N\rightarrow B\rightarrow H\rightarrow 0,$$
where $N$ is the unipotent subgroup.  There is a 
corresponding exact sequence of Lie algebras
$$0\rightarrow \gn\rightarrow \gb\rightarrow \gh\rightarrow 0.$$
If we fix a principal nilpotent element $e$ in $\gn$, any regular element of the 
$\gh$ has a unique representative in $\gb$ of the form $e+ h, h\in \gh$, up to 
the action of the Weyl group. The lift to $S$ resolves this ambiguity, and,
since $\phi$ is everywhere regular, we can choose
a trivialisation of $P$ over $U_0$ such that $\phi(z) = e + h_0(z)$, $z\in S$.
Now let us choose our point $p$ so that $\phi(p)$ is semisimple. On the disk $U_1$,
restricting if necessary, we can conjugate to $\gh$, and so write $\phi$ as 
$h_1(z)\in \gh$. On the overlap of $U_0$ and $U_1$,
$$ h_1(z) = T(z) (e + h_0(z)).$$
The transition matrix $T(z)\in B$ can be split  as $T(z) = T_H(z)T_N(z), T_H\in H,
T_N\in N$. We have then
$$\partial_jT\cdot T^{-1} = \partial_jT_H\cdot T_H^{-1} + 
T_H\partial_jT_N\cdot T_N^{-1}T_H^{-1}.$$
We note that the second term lies in $\gn$, and so gives zero when paired with elements
of $\gh$.
Evaluating the symplectic form, we then have:
\begin{equation}
\label{form-Hitchin}
(<a_1,b_2>_S- <a_2,b_1>_S)=
res_p \left[tr\left[(\partial_1h_1)(\partial_2 T_H\cdot T_H^{-1})-
(\partial_2h_1)(\partial_1T_H\cdot T_H^{-1})\right]\right].
\end {equation}
This last expression is the explicit expression of the Serre duality pairing 
between $H^1(S, {\cal O}\otimes \gh)$ and 
$H^0(S,  K_\Sigma\otimes \gh)$; note that we have ``abelianised'' the pairing by
lifting to the spectral curve.

Now let us choose an equivariant extension of the bundle $P_H$ to a neighborhood
of a spectral curve $S_0$ in $K_\Sigma\otimes \gh$. This extension gives a way,
in effect, of fixing the transition function $T_H$ while varying $\phi$ and the 
spectral curve. One sees from (\ref{form-Hitchin})  that such variations all lie in 
a Lagrangian subvariety. Also, this extension allows us to write
all $H$ bundles on nearby curves as $P_H\otimes P'_H$ where $P'_H$ has degree
zero, and is $W$-invariant. This allows us to write ${\cal M}$ as an integrable 
system of Prym varieties:
$$\Prym\rightarrow U,$$
with the symplectic form defined by Serre duality between the tangent space 
to the zero-section ($H^0(S,  K_\Sigma\otimes \gh)^W$) and the tangent space 
to the fibers ($H^1(S, {\cal O}\otimes \gh)^W)$.

In turn this gives, as we have seen, a $\gh$-valued two-form on the 
corresponding family of Jacobians, again constructed using Serre duality. 
We now want to consider
the pull back of this form to the space $\bbs\rightarrow U$ under the Abel map.
This map usually associates to a point $p$ in a curve $S_u$ the line bundle 
corresponding to the divisor $p-p_0$, where $p_0$ is a base point. 
To obtain an equivariant Abel map, we average over the group:
\begin{equation}
\label{Abel-map}
A(p) = \sum_{w\in W} w(p)-w(p_0)
\end{equation}
Let our base point for $S_u$ be the intersection of $S_u$ with 
a fixed fiber $\pi^{-1}(\lambda)$,
$\pi:K_\Sigma\otimes \gh\rightarrow \Sigma$, and a fixed Weyl chamber.
We now want to compute the symplectic form at $p$.
Let us suppose that we are away from the branch points. 
The projection to $\Sigma$
gives uniform coordinates on all the curves $S_u$, and in effect splits
$T\bbs$ as $TU\oplus T\Sigma$, identifying $T\Sigma$ with $T\SS_u$, and mapping 
$TU$ to the normal bundle 
$N_{S_u}$. With respect to this splitting, the symplectic form is obtained by 
mapping vectors in $T_U$ to $N_{S_u}$, then using the $\gh$-valued  
symplectic form $\Omega_\gh$
on $X=K_\Sigma\otimes \gh$ to map to $K_{S_u}\otimes\gh$, and then pairing 
with $TS_u$. In particular, the fact that we are factoring through $N_{S_u}$
 shows that:

\begin{prop}
\label{thm-rank-two-Hitchin}
The $\gh$-valued symplectic form on $\bbs$ is lifted from
$X=K_\Sigma\otimes \gh$. The Hitchin system for an arbitrary reductive 
group is of rank two,
and the  $(W,V)$-fibration associated to it is the 
quadruple $(\pi: K_S\otimes \gh\rightarrow \Sigma, W, \gh^*, \Omega_\gh)$.
\end{prop}

%***************************************************************
% Coadjoint orbits
%***************************************************************
\subsection{Coadjoint orbits} 
\label{sec-coadjoint-orbits}

Let us consider a reductive complex Lie  group $G$, and  
$\Gotg$ its Lie algebra.
Let $tr$ denote a non-degenerate invariant bilinear on $\Gotg$.
 We can form the loop algebra:
$$
L{\Gotg}=\left\{\sum^k_{i=-\infty}g_i\lambda^i|g_i\in{\Gotg}\right\}
$$
As a vector space, this splits as a sum $L\Gotg_+ \oplus L\Gotg_-$ of 
polynomial and negative-Laurent loops. 
A pairing constructed from $tr$ and taking residues  
identifies (a dense subspace of) $(L\Gotg_+)^*$ with $L\Gotg_-$. 
The dual space $(L\Gotg_+)^*\simeq L\Gotg_-$, with its canonical Lie-Poisson
structure, has many finite dimensional symplectic leaves. 
The symplectic leaves are all coadjoint orbits; when $\Gotg $ is semi-simple, 
the finite dimensional orbits are all of elements of $L \Gotg_-$ of the form:
\begin{equation}
\label{eq-element-in-coadjoint-orbit}
A(\lambda)=N(\lambda)\cdot p(\lambda)^{-1}
\end{equation}
where $ p(\lambda)=\prod^n_{i=1}(\lambda-\lambda_i)$
is some polynomial and $N(\lambda)\in L{\gm g}_+$ is of degree
$<n = deg(p(\lambda))$ \cite{Hu-coadjoint-orbits}. 
One expands around infinity to obtain an 
element of $L \Gotg_-$.
The coadjoint action preserves the form (\ref{eq-element-in-coadjoint-orbit}), 
and the polynomial $p(\lambda)$
is an invariant of the orbit. The conjugacy classes of the polar parts
of $A(\lambda)$ are also invariants of the orbit.  We will consider  only 
these orbits, as we want to consider only the finite-dimensional ones.

As for the Hitchin systems, there is a natural Lagrangian foliation on such 
coadjoint orbits, defined in terms of spectral curves \cite{AvM,RS}.
(Indeed both these systems and Hitchin's fall into a larger class 
\cite{markman-spectral-curves}.)
 Generic elements of 
$L\Gotg_-$  are regular at each $\lambda$, and we will again restrict our 
attention to these.
At each $\lambda$, then, as we are dealing with regular $N(\lambda)$,
one has a unique correspondence under conjugation with
 orbits of the Weyl group $W$ of $G$ in the Cartan subalgebra
$\gh$ of $\Gotg$: for semi-simple elements, this is merely the intersection of
$\gh$ with the 
orbit, and for more general elements, the intersection with the closure
 of the orbit.

The leaves of the Lagrangian foliation then correspond uniquely to 
{\it spectral curves $S^0$}
in $\bbc\times \gh$; to $N(\lambda)$, one associates the points $(\lambda, h)$,
where $h$ lies in the closure of the $G$-orbit of $N(\lambda)$. 
This curve compactifies
to a curve $S$ lying in the total space of the vector bundle 
${\cal O}(n-1)\otimes \gh$
over $\bbp_1$. We will suppose that the curves are smooth; this in fact 
implies that  $N(\lambda)$ is everywhere regular. 
The curve $S$ admits a Galois action of the Weyl group, 
with quotient $\bbp_1$. 
At the zeroes of $p$ , the conjugacy class is fixed along the orbit, 
and so the fiber of the spectral curve is also fixed.

We can reduce the orbits by the action of conjugation by $G$. 
The moment map for this action consists of taking the leading order 
(order $(n-1)$)  term at $\lambda = \infty$.
The reduction then fixes the leading order term, and so the spectral 
curve at infinity. For the reduced orbits, we then obtain spectral curves  
$S$ which all meet in the same $|W|$-orbits in $\gh$ over  
$\lambda_i$ and at $\lambda=\infty$, that is (n+1) points; we will suppose 
that these orbits are generic, of cardinality $|W|$.

We can again define an additional datum for each element of the reduced orbit. 
At generic (semi-simple) points of the spectral curve, 
the lift of the trivial principal $G$-bundle on
$\bbp_1$ to the spectral curve, has a natural reduction $P_H$ to the Cartan 
subgroup $H$, consisting of bases  which conjugate  $N(\lambda)$ to the 
corresponding element of $\gh$. This $H$-bundle
extends to the non-generic points, that is,
 where the curve meets the walls of the Weyl chamber. 
Both curve and bundle are 
invariant under the action of $W$.

This correspondence gives a local isomorphism
(\cite{faltings,scognamillo-isogenies,Hu-Kostant-Kirillov-form})
$${\rm reduced\ orbits}\simeq \{{\rm Inv't\ spectral\ curves,\ inv't\ }
 H{\rm -bundles\ over\ the\ curves}\},
$$
Choosing, as for the Hitchin systems, an equivariant extension of the bundle 
$P_H$ at a spectral curve $S$ to a neighborhood of the curve in 
${\cal O}(n-1)$ 
gives a local identification of this system with 
an integrable system  $\Prym$ of generalised Prym varieties,again over a base
$U$ parametrising the possible spectral curves:
$\Prym \rightarrow U$.

Under the isomorphism of the orbit with the fibration by Pryms, we have a description 
of the tangent space to the reduced orbit as an exact sequence
\begin{equation}
\label{6.5}
0\to H^1(S,{\cal O}\otimes{\gm h})^W\to T({\rm Orbit})\to R\to 0,
\end{equation}
where $R$ is a subspace of $H^0(S,N)^W$. An equivariant extension of the bundle from $S$ to a neighborhood of $S$ in $X$ gives a way of 
moving the curve while fixing the bundle, and so of splitting (\ref{6.5}). 
If we do this, we have an isomorphism:
\begin{equation}
\label{eq-decomposition-of-tangent-space}
T(orbit)\simeq H^1(S,{\cal O}\otimes{\gm h})^W\oplus R,
\end{equation}
 Now the curves $S$ lie in the space ${\cal O}(n-1)\otimes \gh$. 
The coadjoint action fixes the spectral curves over the zeroes of 
$p(\lambda)$, and the reduction fixes it at infinity. 
$R$ is in fact the subspace of $H^0(S,N)^W$ of sections which vanish at 
these points. There is a natural  meromorphic $\gh$-valued symplectic form 
$\tilde\Omega_\gh$ on ${\cal O}(n-1)\otimes \gh$. 
The choice of coordinate $\lambda$ on $\PP^1$ and the polynomial $p(\lambda)$
determine an isomorphism of line bundles 
${\cal O}(n-1)\cong T^*\PP^1(\infty+\sum_{i=1}^n\lambda_i)$ and
on $T^*\PP^1(\infty+\sum_{i=1}^n\lambda_i)\otimes \gh$ we have a canonical 
meromorphic $2$-form. Explicitly, 
let 
%$\lambda$ be a coordinate on $\bbp^1$, and 
$z$
be the fiber coordinate in the standard trivialisation of ${\cal O}(n-1)$.
If $e_i$ is a basis of 
$\gh$, with $z^i$ the corresponding coordinate on ${\cal O}(n-1)\otimes \gh$,
then the form is  given on ${\cal O}(n-1)\otimes \gh$ by 
$\tilde\Omega_{\gh} = p(\lambda)^{-1}
\sum_i (d\lambda\wedge
dz^i) \otimes e_i$. 
The poles of this form coincide with the zeroes of the elements of $R$,
and we contract with it to obtain an isomorphism:
$$R\simeq H^0(S, K\otimes \gh)^W.$$
With this isomorphism, 
(\ref{eq-decomposition-of-tangent-space}) becomes:
$$
T(orbit)\simeq H^1(S,{\cal O}\otimes{\gm h})^W\oplus H^0(S, K\otimes \gh)^W,
$$
The two summands are Serre duals to each other, allowing a natural 
definition of a skew form on $\Prym$. This skew form is independent of the 
extension chosen, and one has:

\begin{prop} \cite{Hu-Kostant-Kirillov-form} 
On the reduced orbits, this skew form is the reduced 
Kostant-Kirillov form.
\end{prop}

This system is an example of the rank two systems discussed above. 
The variety $X$ is obtained from ${\cal O}(n-1)\otimes \gh$ 
by blowing up at the points on the spectral curves lying over 
 the  zeroes $\lambda_1,..,\lambda_n$
of $p(\lambda)$, and over $\infty$,  recalling that these are constant along
the coadjoint orbit. (One blows up taking multiplicities of the 
$\lambda_i$ into account).   The variety $X$ is foliated by projection to $\bbp^1(\bbc)$: this
will be the $\phi_0$-foliation.

The lift $\Omega_\gh$ of the $\gh$-valued 2-form $\tilde\Omega_\gh$ to $X$ is non-singular.
 The Weyl group acts naturally on $X$, and preserves the $\gh$ valued 
form $ \Omega_\gh$ under the simultaneous action of $W$ on $X$ and $\gh$.  

Repeating the argument for  the Hitchin systems, we have 

\begin{prop}
\label{6.} 
The coadjoint orbit systems have rank two, and the $(W,V)$-Lagrangian 
fibrations associated to them are the four-tuples $(\pi:X\rightarrow
\bbp^1, W, \gh^*, \Omega_\gh)$ constructed above. 
\end{prop}

%The following paragraph is excluded (Eyal, March 1, 98)
%This case illustrates a subtlety which deserves some explanation. 
%Theorem 1.1 is essentially local in nature, where local means in a 
%neighbourhood of one of the Prym varieties. These parametrise a component
%of the  degree 0 invariant $H$-bundles on the curve. However, the 
%natural $H$-bundles associated to the 
%coadjoint orbits do not lie in the Prym variety, and indeed are not
% of degree zero. The identification with a system of Prym varieties
%relies on a choice of a bundle $P_H$ on a neighbourhood of the spectral
% curve. More invariantly, if $U$ parametrises the possible spectral curves, 
%and $\Pr\rightarrow U$ is the associated family of Pryms, one has 
%an identification of the coadjoint orbit with an open set of a 
%$\Pr$-torsor ($\Pr$ is a group scheme).

%***************************************************************************
% The $(W,V)$-Lagrangian fibration of a coadjoint orbit
%***************************************************************************
\subsection{
$(W,\LieAlg{h})$-Lagrangian fibrations as invariants of coadjoint orbits
}
\label{sec-the-lagrangian-fibration-of-a-coadjoint-orbit}

Given a coadjoint orbit $R$ of $L\LieAlg{g}$ determined by an element 
$A(\lambda)=N(\lambda)/p(\lambda)$ as in (\ref{eq-element-in-coadjoint-orbit})
it is natural to ask: 
{\em Is the reduced coadjoint orbit a Lagrangian fibration whose generic 
fiber is isogenous to the Prym variety associated to a smooth $W$-Galois 
spectral cover of $\PP^1$?} 
If $R$ is an orbit for which the answer is affirmative, we say that 
$R$ {\em has a smooth spectrum.} 
An analogous question can be formulated for the moduli spaces of 
meromorphic Higgs pairs over a curve of higher genus as in 
\cite{markman-spectral-curves}. 
We know that if 
the leading coefficients of the polar part of $A(\lambda)$ are regular and 
semi-simple, then $R$ has a smooth spectrum. 
For a general orbit, the question seems rather subtle. 
There are examples  of orbits $R$ with singular spectrum. 
For example, if $\LieAlg{g}=\LieAlg{gl}_2$ and the Laurent tail of 
$A(\lambda)$ at $\lambda=0$ has the form
\begin{equation}
\label{eq-coadjoint-orbit-with-singular-spectrum}
A(\lambda) \ \equiv \   
\left(
\begin{array}{cc}
0 & \lambda^{-k}
\\
1 & 0
\end{array}
\right) \ \ \ (\mbox{mod} \ \lambda), \ \ \ \ \ k\geq 2,
\end{equation}
then the generic fiber will be the Picard variety of 
a hyperelliptic curve with a planar singularity of analytic type $y^2=x^k$. 

Theorems \ref{thm-quotient-of-universal-curve}, 
\ref{thm-rigidity-in-introduction}, and 
\ref{thm-from-hilbert-schemes-to-acihs-in-introduction} 
suggest a geometric characterization of the coadjoint orbits with smooth
spectrum. 
The main point is that the $(W,\LieAlg{h})$-lagrangian fibration $X_R$, 
associated to an orbit $R$ with smooth spectrum, 
can be constructed from a simple infinitesimal invariant of $R$. 
In fact, regardless of the smoothness of the spectrum of $R$, 
we can construct a variety $X_R$ with a $W$-action as follows. 
Let $D$ be the divisor $\infty\!+\!\sum_{i=1}^n\lambda_i$ (see section
\ref{sec-coadjoint-orbits}). 
$X_R$ is isomorphic to $X_D:=T^*\PP^1(D)\otimes\LieAlg{h}$
away from the fibers over points in $D$. 
$X_R$ can be obtained from $X_D$ via a complicated blow-up.
It is easier to  construct first the quotient $Y_R := X_R/W$ which is an 
affine bundle over $\PP^1$. Let $Y_D$ be the vector bundle $X_D/W$. 
$Y_D$ is isomorphic to $\oplus_{i=1}^{\dim(\LieAlg{h})}(T^*\PP^1(D))^{d_i}$
where the $d_i$'s are the degrees of the Casimirs of $\LieAlg{g}$. 
$Y_R$ is the affine bundle whose sections are sections 
of $Y_D$ coming from characteristic polynomials of 
elements of $R$. More precisely, 
the coadjoint orbit $R$ determines a pair $(Y'_R,\bar{u}_R)$ 
of a subsheaf $Y'_R$ of $Y_D$ and a section $\bar{u}_R$ of the torsion sheaf
$Y_D/Y'_R$ (supported on a subscheme of $\PP^1$ 
whose set-theoretical support is $D$). 
$Y_R$ is the affine bundle whose sections come from sections of $Y_D$ 
which project to $\bar{u}_R$. 
$Y_R$ is a $Y'_R$-torsor and the two are isomorphic upon a
choice of a section of $Y_R$. 

Consider, for example, the case of a simple root  $\lambda_i$ of
$p$. If $N(\lambda_i)$ is regular and semi-simple, then 
$Y'_R$ is equal to $Y_D(-\lambda_i)$ locally around $\lambda_i$ and 
$\bar{u}_R$ is the characteristic polynomial of $N(\lambda_i)$. 
If, on the other hand, $N(\lambda_i)$ is nilpotent, then around $\lambda_i$ 
the section $\bar{u}_R$ vanishes and $Y_R$  is the subsheaf $Y'_R$ of $Y_D$. 
The length of the stalk of the 
torsion sheaf $Y_D/Y'_R$ at $\lambda_i$ records the order of vanishing of the 
Casimirs along the coadjoint orbit of $N(\lambda_i)$ in $\LieAlg{g}^*$. 
For example, if $N(\lambda_i)$ is  regular nilpotent, then 
$Y'_R=Y_D(-\lambda_i)$ around $\lambda_i$. If $N(\lambda_i)=0$, 
then  $Y'_R= \oplus_{i=1}^{\dim(\LieAlg{h})}(T^*\PP^1)^{d_i}$ around 
$\lambda_i$. 

If $R$ has a smooth spectrum and the vector bundle
(\ref{4.6}') is generated by its global sections, 
then the $(W,\LieAlg{h})$-lagrangian fibration $X_R$ is determined as the 
smooth locus of the normalization of the  fiber product
\[
X_D\times_{Y_D}Y_R
\]
(the fiber product is singular along a hypersurface in general). 
Often, the coadjoint orbit $R$ determines a singularity for {\em all}
spectral curves in $X_D$ and the morphism
$X_R\rightarrow X_D$ resolves the singularity of the generic curve. 
By Theorems \ref{thm-quotient-of-universal-curve} and 
\ref{thm-rigidity-in-introduction},
the natural $\LieAlg{h}$-valued meromorphic $2$-form 
$\tilde{\Omega}_\LieAlg{h}$ on $X_D$ pulls back to a holomorphic 
$2$-form on $X_R$ which satisfies the minimal degeneracy condition 
of Definition \ref{def-V-Lagrangian-fibration} (see also part
\ref{rem-item-minimal-degeneracy-is-satisfied-by-quotient-of-prym-acihs} 
of remark \ref{rem-after-def-of-W-V-lagrangian-fibration}). 
%{\em We get a necessary
%condition for $R$ to have a smooth spectrum. We expect that 
%the condition is also sufficient.} \ In other words, we expect that $X_R$ is 
%never a $(W,\LieAlg{h})$-lagrangian fibration if $R$ has a singular 
%spectrum. 
This proves the implication 
\ref{conj-item-R-has-a-smooth-spectrum} $\Rightarrow$ 
\ref{X-R-is-a-lagrangian-fibration} in the following conjecture.
Theorem \ref{thm-compatibility-of-2-forms} supports the reverse implication. 

\medskip
\noindent
{\bf Conjecture:}
{\em
The following are equivalent:
\begin{enumerate}
\item
\label{conj-item-R-has-a-smooth-spectrum}
The orbit $R$ has a smooth spectrum and the vector bundle (\ref{4.6}')
of its generic spectral curve $S_u$ 
is generated by its global sections. 
\item
\label{X-R-is-a-lagrangian-fibration}
$X_R$ is a $(W,V)$-lagrangian fibration and $Y'_R$ is generated by its global 
sections.
\end{enumerate}
}

\begin{example}
{\rm
It is easy to see that for $\LieAlg{g}=\LieAlg{sl}_2$ and $R$ 
the orbit with singular spectrum given by
(\ref{eq-coadjoint-orbit-with-singular-spectrum}), the surface $X_R$ 
is {\em not} symplectic. 
%In this case, $\LieAlg{h}$ is the non-trivial character of $W=\Integers/2$. 
We get around $\lambda_i$ the equalities: 
\begin{eqnarray*}
Y_D & = & [T^*\PP^1]^{\otimes 2}(2k\lambda_i)
\ \ \ \ \mbox{and} \\ 
Y_R & = & Y'_R \ = \ Y_D(-k\lambda_i)\ = \ [T^*\PP^1]^{\otimes 2}(k\lambda_i)
\end{eqnarray*}
while 
$X_R =  T^*\PP^1(\frac{k}{2}\lambda_i)$ if $k$ is even and has a fiber 
of multiplicity $2$ over $\lambda_i$ if $k$ is odd. In the odd case, $X_R$ 
is obtained by 1) blowing-up $T^*\PP^1(\lceil\frac{k}{2}\rceil\lambda_i)$
at the zero point in the fiber over $\lambda_i$, 2) blowing-up again
at the intersection point of the two components in the fiber, and 
3) excluding the two previous components in the fiber. 
In both cases, $\Omega_{\LieAlg{h}}$ is a section of 
$\omega_{X_R}(\lfloor\frac{k}{2}\rfloor f_i)$ where $f_i$ is the 
fiber of $X_R$ over $\lambda_i$.  
Thus, the $2$-form $\Omega_\LieAlg{h}$ is meromorphic when $k\geq 2$. 
Contraction of sections of the normal bundle of a (smooth) curve in
$X_R$ with $\Omega_\LieAlg{h}$ results in a meromorphic $1$-form on the curve. 
This corresponds to the fact that the dualizing line-bundle of the
singular curve pulls-back to the normalization as a sheaf of meromorphic
$1$-forms. 
}
\end{example}

%We do not know if the map $R\mapsto (Y'_R,\bar{u}_R)$ is injective in 
%general. It seems almost certain however that it is injective on the subset 
%of orbits with smooth spectrum. 

%******************************************************************
% Special cases
%******************************************************************
\subsection{Special cases}

For both the Hitchin system and the coadjoint orbits, further reductions are 
possible for particular choices of the reductive group: as both cases are 
quite similar, we will treat the Hitchin case.

%******************************************************************
% $GL(n, \bbc)$
%******************************************************************
\subsubsection{$GL(n, \bbc)$}
In this case, one has spectral curves sitting inside 
$K_\Sigma\otimes \bbc^n$, with the permutation group ${\cal S}_n$ acting on 
$\bbc^n$. Projecting from $\bbc^n$ to the first coordinate maps the $n!$-fold covering
$S$ of $\Sigma$ to an $n$-fold covering $\hat S$ lying in $K_\Sigma$. In a parallel
fashion, one takes the character $\rho = (1,0,0,..,0)$ and forms the line bundle
$L= P_H\times_\rho \bbc$, which descends to $\hat S$. (The curve $S$ corresponds to 
all orderings of eigenvalues, and the map to $\hat S$ picks out one of these,
the first. Similarly, $P_H$ corresponds to the set of frames consisting of 
eigenvectors, and one just picks out the first to form $L$). One can reconstitute 
the pair $S, P_H$ from $\hat S, L$. The map $(S, P_H)\mapsto (\hat S, L)$ 
converts the rank two system of Pryms to a rank two system of Jacobians, whose
corresponding
 surface is $K_\Sigma$.
%******************************************************************
% $SL(n, \bbc)$
%******************************************************************
\subsubsection{$SL(n, \bbc)$}
Here one again has the permutation group ${\cal S}_n$
acting on spectral curve, which are now embedded in $K_\Sigma\otimes \bbc^{n-1}$,
where $\bbc^{n-1}\in\bbc^n$ is the set of elements whose coordinates sum to zero.
 Taking this embedding and mapping to the first coordinate as before, one 
obtains a curve $\hat S$ lying in $K_\Sigma$ such that over each point of $\Sigma$,
the $n$ points of $\hat S$ lying above that point  sum to zero in the fiber of 
$K_\Sigma$. In turn, the line bundle $L$ that one obtains by following the procedure 
outlined above has the property that the determinant of $E=\hat \pi_*L$ is constant,
 $\hat \pi:\hat S\rightarrow \Sigma$. The system has a canonical extension to a rank
two system of Jacobians, obtained by dropping both the sum constraint on the curve,
and the constraint on the bundle $E$. 
This amounts to passing from $SL(n, \bbc)$ 
to $GL(n, \bbc)$. 

%******************************************************************
% $SO(2n,\bbc)$ or $Sp(2n, \bbc)$
%******************************************************************
\subsubsection{$SO(2n,\bbc)$ or $Sp(2n, \bbc)$}
In this case, one has an extension 
of $(\bbz/2)^n$ by ${\cal S}_n$, acting on curves embedded in 
$K_\Sigma\otimes \bbc^n$. Projecting to the first coordinate yields curves 
$\hat S$ in $K_\Sigma$ which are $2n$-fold covers of \
$\Sigma$ under the projection, and which are invariant under the involution $I$
on the fibers of $K_\Sigma\rightarrow \Sigma$ given by the action of $-1$.
 On the level of line bundles, one obtains elements lying in the classical 
Prym varieties of $\hat S$. The rank two system of $W$-Pryms 
(with corresponding curves $S$) is equivalent under this map to a rank
two system of $\bbz/2$-Pryms (with corresponding curves $\hat S$).

%*******************************************************
% Principal bundles over Poisson elliptic surfaces:
%*******************************************************

\section{Principal bundles over Poisson elliptic surfaces}
\label{sec-elliptic-K3-systems}

Examples of $(W,V)$-Lagrangian fibrations 
$\pi:X\rightarrow \Sigma$ with compact fibers 
%whose fibers are abelian varieties, 
arise naturally in the  study of
principal bundles over Poisson elliptic surfaces. 
The moduli space of principal $G$-bundles on the surface turns out to be
an associated prym integrable system. These moduli spaces 
play a  central role in the duality
between heterotic string and $F$ theory compactifications. 

%*******************************************************
% $(W,V)$-Lagrangian fibrations with compact fibers
%*******************************************************
\subsection{$(W,V)$-Lagrangian fibrations with compact fibers}
\label{sec-compact-fibers}

Let $(Z,\psi)$ be a smooth projective Poisson surface, 
$\psi \in H^0(Z,\Wedge{2}TZ)$ the Poisson structure, 
$p:Z \rightarrow \Sigma$ an elliptic fibration, and
$\sigma : \Sigma \hookrightarrow Z$ a section. 
%We require that the degeneracy divisor of $\psi$ is a linear combination
%of elliptic fibers. 
For simplicity, we  assume
that the fibers of $p$ are all irreducible. It follows from the existence of 
the section $\sigma$ that all fibers are also reduced. 
% $p:Z \rightarrow \Sigma$ is a Weierstrass model. 
Every such surface $S$ is a double cover of a ruled surface
\[
\eta \ : \ Z \ \rightarrow Z/\iota
\]
where $\iota$ is the elliptic involution
(see \cite{friedman-morgan-book} Ch. 1 Theorem 4.4). 
Denote by $L$ the conormal bundle of $\sigma$ in $Z$. 
$L$ is the line bundle over $\Sigma$ 

\begin{equation}
\label{eq-L}
L \ := \ [R^1p_*\StructureSheaf{Z}]^*. 
\end{equation}

\noindent
There exist sections $g_2$ of $H^0(\Sigma,L^{\otimes 4})$
and $g_3$ of $H^0(\Sigma,L^{\otimes 6})$ such that $Z$ is the zero 
divisor in the threefold 
$\PP[L^{\otimes 2}\oplus L^{\otimes 3}\oplus \StructureSheaf{\Sigma}]$ 
of the section
\[
y^2z - 4x^3 - g_2xz^2 - g_3z^3
\]
of $\StructureSheaf{}(3)\otimes L^6$ where $\StructureSheaf{}(1)$ 
is the tautological line bundle on the $\PP^2$ bundle. Above, 
$x$, $y$, and $z$ are the tautological sections of 
$\StructureSheaf{}(1)\otimes L^2$,
$\StructureSheaf{}(1)\otimes L^3$, and 
$\StructureSheaf{}(1)$ respectively. 
The projection 
\begin{equation}
\label{eq-elliptic-involution}
\eta \ : \ Z \rightarrow \PP[L^{\otimes 2}\oplus \StructureSheaf{\Sigma}]
\end{equation}
is the quotient by the elliptic involution. 
The branch locus of $\eta$ belongs to the linear system 
$\linsys{\StructureSheaf{}(4)\otimes L^{\otimes 6}}$ on 
$\PP[L^{\otimes 2}\oplus \StructureSheaf{\Sigma}]$. 
The branch locus decomposes as the disjoint union
$\bar{\sigma} + \bar{\Delta}_\eta$ where $\bar{\sigma}$ is the section at 
infinity $z=0$ and $\bar{\Delta}_\eta$ is a trisection of the ruling
defined by $4x^3 - g_2xz^2 - g_3z^3$. The components 
$\bar{\sigma}$ and 
$\bar{\Delta}_\eta$ are in the linear systems 
$\linsys{\StructureSheaf{}(1)}$ and 
$\linsys{\StructureSheaf{}(3)\otimes L^{\otimes 6}}$ respectively. 

There are four families of examples of such surfaces depending on
the degree of $L$ and the genus of the base curve $\Sigma$ 
(see \cite{friedman-morgan-book} 
Ch. 1 Lemma 3.18, Proposition 3.23, and Theorem 4.3).
The canonical bundle of $Z$ is 
\[
K_{Z} \ = \ p^*(K_\Sigma \otimes L). 
\]
In general, $\deg(L)$ in non-negative. 
The assumption that $Z$ is a Poisson surface implies that 
$\deg(L)$ is either $0$, $1$, or $2$. Moreover, the base curve 
$\Sigma$ is either $\PP^1$ or an elliptic curve. 

\begin{enumerate}
\item
\label{item-deg-L-0-Sigma-elliptic} 
If $\deg(L) =0$ and $\Sigma$ is an elliptic curve, then $L$ is 
the trivial  line bundle (the non-degeneracy of the section $\psi$ 
along the section $\sigma$ rules out the 
possibility that $L$ is a line bundle of order $2$ on $\Sigma$). 
The symplectic surface $Z$ is a product of two elliptic curves. 
\item
\label{item-deg-L-0-Sigma-P-1} 
If $\deg(L) =0$ and $\Sigma=\PP^1$, then $Z$ is the product $E\times \PP^1$ 
of an elliptic curve $E$ and $\PP^1$. In this case, $\psi$ is the pull-back 
of a section $\bar{\psi}$ of 
$\omega_{\PP^1}^{-1}\cong \StructureSheaf{\PP^1}(2)$. 
We have two sub-cases: $\bar{\psi}$ has either two simple zeroes, 
or a double zero. 
\item
\label{item-deg-L-1}
If $\deg(L) =1$, then $\Sigma=\PP^1$ and $Z$ is the blow-up of
$\PP^2$ at $9$ points which are the base points of a pencil of plane 
cubic curves. $Z$ is a double cover of the Hirzebruch surface ${\Bbb F}_2:= 
\PP[\StructureSheaf{\PP^1}(2)\oplus \StructureSheaf{\PP^1}]$. 
The Poisson structure $\psi$ is degenerate along a single elliptic fiber. 
Denote by $\bar{f}$ the class of the fiber of ${\Bbb F}_2 \rightarrow \PP^1$.
Then the trisection  $\bar{\Delta}_\eta$ is a curve of genus 4
in the linear system
$\linsys{\StructureSheaf{{\Bbb F}_2}(3\bar{\sigma}+6\bar{f})}$. 
\item
\label{item-deg-L-2}
If $\deg(L) =2$, then $\Sigma=\PP^1$ and $Z$ is a (symplectic) 
K3 surface which is 
a double cover of the Hirzebruch surface ${\Bbb F}_4:= 
\PP[\StructureSheaf{\PP^1}(4)\oplus \StructureSheaf{\PP^1}]$. 
We have the equality 
\[
\omega_{{\Bbb F}_4}^{-2} \ = \ 
\StructureSheaf{{\Bbb F}_4}(4\bar{\sigma}+12\bar{f}).
\]
The trisection  $\bar{\Delta}_\eta$ is a curve of genus 10 
in the linear system
$\linsys{\omega_{{\Bbb F}_4}^{-2}(-\bar{\sigma})}$. 
%The ramification locus of
%$\eta$ (\ref{eq-elliptic-involution}) 
%is the disjoint union $\sigma + \Delta_\eta$
%where $\Delta_\eta$ is in the linear system 
%$\linsys{\StructureSheaf{Z}(3\sigma+6f)}$.
\end{enumerate}

\bigskip
Let $G$ be a semi-simple simply-connected complex Lie group, 
$H$ a maximal torus, $\Lambda := \Hom(H,\ComplexNumbers^\times)$
its weight (and root) lattice and $\chi := \Lambda^*$ the dual lattice. 
Denote by $Z^0$ the locus in $Z$ where $p$ is a smooth fibration. 
We consider  $Z^0$  as a group scheme,  taking our section $\sigma$ as the 
zero-section, and form the quasi-projective 
variety 

\begin{equation}
\label{eq-the-lagrangian-fibration-of-an-elliptic-surface}
X \ := \ Z^0 \otimes_{\Integers} \chi. 
\end{equation}

\noindent
If the surface $Z$ is symplectic, 
there is a natural $\LieAlg{h} := \chi\otimes_{\Integers}\ComplexNumbers$ 
valued $2$-form $\Omega_\LieAlg{h}$ on $X$ and 
$\pi:X\rightarrow \Sigma$ is a $(W,\LieAlg{h})$-Lagrangian fibration. 
If $(Z,\psi)$ is a Poisson surface and $\psi$ has  simple zeroes 
along the fiber over  $a_1$ or the  two fibers over  $a_1,a_2$ in 
$\Sigma$ (cases 3,2 above), 
then there is a natural meromorphic $\LieAlg{h}$-valued 
$2$-form $\Omega_\LieAlg{h}$ on $X$ with simple poles along 
the fibers over $a_i\in \Sigma$. 
Denote by $R$ a choice of a generic
$W$-orbit in each of the polar fibers of $X$ over $\Sigma$. 
The blow-up $X_R$ of $X$ along the points ``in'' $R$ is then 
a $(W,\LieAlg{h})$-Lagrangian fibration 
(after excluding the proper transform of the polar fibers). 
$Z^0$ may be thought off as a non-linear version of the line-bundle
$T^*\Sigma(\sum_{i=1}^ka_i)$ and $X_R$ as a non-linear analogue 
of the $(W,\LieAlg{h})$-Lagrangian fibration
associated to semi-simple and regular coadjoint orbits $R_{a_i}$,
$1\leq i \leq k$,  of 
% the positive half of the untwisted loop-group of 
$G$.
It is interesting to note \cite{morrison-vafa} that elliptic K3 surfaces as in 
case \ref{item-deg-L-2} admit stable degenerations to the union $Z_1\cup Z_2$
of two Poisson surfaces as in case \ref{item-deg-L-1} intersecting along a 
common polar elliptic fiber $Z_1\cap Z_2$.

% $(W,\LieAlg{h})$-Lagrangian fibrations analogous to
%more general choices of coadjoint orbits can be constructed as well. 
%When $\psi$ has higher order poles along a fiber, we get 
%$(W,\LieAlg{h})$-Lagrangian fibrations analogous to choices of 
%coadjoint orbits $R_{a_i}$ of 
%the positive half of the untwisted loop-group of 
%$G$.

\bigskip
Let $q : X \rightarrow Y$ be the quotient of $X$ by the natural $W$-action. 
Theorem \ref{thm-compatibility-of-2-forms} associates 
Prym integrable systems to (Zariski open
subsets of) Hilbert schemes $U$ of sections of $Y \rightarrow \Sigma$
(when $Z$ is symplectic, otherwise work with 
%Similarly, we have 
the quotient $q: X_R \rightarrow Y_R$).
% in case $Z$ is Poisson but not symplectic. 
A theorem of Looijenga \cite{looijenga} 
and Bernshtein-Shvartsman \cite{Bernshtein-Shvartsman}
implies that the fiber
\[
Y_a \ := \ [E_a\otimes_{\Integers} \chi]/W, \ \ \ a\in \Sigma
\]
is a weighted projective space when $E_a$ is a smooth elliptic fiber of 
$p:Z\rightarrow \Sigma$ 
(see also \cite{friedmann-morgan-witten2} Theorem 4.3). 
Notice that $E_a\otimes_{\Integers} \chi$ is the identity component of 
the moduli space of principal $H$-bundles over $E_a$, 
while $Y_a$ is the moduli space of equivalence classes, 
under the action of the inner automorphisms of the root system,  
of the $H$-bundles over $E_a$. 
In the next section a reduction theorem will imply 
that $Y$ is also (a subset of)  the 
{\em relative} moduli space $\M_{Z/\Sigma}\rightarrow \Sigma$ 
of principal $G$-bundles
on the elliptic fibration $Z\rightarrow \Sigma$. The moduli space $\M$ of 
semi-stable 
principal $G$-bundles on $Z$ will turn out to be a Prym integrable system.

%*******************************************************
% Principal bundles on elliptic fibrations
%*******************************************************
\subsection{Principal bundles on elliptic fibrations}
\label{sec-principal-bundles-on-elliptic-fibrations}
We describe briefly, 
following \cite{donagi-principal-bundles-elliptic-fibrations,
friedmann-morgan-witten1,friedmann-morgan-witten2}, 
the abelianization of the moduli space $\M$ of 
principal $G$-bundles on $Z$. By  abelianization we refer to a rational 
morphism $h: \M \rightarrow U$ 
%from $\M$ to a Hilbert scheme of sections 
%of $Y\rightarrow \Sigma$ 
whose generic fiber is a generalized Prym.

%*******************************************************
% Reduction 
%*******************************************************
\subsubsection{Reduction of a $G$-bundle on an Elliptic curve to a $H$-bundle}
\label{sec-reduction} 

A vector bundle $V$ on a curve $E$ is {\em stable (resp. semi-stable)}
if every proper non-zero subbundle $V'$ satisfies the inequality
%\[
$
\frac{deg(V')}{rank(V')} \ < \  \frac{deg(V)}{rank(V)}
\ \ \ \ (resp. \leq).
$
%\]
A principal $G$-bundle $P$ on $E$  is {\em semi-stable} if 
the associated vector bundle $ad(P)$ is semi-stable. 
In the construction of the moduli space of semi-stable 
$G$-bundles on $E$ one has to parametrize equivalence classes coarser then 
isomorphism classes. 
The difficulty arises from jump phenomena: the existence of flat
families of principal bundles ${\cal P}$ over $E\times B$, 
$B$ an irreducible parameter scheme, such that the restriction of ${\cal P}$
to a generic fiber $E\times \{b\}$, $b \in B$, are all isomorphic to a fixed
semi-stable bundle $P_1$ but  the restriction to some fiber $E\times \{b_0\}$ 
is another semi-stable  $G$-bundle $P_0$. 
Two such bundles $P_0$ and $P_1$ are called $S$-equivalent. 
Consider the case  where $G$ is $SL_2$ and $E$ is elliptic. Choose a
line bundle $F$ of order two on $E$. 
Then $\Ext^1(F,F)$ is one-dimensional and there is a universal rank $2$ 
vector bundle ${\cal V}$ on $E\times \Ext^1(F,F)$ with trivial determinant 
whose restriction to $E\times \{b\}$, $b\neq 0$, 
is the unique non-trivial extension of $F$ by $F$ while its restriction to 
$E\times \{0\}$  is the direct sum $F\oplus F$. 
The generic $S$-equivalence class of a rank $2$ semi-stable vector bundle 
with trivial determinant is the isomorphism class $F\oplus F^{-1}$
for some $F\in \Pic^0(E)$ but there are four coarser $S$-equivalence classes 
corresponding to the four line bundles of order $2$. 

Assume now that the curve is a smooth elliptic curve $E$. 
The main reduction theorem for $G$-bundles on $E$ is:

\begin{thm}
\label{thm-reduction} 
(\cite{friedmann-morgan-witten2} Proposition 3.9)
Let $P$ be a semi-stable principal $G$-bundle over $E$. 
\begin{enumerate}
\item
\label{thm-item-reduction-to-H-bundle}
$P$ is $S$-equivalent to a principal $G$-bundle $P'$ which 
admits a reduction to a $H$-bundle unique up to the $W$-action
on $H$.
\item
\label{thm-item-reduction-to-regular-bundle}
$P$ is $S$-equivalent to a principal $G$-bundle $P''$ 
which is {\em regular} in the sense that 
\[
\dim H^0(E,ad(P)) \ = \ rank (G). 
\]
$P''$ admits a reduction to a $rank(G)$-dimensional
connected commutative subgroup of $G$ 
which is a regular centralizer and is unique up to conjugation. 
\end{enumerate} 
\end{thm}

The case of $SL_n$ is due to Atiyah \cite{atiyah}. 
A semi-stable vector bundle $V$ with trivial determinant line-bundle 
over an elliptic curve $E$ is $S$-equivalent 
to a direct sum of line-bundles of degree zero whose tensor product 
is the trivial line bundle
\[
V \ \equiv \ V':=\oplus_{j=1}^n F_j, \ \ \ \ F_j \in \Pic^0E, \ \ \ \ 
\mbox{and} \ \ \ \ \otimes_{j=1}^n F_j = \StructureSheaf{E}. 
\]
It is also $S$-equivalent to a direct sum of indecomposable subbundles 
\[
V \ \equiv \ V'':=\oplus_{j=1}^k V_j, \ \ \ \ 
\otimes_{j=1}^k \det(V_j) = \StructureSheaf{E}  \ \ \ \ 
\mbox{and} \ \ \det(V_j) \ \ \mbox{are pairwise non-isomorphic}.
\]
$V''$ is a regular bundle. 

Theorem \ref{thm-reduction} part 
\ref{thm-item-reduction-to-H-bundle} introduces a bijection 
between the moduli space $\M_E$ of $S$-equivalence classes of 
semi-stable $G$-bundles on $E$ and the quotient
$[E\otimes_{\Integers}\chi]/W$. It can be shown  that this bijection is
in fact an isomorphism of algebraic varieties
(\cite{friedmann-morgan-witten2} Theorem 4.3). 
For $SL_n$ it is easy to see that both $[E\otimes_{\Integers}\chi]/W$
and $\M_E$ are naturally isomorphic to the 
linear system $\linsys{\StructureSheaf{E}(np_0)}$, 
where $p_0$ is the marked point on $E$ (which is used to identify $E$
with $\Pic^0E$).

%*******************************************************
% Relative moduli spaces 
%*******************************************************
\subsubsection{Relative moduli spaces of $G$-bundles}
\label{sec-relative-modili-spaces}

We will assume from now on that $G$ is simple and exclude the
exceptional group $E_8$. 
%If $G$ is $E_8$ we will
%impose the restriction that the singularities of fibers of the elliptic
%fibration $p:Z\rightarrow \Sigma$ are simple nodes. 
Denote by 
$\pi:\M_{Z/\Sigma} \rightarrow \Sigma$ the {\em relative}
moduli space of semi-stable principal $G$-bundles over the elliptic fibration
$p:Z\rightarrow \Sigma$. 
A point in $\M_{Z/\Sigma}$ over $a\in \Sigma$ parametrizes an $S$-equivalence 
class of semi-stable principal $G$-bundles on the elliptic curve $E_a$. 
For example, if $G=SL_n$, then $\M_{Z/\Sigma}$ is naturally isomorphic 
to the $\PP^{n-1}$-bundle 
$
\PP\left[
\StructureSheaf{\Sigma}\oplus L^{-1}  \oplus \cdots \oplus L^{1-n}
\right]
$
where $L$ is the line bundle introduced in  (\ref{eq-L}). 
For all simple simply-connected groups $G$, other than $E_8$, 
$\M_{Z/\Sigma}$ is naturally isomorphic to the bundle of weighted projective 
spaces obtained from the vector bundle
\[
\StructureSheaf{\Sigma} \oplus L^{-d_1}   \oplus \cdots \oplus L^{-d_r}
\]
by a $C^{\times}$-action with weights $g_0, \dots ,g_r$ 
(\cite{friedmann-morgan-witten1},\cite{friedmann-morgan-witten2} Theorem 4.4). 
Above, $\{d_1, \dots ,d_r\}$ are the degrees of the homogeneous generators of 
the algebra of invariant polynomials on $\LieAlg{g}$. 
The weight $g_0$ is $1$ and the 
weights $\{g_1, \dots ,g_r\}$ 
are coefficients in the linear combination expressing the co-root
dual to the highest root in terms of the co-roots dual to the simple roots. 

The isomorphism $[E_a \otimes_{\Integers}\chi]/W \cong \M_{E_a}$
fits nicely in families and extends to the singular fibers to
give: 

\begin{thm}
\label{thm-W-quotient-is-the-relative-moduli-of-bundles}
(\cite{friedmann-morgan-witten2} Theorem 4.4)
The quotient $Y$ of $X$ by the natural $W$-action 
embeds as a Zariski open subset of $\M_{Z/\Sigma}$. 
\end{thm}

%Theorem \ref{thm-W-quotient-is-the-relative-moduli-of-bundles} 
%is a relative version of a reduction theorem which states that every 
%$S$-equivalence class of a semi-stable principal $G$-bundle on a 
%smooth  elliptic curve $E$ contains a representative which admits
%a reduction to a principal $H$-bundle over $E$. The principal $H$-bundle is 
%unique up to the action by $W$, the group of inner automorphisms 
%of the root system. 

%***********************************************************************
% A prym integrable system of $G$-bundles on an Elliptic Poisson surface 
%***********************************************************************
\subsubsection
{A prym integrable system of $G$-bundles on an Elliptic surface}
\label{sec-prym-integrable-system-of-G-bundles-on-elliptic-surface} 

We denote by $\M$ the moduli space of semi-stable principal $G$-bundles over 
the surface $Z$ 
(the definition of semi-stability depends on a choice of an ample
line bundle $H$ on $Z$ which we implicitly assume). 
A semi-stable $G$-bundle $P$ on $Z$ determines 
a section $u(P)$ of $\M_{Z/\Sigma}$ because its restriction
to a generic elliptic fiber of $Z\rightarrow \Sigma$ is also semi-stable. 
Consider the Zariski open (possibly empty) subset $\M^0$ of $\M$
parametrizing $G$-bundles $P$ satisfying:
\begin{enumerate}
\item
The restriction of $P$ to every fiber of $Z\rightarrow \Sigma$ is 
semi-stable.
% and regular (in the sense of part 
%\ref{thm-item-reduction-to-regular-bundle} of Theorem \ref{thm-reduction}). 
\item
\label{item-condition-on-section-in-U}
The section $u(P)$ is contained in the smooth locus of 
$Y \subset \M_{Z/\Sigma}$ and is transverse to the branch locus of 
$q:X\rightarrow Y$. 
\end{enumerate}

\noindent
Let $U$ be the  Zariski open subset of the Hilbert scheme of sections of 
$Y\rightarrow \Sigma$ which satisfy \ref{item-condition-on-section-in-U}
above. Note that both $\M^0$  and $U$ have, in general, infinitely many 
components. 
We get a natural morphism $\M^0 \rightarrow U$, associating to 
each bundle a smooth $W$-Galois cover $q_u : S_{u(P)}\rightarrow \Sigma$.
Following \cite{ donagi-principal-bundles-elliptic-fibrations,
friedmann-morgan-witten3}, the fiber is a generalized (twisted) Prym.
The idea, roughly, is that the spectral curve tells us what is the restriction
of the bundle to each fiber; the element of the generalised Prym tells us 
how these are glued together, in essence along the section $\sigma(\Sigma)$.

We briefly, and again roughly, summarise how this goes. 
The curve $S_u$ over each point $a$ in $\Sigma$, is a $W$-orbit of reductions
to $H$ of (an element in the $S$-equivalence class) of the restriction 
of the bundle to $E_a$. Another way of seeing this is that the 
curve parametrises compatible reductions to a Borel subgroup containing
 $H$ of the bundle $P$ over $E_a$, where the compatibility is that the 
associated reduction of $Ad(P)$
 should contain the global automorphisms of $P$ over $E_a$. There is then
a tautological reduction of the lift of $P$ to $S_u$, giving a 
Borel bundle $P_B$. This then projects to an $H$-bundle $P_H$ under the 
natural projection $B\rightarrow H$. Twisting $P_H$ by a divisor 
associated to the ramification points of $S\rightarrow\Sigma$, as for the 
Hitchin system, gives an element of the Prym variety of $S$. 

Conversely, from the element of the Prym variety, one can get back
$P_H$, and then rebuild the bundle $P$ over $\sigma(\Sigma)$, as for 
the Hitchin systems. The   curve $S_u$ then tells us how to 
reconstruct the bundle over all of $Z$.

%*******************************************************
% bibliography
%*******************************************************

\bigskip
\noindent J. C. Hurtubise: Department of Mathematics and Statistics, 
McGill University. 
(hurtubis@math.mcgill.ca)

\medskip
\noindent E. Markman: Department of Mathematics and Statistics, 
University of Massachusetts, Amherst,
(markman@math.umass.edu)

\end{document}